\documentstyle[12pt]{book}
\title{Geodesic completeness for meromorphic metrics; the case of coercive ones}
\author{Claudio Meneghini}
\begin{document} 
\bibliographystyle{plain} 
\maketitle

\def\quan{\vrule height6pt width6pt depth0pt}
\def\QUAN{\nobreak $\ $\quan\hfill\vskip0.6truecm\par}
\def\BETA{\mathop{\beta}\limits}
\def\GAMMA{\mathop{\gamma}\limits}
\def\VI{\mathop{v}\limits}
\def\UI{\mathop{u}\limits}
\def\VII{\mathop{V}\limits}
\def\WI{\mathop{w}\limits}
\def\ZETA{\mathop{Z}\limits}

\newtheorem{definition}{Definition}[chapter]
\newtheorem{lemma}[definition]{Lemma}
\newtheorem{proposition}[definition]{Proposition}
\newtheorem{theorem}[definition]{Theorem}        
\newtheorem{corollary}[definition]{Corollary}  
\newtheorem{remark}[definition]{Remark}  

\font\sdopp=msbm10
\def\ESSE {\sdopp {\hbox{S}}}
\def\ERRE {\sdopp {\hbox{R}}}
\def\CI {\sdopp {\hbox{C}}}
\def\ENNE{\sdopp {\hbox{N}}}
\def\ZETA{\sdopp {\hbox{Z}}}
\def\PI {\sdopp {\hbox{P}}}
\def\M{\hbox{\boldmath{}$M$\unboldmath}} 
\def\N{\hbox{\boldmath{}$N$\unboldmath}} 
\def\P{\hbox{\boldmath{}$P$\unboldmath}} 
\def\Y{\hbox{\boldmath{}$Y$\unboldmath}} 
\def\tr{\hbox{\boldmath{}$tr$\unboldmath}} 
\def\f{\hbox{\boldmath{}$f$\unboldmath}} 
\def\u{\hbox{\boldmath{}$u$\unboldmath}} 
\def\v{\hbox{\boldmath{}$v$\unboldmath}} 
\def\U{\hbox{\boldmath{}$U$\unboldmath}} 
\def\V{\hbox{\boldmath{}$V$\unboldmath}} 
\def\W{\hbox{\boldmath{}$W$\unboldmath}} 
\def\id{\hbox{\boldmath{}$id$\unboldmath}} 
\def\alph{\hbox{\boldmath{}$\aleph$\unboldmath}} 
\def\bet{\hbox{\boldmath{}$\beta$\unboldmath}} 
\def\gam{\hbox{\boldmath{}$\gamma$\unboldmath}} 
\def\U{\mathop{u}\limits}
\def\f{\hbox{\boldmath{}$f$\unboldmath}} 
\def\g{\hbox{\boldmath{}$g$\unboldmath}} 
\def\h{\hbox{\boldmath{}$h$\unboldmath}} 

\def\TTT{\sl}
\sloppy
{
\TTT
\section{Foreword}
This thesis is concerned with extending the idea of geodesic complete\-ness from Riemannian (or pseudo-Riemannian) context to the frame\-work of complex geometry: we take, however a completely holomorphic point of view; that is to say, a 'metric' will be a (meromorphic) symmetric section of the twice covariant holomorphic tensor bundle. Moreover, the use of smooth objects will be systematically avoided.

We shall hint at the need of reformulating some aspects of the theory of differential equations in the complex domain, originating in the interpenetration betwixt differential and algebraic aspects when analytical continuation is pursued.

As one could expect, the notion itself of path should be reformulated: a section will be devoted to bring forward a deeper discussion of this point of view, whilst we now confine ourselves to remarking that geodesics, and, more generally, paths, will be defined on Riemann surfaces which are domains over regions in the complex plane. In other words we adopt the point of view according to which a curve is an analytical continuation of an (initial) germ.

Of course,even in this broader context, geodesics will be eventually defined to be auto-parallel paths, but we are urged to focus our attention on the fact that the Levi-Civita connection will be meromorphic if the metric from which it is induced is allowed to have poles or even simply to lower somewhere in its rank.

We shall study rather more deeply a class of manifolds, namely warped products of Riemann surfaces; some hypotheses concerning their metrics will be done, but we shall show that the range of applicability of the yielded completeness theorems will not be exceedingly restricted.
}
\vfill\eject
\tableofcontents
\chapter{Analytical continuation}
{
\TTT
Poles,
essential and logarithmic singularities
}
\vfill\eject
\section{Introduction and main definitions}
We begin this section by fathoming analytical continuation techniques when applied to (initially holomorphic) functions defined on 'regions' (i.e. connected open sets) in the extended complex plane $\PI^1$ and taking values in arbitrary complex manifolds.

The definition of holomorphy in this framework is well known, and quite natural; the same holds for the treatment of the 'multi-valuedness' of analytical continuations, which leads us to functions defined on Riemann surfaces over open sets in $\PI^1$.

The situation is very different if we think of meromorphic functions, indeed there is no canonical way to 'divide' points
in general complex manifolds; analogous warning should be kept into account if we are concerned with generalizing essential singularities, thought of as points 'at' which analytic functions admit Laurent expansions with infinite order characteristic part.

Therefore, we have made the choice of focusing our attention on the following two facts:
\begin{itemize}
\item[(1)]
approaching a pole, a (complex-valued) analytic function 'tends to infinity', that is to say, every complex number is eventually left out of the images of decreasing neighbourhoods of the pole itself; conversely, holomorphic functions in a punctured disc $D\setminus\{p\}$ sharing the above property, must have, by Picard's theorem, a pole at $p$.
\item[(2)] 
Approaching an essential singularity, a (complex-valued) analytic function comes arbitrarily near any complex number (note that this is a weaker statement than Picard's theorem); the converse is also true, again by Picard's theorem.
\end{itemize}
We manage to adapt all that to functions taking values in higher dimensional manifolds: the only new feature we introduce is allowing some 'components' to be in fact holomorphic near a singularity in the domain of definition.

\begin{definition}\rm
Let  $\M$ be a connected complex manifold, $S$ a Riemann surface, $\Sigma$ a discrete subset of $S$, $F\in{\cal
O}(S\setminus\Sigma, \M)$ and $p\in\Sigma$; moreover let $\{V_k\}_{k\geq K}$ be a sequence of decreasing open
neighourhoods of $p$ in $S$, making up a fundamental system of neighbourhoods of $p$ itself, such that
$\bigcap_{k\geq
K}V_k=\{p\}$.
We shall consider  the sequences $\{V_k\}_{k\geq K}$ and $\{V_k\}_{k\geq L}$ equivalent if $K\geq L$.

Then $p$ is 
\begin{itemize}
\item
a {\bf removable singularity} for $F$ if $F$ itself is analitycally continuable up to $p$;
\item
a {\bf pole} for $F$ if there exist:
\begin{itemize}
\item {\TTT an open set } $\Omega\subset\M$;
\item {\TTT complex submanifolds }$\N\subset\Omega$ {\TTT and }$\P\subset\Omega$
\end{itemize}
such that
\begin{itemize}
\item
$dim(\P)\geq 1$;
\item
$\Omega$ {\TTT and }$\N\times\P$ {\TTT are biholomorphic};
\item
 {\TTT for every } $k$, $F\left(V_k\setminus\{p\}   \right)\subset\Omega$;
\item
$pr_1\circ F:V_k\longrightarrow\N$ 
{\TTT has a removable singularity at } $p$;
\item
$\bigcap_{k=\geq K}\overline{pr_2\circ F(V_k)}=
\emptyset$.
\end{itemize}
Here $pr_1$ and $pr_2$ are the canonical projections of the cartesian product $\N\times\P$ on its factors.
\item
an {\bf essential singularity} for $F$ if there exists a n-dimensional complex submanifold ($1\leq n \leq m$) $\N\subset\M$ such that $\bigcap_{k=1}^{\infty}\overline
{F\left(V_k\setminus\{p\}\right)}=\N$.
\end{itemize}
\label{singolarita'}
\end{definition}
\begin{definition}\rm
A function $F:S\longrightarrow\M$ is {\bf meromorphic} provided that its only singularities 
in $S$ are isolated
poles.
\label{merom}
\end{definition}
If $\M=\CI^m$, then definition \ref{singolarita'} is in fact a classification of isolated singularities: by Riemann's uniformization theorem, the universal covering  $\pi:\Sigma\longrightarrow S$ of $S$ is such that $\Sigma$ is a region of the extended complex plane: let $\{q_l\}_{l\in L}$ be the set of all preimages of $p$, where $L$ is a suitable index set;
for each $l\in L$, let $D_l$ be a neighbourhood of $q_l$ such that $\pi\vert_{D_l}$ is biholomorphic.

Then, for each $l\in L$, the $m$ component functions
$$
\Phi^{kl}=pr_k\circ F\circ\left(\pi\vert_{D_l}
\right)^{-1},\qquad k=1...m,
$$
are holomorphic on the punctured disc $D_l\setminus\{q_l\}$, hence, if $p$ is not a removable singularity, each $\Phi^{kl}$ admits a Laurent development around $q_l$.

Depending now on there being some or no infinite-order principal parts amongst all these developments, then $p$ will be an essential singularity or a pole.

If $\M$ is a more general complex manifolds, we can only affirm that the three types of singularities enumerated in definition \ref{singolarita'} are mutually exclusive.
\vskip0.5truecm

If $F$ has singularities only like in definition 
\ref{singolarita'}, we shall say that $F$ is {\bf analytic with isolated singularities} (or, shortly, 'analytic': note the difference with 'holomorphic') on $S$ and denote the space of such functions by the symbol ${\cal A}_{\Sigma}(S,\M)$.

\begin{definition}\rm
Let $\M$ be a connected complex manifold: 
\begin{itemize}
\item
a {\bf path element } in $\M$, {\bf a $\M$-valued holomorphic mapping element}, or even simply a {\bf $\M$-valued function element } is a pair $\left(U,f    \right)$, where $U$ is a connected open set in $\PI^1$ and $f$ a holomorphic function defined on $U$ and taking values in $\M$; 
\item
a {\bf $z_0$-starting path element }  in $\M$, 
{\bf a $z_0$-starting holomorphic mapping element}, or a {\bf $z_0$-starting function element } in $\M$ is a triple $\left(z_0, U,f    \right)$, where $U$ is a connected open neighbourhood of $z_0\in\CI^N$ and $f$ a holomorphic function defined on $U$ and taking values in $\M$; 
\item
a {\bf germ $\f_{z_0}$} of holomorphic mapping at $z_0$ is an equivalence class of $z_0$-starting path elements, with respect to the equivalence relation of coinciding
within a neighbourhood of $z_0$. 
\end{itemize}
\label{elem}
\end{definition}
\begin{definition}\rm
Two $\M$-valued holomorphic mapping  elements $\left(U,f    \right)$ and $\left(V,g    \right)$ {\bf are connectible} if there exists a finite sequence 
$$ \left(U_0,f_0    \right),\,\left(U_1,f_1    \right),\,....\,
\left(U_n,f_n    \right)$$
such that
\begin{itemize}
\item
$\left(U_0,f_0    \right)=\left(U,f    \right)$,\ $\left(U_n,f_n    \right)=\left(V,g    \right)$;
\item 
for every $j=0,....,n-1$,\ there holds
$$
\cases
{
U_j\bigcap U_{j+1}\not= \emptyset\cr
f_{j+1}\vert_{U_j\bigcap U_{j+1}}=f_{j}\vert_{U_j\bigcap U_{j+1}}.
}
$$
\end{itemize}
Two $\M$-valued  $z_0$-starting holomorphic mapping  elements $\left(z_0,U,f    \right)$ and $\left(z_0,V,g    \right)$ {\bf are connectible} if so are $\left(U,f    \right)$ and $\left(V,g    \right)$; two germs of holomorphic mappings $\M$ {\bf are connectible} if they admit two connectible representatives.
\label{connect_elem}
\end{definition}
\begin{remark}\rm
We shall also say that any two objects in definition \ref{connect_elem} are connectible when they spot two connectible holomorphic function elements. 
\end{remark}
\begin{definition}\rm
A {\bf standard analytical continuation} (resp. a {\bf regular analytical continuation}) of a $\M$-valued path element $\left(U,f    \right)$, or of a  $z_0$-starting path element $\left(z_0,U,f    \right)$  is a quintuple $$Q_{\M}=\left(S,\pi,j,F,\M    \right),$$ where
\begin{itemize}
\item
$S$ is a connected Riemann surface over a region of the extended complex plane $\PI^1$;
\item
$\pi\,\colon\, S\longrightarrow \CI$ is a nonconstant holomorphic mapping (resp. an everywhere maximum-rank holomorphic mapping) such that $U\subset \pi(S)$;
\item
$j\,\colon\, U\longrightarrow S$ is a holomorphic (hence open) immersion such that $\pi\circ j=id\vert_{U}$;
\item
$F\,\colon\, S\longrightarrow \M$ is a holomorphic mapping such that $F\circ j=f$.
\end{itemize}
An {\bf analytical continuation of a germ} is an analitycal continuation of any one of its representatives:
of course, this definition does not depend on the choice of the representative.
\label{analytic_cont_1}
\end{definition}

Since each standard analytical continuation fails to be regular at most at a discrete set, we could see it as a regular analytical continuation plus some points $p's$.
Due to the fact that $\pi$ is not regular at each $p$, there is no function element containing $\pi(p)$ which could be connected with the other ones; notwithstanding, $F$ may be holomorphically extended at $p$.

Every function element $\left(U,f    \right)$ admits at least the trivial analytical continuation $$\left(U,id_U,id_U,f,M    \right),$$ which may be identified with
$\left(U,f    \right)$ itself; we shall also say that the function element $\left(W,h    \right)$ is a (regular) analytical continuation of $\left(U,f    \right)$ provided that so is $\displaystyle\left(W,id_W,id_W,h,M    \right)$, or, in other words, if $\displaystyle W\supset U$ and $\displaystyle h\vert_U=f$.
\begin{definition}\rm
Let $\gamma:[0,1]\longrightarrow \PI^N$ be an embedded rectifiable analytic arc such that $\gamma(0)=z_0$: then
a {\bf regular analytical continuation} of a germ
$\f_{z_0}$ of holomorphic mapping {\bf along} $\gamma$
(resp. {\bf along} $\gamma\vert_{[0,1)}$)
is a regular analytical continuation $\displaystyle \left(S,\pi,j,F,\M \right)$
such that $\displaystyle \pi(S)\supset\gamma([0,1])$ (resp. $\displaystyle \pi(S)\supset\gamma([0,1))$ ).
\label{metonymy}
\end{definition}

\begin{definition}\rm
A {\bf morphism} between two analytical continuations
$$
Q_{\M}=\left(S,\pi,j,F,\M    \right);\quad Q_{\M}^{\prime}=\left(S^{\prime},\pi^{\prime},j^{\prime},F^{\prime},\M    \right)
$$
of the same element $\left(U,f    \right)$ is a holomorphic mapping $h\,\colon\, S\longrightarrow S^{\prime}$ such that $h\circ j=j^{\prime}$.
\label{morphism}
\end{definition}
\begin{remark}\rm
A morphism between two analytical continuations is a nonconstant (in particular open) mapping, uniquely determined on $j(U)$, hence everywhere on $S$, by $j^{\prime}\circ j^{-1}$.
Moreover, there holds $\pi^{\prime}\circ h=\pi$ and $F^{\prime}\circ h=F$ on $j(U)$, hence everywhere on $S$.

\end{remark}

The only existing morphism betwixt one analytical continuation and itself is the identity mapping; the composition of two morphisms is still a morphism; if a morphism admits a holomorphic inverse mapping, this is again a morphism: in such a case we talk about {\bf isomorphisms} of analytical continuations.    
A morphism $h$ transforms a Riemann domain $S$ into a 'larger' one if $h$ is injective, or into a less sheeted one on $U$ than $S$ otherwise.
Therefore, going forward along morphism we should find larger and larger analytical continuations, or less and less branched ones.
\begin{definition}\rm
A continuous mapping $f:X\longrightarrow Y$ ($X$ and $Y$ $T_2$ topological spaces) is {\bf proper} provided that:
\begin{itemize}
\item[(a)] it is closed;
\item[(b)] for every $y\in Y$, $p^{-1}(y)$ is compact.
\end{itemize}
\end{definition}
We recall that the above definition implies that, for every compact set $K\subset Y$, $p^{-1}(K)$ is compact. 
\begin{definition}\rm
A {\bf Riemann domain} over a region of $\CI^N$ is a pair $\left(R,p    \right)$ consisting of a complex manifold $R$ and
of a everywhere maximum-rank holomorphic mapping $p:R\longrightarrow\CI^N$;
\begin{itemize}
\item  
$R$ is {\bf proper} provided that so is $p$;
\item 
$R$ is {\bf finite} provided that
$$
\cases{
p\ \hbox{\rm is proper}\cr
\hbox{for every } Z\in\CI^N,\ p^{-1}(Z)\ \hbox{\rm is a finite set}.
}
$$
\end{itemize}
\label{rido}
\end{definition}
Note that there does not necessarily hold $p(R)=\CI^2$: we shall sometimes talk about the Riemann domain $R$, understanding the projection mapping $p$.
\begin{definition}\rm
A {\bf regular analytical continuation} of a $\M$-valued holomorphic mapping  element $\left(U,f    \right)$, or of a $\M$-valued  $z_0$-starting holomorphic mapping element $\left(z_0,U,f    \right)$ is a quintuple $$Q_{\M}=\left(S,\pi,j,F,\M    \right),$$ where
\begin{itemize}
\item
$S$ is a connected Riemann domain over a region of $\CI^N$;
\item
$\pi\,\colon\, S\longrightarrow \CI$ is an everywhere maximum-rank holomorphic mapping such that $U\subset \pi(S)$;
\item
$j\,\colon\, U\longrightarrow S$ is a holomorphic immersion such that $\pi\circ j=id\vert_{U}$; this implies that $j$ is an open mapping;
\item
$F\,\colon\, S\longrightarrow \M$ is a holomorphic mapping such that $F\circ j=f$.
\end{itemize}
An {\bf analytical continuation of a germ} is an analitycal continuation of any one of its representatives.
\label{analytic_cont_N}
\end{definition}
\vfill\eject
\section{The regular Riemann surface of an element}
\begin{definition}\rm
A regular analytical continuation 
$$
\widehat Q_{\M}=\left(\widehat S,\widehat\pi,\widehat j,\widehat F,\M    \right)
$$
is the {\bf regular Riemann surface}, or the {\bf  maximal regular analytical continuation}  of the element $\left(U,f    \right)$ if for every regular ana\-lytical continuation $Q_{\M}=\left(S,\pi,j,F,\M    \right)$ of $\left(U,f    \right)$ there exists a morphism $h\,\colon\, S\longrightarrow \widehat S$.
\label{reg_max_sup}
\end{definition}
Two regular Riemann surfaces of the same function element must be isomorphic, since they admit morphism one into each other, hence the maximal regular  analytical continuation of a function element is unique up to isomorphisms.
 
\begin{theorem}\TTT
Every path element $\left(U,f    \right)$ admits a regular Riemann surface.
\label{ex_max_sup}
\end{theorem}
{\bf Proof:} let 
\begin{itemize}
\item
${\cal U}=\{\left(U_i,f_i    \right)\}_{i\in I}$ be the class of all mapping elements connectible with $\left(U,f    \right)$;
\item
$S_0=\coprod_{i\in I}U_i$ as a topological space; 
\item
$\pi_0\,\colon\, S_0\longrightarrow \PI^1$ defined by setting $\pi_0=\coprod_{i\in I}id\vert_{U_i}$;
\item
$j_0\,\colon\, U\longrightarrow S_0$ the natural immersion;
\item
$F_0=\coprod_{i\in I}f_i$.
\end{itemize}

Let's introduce an equivalence relation in $S_0$: $z_1\in U_{i_1}$ and $z_2\in U_{i_2}$ are told to be equivalent if
and only if $\pi_0(z_1)=\pi_0(z_2)$ and $f_{i_1}=f_{i_2}$ in a neighbourhood of $\pi_0(z_1)=\pi_0(z_2)$ in $U_{i_1}\bigcap U_{i_2}$.

Name $\widehat S$ the quotient space and $q\,\colon\,S_0\longrightarrow \widehat S$ the quotient mapping: a basis for the topology of $\widehat S$ is yielded by the family $[U_i]=\{q\left(u_i    \right)\}$. Now define 
\begin{itemize}
\item
$\widehat j\,\colon\,U\longrightarrow\widehat S$ by setting $\widehat j=q\circ j_0$;
\item
$\widehat\pi\,\colon\,\widehat S\longrightarrow \PI^1$ by setting $\widehat\pi \left(q(z)    \right)=\pi_0(z)$;
\item
$\widehat F\,\colon\, \widehat S\longrightarrow \M$ by setting 
$\widehat F\left(z_i    \right)=f_i\left(z_i    \right)$.
\end{itemize}

The above mappings are well defined and continuous; the topological space $\widehat S$ is Hausdorff: indeed if $q\left(z_i    \right)\not=q\left(z_j    \right)$ and $\pi_0\left(z_i    \right)=\pi_0\left(z_j    \right)$, consider a connected neighbourhood $V$ of $\pi_0\left(z_i    \right)=\pi_0\left(z_j    \right)$, such that $f_i$ and $f_j$ are defined and different on $V$.

Let $V_i$ and $V_j$ be the disjoint copies of $V$ in $U_i$ and $U_j$ in $S_0$: we claim that $\displaystyle{q\left(V_i    \right)\bigcap q\left(V_j    \right)=\emptyset}$.
Indeed, if there were two points $w_i\in V_i$ and $w_j\in V_j$ with $\displaystyle{q\left(w_i    \right)=q\left(w_j    \right)}$ it would be $\displaystyle{f_i=f_j}$ in a neighbourhooud of $\displaystyle{\pi_0\left(w_i    \right)=\pi_0\left(w_j    \right)}$, hence in $V$, which is a contradiction.

The space $\widehat S$ is connected, since for every pair of points $p_1\in [U^{\prime}]$ and $p_2\in
[U^{\prime\prime}]$, there exists a chain $\displaystyle{{\cal K}=\{U_{i_0},U_{i_1} ..... U_{i_n}\}}$ of nonempty connected open
sets such that:
\begin{itemize}
\item
for every $k=0,....,n-1$, $U_{i_k}\bigcap U_{i_{k+1}}\not=\emptyset$;
\item
${\cal K}$ connects
$U_{i_0}=U^{\prime}$ with $U_{i_n}=U^{\prime\prime}$.
\end{itemize}

Hence the open set 
$\displaystyle{[U_{i_0}]\bigcup\cdots \bigcup[U_{i_n}]}$ is connected and contains $p_1$ and $p_2$.
Since $q$ is a local homeomorphism between $U_i$ and $q\left(U_i    \right)$, the space $\widehat S$ is a connected topological surface; since also $\widehat\pi\,\colon\,\widehat S\longrightarrow \PI^1$ is a local homemorphism, then by Poincar\'e-Volterra's theorem, $\widehat S$ is second countable too.

Now consider the atlas $$\displaystyle{\left\{\left([U_i],\widehat\pi\vert_{[U_i]}    \right)    \right\}_{i\in I}}:$$ it defines a complex structure on $\widehat S$, since for every couple of overlapping charts $\displaystyle\left([U_i],[U_j]\right)$, the transition mapping $\displaystyle{\widehat\pi\vert_j\circ\widehat\pi\vert_i^{-1}}$ is the identity mapping on an open subset in $\displaystyle{U_i\bigcap U_j}$.

By construction, the mappings $\widehat\pi, \widehat j, \widehat F$ are holomorphic with respect to this structure, hence $\displaystyle{
\widehat Q_{\M}=
\left( \widehat S, \widehat\pi, \widehat j, \widehat F,\M    \right)}$ is a regular analytical continuation of $\left(U,f    \right)$.

Now we turn to show that $\widehat Q_{\M}$ is maximal, or which is the same, it is the regular Riemann surface of $\left(U,f    \right)$.
Let $\displaystyle{Q_{\M}=\left(S,\pi,j,F,\M
\right)}$ be any regular analytical continuation of $\left(U,f    \right)$, $V$ any open set in $S$ such that $\pi$ is an immersion of $V$ in $\pi\left(V    \right)$: then the pair $\displaystyle{\left(\pi(V),F\circ\pi\vert_V^{-1}    \right)}$ is a holomorphic mapping element which may be connected with $\left(U,f    \right)$; define 
$\displaystyle{h_V: V\longrightarrow\widehat S}$ by setting $\displaystyle{h_V=q\circ\pi\vert_V}$.

If $V^{\prime}$ is another open set sharing the properties of $V$, there holds $h_V=h_{V^{\prime}}$ on $V\bigcap V^{\prime}$, hence the definition of $h$ on $V$ may be extended to $S$, since an open covering of $S$ by open sets like $V$ may be found; finally, $h$ is holomorphic and $\displaystyle{h\circ j=\widehat j}$.
\QUAN
\begin{corollary}\rm
Every germ of holomorphic mapping admits a maximal regular analytical continuation.
\end{corollary}
\vfill\eject

\section{Other Riemann surfaces of a path element}
\subsection{Superstructural singularities}
Let now 
$$\breve Q_{\M}=\left(\breve S,\breve\pi,\breve j,\breve F,\M    \right)$$
be a one dimensional regular maximal analytical continuation of the element $(U,f)$.
\begin{definition}
\rm 
A {\bf superstructural singularity} $p$ of $\breve Q_{\M}$ is a decreasing sequence of open sets $\{V_k\}_{k\geq K}\subset \breve S$
such that there exist a positive integer $n$ and a point $\displaystyle z_0\in\PI^1$, $\left(z_0=z_0(p)\right)$ depending on $p$, such that:
\begin{itemize}
\item (SPS1)
for every $k\geq K$ $V_k$ is a connected component of 
$${\breve \pi^{-
1}\left(D\left(z_0,\frac{1}{k}\right)\setminus \{z_0\}\right)};$$
\item (SPS2) for every $k\geq K$
$$\displaystyle{\breve\pi\vert_{V_k}:V_k\longrightarrow
\left(D\left(z_0,\frac{1}{k}\right)
\setminus \{z_0\}\right)}$$ is a n-sheeted covering;
\item (SPS3)
$\displaystyle{\bigcap_{k\geq K}\overline{V_k}=
\emptyset}$
\end{itemize}
\label{sps}
\end{definition}
We consider the sequences $\{V_k\}_{k\geq K}$ and $\{V_k\}_{k\geq H}$ equivalent if $H\geq K$ and say that they spot the same superstructural singularity in $$\breve Q_{\M}=\left( \breve S, \breve\pi, \breve j, \breve F,\M    \right).$$
\begin{definition}\rm
The {\bf superstrucural completion} of $\breve S$ is
$$
\widetilde S=\breve S\bigcup \left\{\hbox{superstructural singularities of }    \left( \breve S, \breve\pi, \breve j, \breve F,\M    \right)\right\}
$$
as a set, endowed with the following topology: open sets should be the open sets in $\breve S$ and a fundamental neighbourhood system of a superstructural singularity $p_{\{V_k\}_{k\geq K}}\in\widetilde S\setminus\breve S$ should be yielded by the sets $\widetilde V_k=V_k\bigcup \{p\}$.
\label{supcomp2}
\end{definition}
Set now (see definition \ref{sps})
$$
\widetilde{\pi}(s)=
\cases{\breve\pi(s) & if $s\in \breve S$\cr
z_0(s) & if $s$ is a superstructural singularity.
}
$$
\begin{lemma}
\rm
$\widetilde S$ is a connected Riemann surface.
\label{supcomp}
\end{lemma}
{\bf Proof:} if we endow $\widetilde S$ with the above topology, $\breve S$ is dense in $\widetilde S$, hence $\widetilde S$  is connected too.
Moreover $\widetilde S$ is Hausdorff, because if $p\left(\{V_k\}    \right)$ and $q\left(\{W_k\}    \right)$ are different superstructural singularities, and $\widetilde \pi(p)\not=\widetilde \pi(q)$, then $$D\left(\widetilde \pi(p),\frac{1}{k}    \right)\bigcap D\left(\widetilde \pi(q),\frac{1}{k}    \right)=\emptyset$$ if $k$ is large enough, hence $\widetilde V_k \bigcap \widetilde W_k =\emptyset$; if, instead, $\widetilde \pi(p)=\widetilde \pi(q)=z_0$, since $V_k$ and $W_k$ are different connected components of $\breve \pi^{-
1}\left(D\left(z_0,\frac{1}{k}\setminus \{z_0\}\right)\right)$, then $V_k\bigcap
W_k$, hence $\widetilde V_k \bigcap \widetilde W_k =\emptyset$ too.

The n-sheeted covering 
$$
\breve \pi\vert_{V_k}:V_k\longrightarrow D\left(z_0,\frac{1}{k}\right)
\setminus \{z_0\}
$$
is equivalent to the covering 
$$
\ell\,\colon\, D\left(0,\sqrt[n]{\frac{1}{k}}\right)\setminus \{0\} \longrightarrow D\left(z_0,\frac{1}{k}\right)\setminus \{z_0\}
$$ 
defined by setting $\ell(z)=z^n+z_0$, in the sense that there exists a homeomorphism
$$
g\,\colon\,D\left(0,\sqrt[n]{\frac{1}{k}} \right)\setminus \{0\} \longrightarrow V_k
$$ 
such that $\breve \pi\circ g=\ell$.
Then each $\widetilde V_k$ is homeomorphic to a disc by means of the homeomorphism $$\widetilde g\,\colon\,D\left(0,\sqrt[n]{\frac{1}{k}}\right)\longrightarrow \widetilde V_k$$ defined by setting
$$
\cases{
\widetilde g(z)=g(z) & if $z\not=0$\cr
\widetilde g(0)=p,
}
$$
thus $\widetilde S$ is a topological surface.

Now let's introduce a complex structure on $\widetilde S$: if $\cal A$ is the complex structure of $\breve S$, the family
$$
\widetilde{\cal A}={\cal A}\bigcup\left\{\left(\widetilde V_{p,k},\widetilde g_{p,k}^{-1}\right)\right\}_{p\in \widetilde S\setminus\breve S,k\geq K_p}
$$
defines a complex atlas on $\widetilde S$ which makes it into a Riemann surface.
\QUAN
\begin{remark}\rm
By construction, the mapping $\widetilde \pi\,\colon\, \widetilde S\longrightarrow\CI$ is holomorphic; moreover, if $\widetilde  j: U\longrightarrow \widetilde  S$ is the canonical immersion $\breve j$ of $U$ in $S$ followed by the identity mapping $id:\breve S\longrightarrow\widetilde S$, $\widetilde j$ is holomorphic too.
\label{breve}
\end{remark}
\begin{definition}\rm
A superstructural singularity of $\breve Q_{\M}=\left(\breve S,\breve\pi,\breve j,\breve F,\M    \right)$  is
\begin{itemize}
\item (RMS)
a {\bf removable singularity} for $\breve F$ if there exists $\eta\in\M$
such that $\bigcap_{k}\overline{\breve F(V_k)}=\{\eta\}$;
\item (P)
a {\bf pole} for $\breve F$ if there exist:
\begin{itemize}
\item {\TTT an open set } $\Omega\subset\M$;
\item {\TTT complex submanifolds }$\N\subset\Omega$ {\TTT and }$\P\subset\Omega$
\end{itemize}
such that
\begin{itemize}
\item
$dim(\P)\geq 1$;
\item
$\Omega$ {\TTT and }$\N\times\P$ {\TTT are biholomorphic};
\item
 {\TTT for every } $k$, $\breve F\left(V_k\setminus\{p\}   \right)\subset\Omega$;
\item
$pr_1\circ \breve F:V_k\longrightarrow\N$ 
{\TTT has a removable singularity at } $p$;
\item
$\bigcap_{k\geq K}\overline{pr_2\circ \breve F(V_k)}=
\emptyset$.
\end{itemize}
\item (ESS) an {\bf essential singularity} for 
$\breve F$ if there exists a n-dimensional complex submanifold ($1\leq n \leq m$) $\N\subset\M$ such that $\bigcap_{k=1}^{\infty}\overline{F(V_k)}=\N$.
\end{itemize}
\label{singol2}
\end{definition}
\subsection{The standard Riemann surface
of a path element }
\begin{definition}
\rm
A standard analytical continuation
$$
\breve Q_{\M}=\left(\breve S,\breve\pi,\breve j,\breve F,\M    \right)
$$
is the {\bf standard Riemann surface}, or the {\bf  maximal standard analytical continuation}  of the element $\left(U,f    \right)$ if for every standard analytical continuation $$\displaystyle{Q_{\M}=\left(S,\pi,j,F,\M    \right)}$$ of $\left(U,f    \right)$ there exists a morphism $h\,\colon\, S\longrightarrow \breve S$.
\label{max_sup_1}
\end{definition}

\begin{theorem}
\TTT
Every path element $\left(U,f    \right)$ admits a standard  Riemann surface.
\label{ex_maxx_sup}
\end{theorem}
{\bf Proof:} set 
$$
\breve S=\widetilde S\setminus\{{\hbox{poles and essential singularities of }}\breve Q_{\M} \};
$$
then $\breve S$ is a dense open subset  in $\widetilde S$, hence it is a connected Riemann surface.
Set $\breve\pi=\widetilde\pi\vert_{\breve S}$, $\widetilde\pi$ being defined in remark \ref{breve} and $\breve j=id_{\widehat S\longrightarrow\breve S}\circ\widehat j$; $\breve\pi$ and $\breve j$ are holomorphic and $\breve\pi\circ\breve j=id\vert_{U}$.

Now, for every removable singularity $p$ there exists a point $\eta\in\M$ such that the image of a neighbourhood of $p$ is contained in a local chart $({\cal W},\Psi) $ around $\eta$, hence, by Riemann's removable singularity theorem, the holomorphic mapping $$\Psi\circ\widehat F\circ g_{p,k}\,\colon\,D\left(0,\sqrt[n]{\frac{1}{k}}   \right)\setminus \{0\} \longrightarrow \CI $$
could be holomorphically extended up to the point $0$, hence $\widehat F\,\colon\,V_k\longrightarrow\M$ could be holomorphically extended up to $p$.

Therefore we are able to define a holomorphic mapping $\displaystyle{\widetilde  F\,\colon\, \widetilde V_k\longrightarrow\M}$, hence a holomorphic mapping $\displaystyle{\breve F\,\colon\,
\longrightarrow\M}$ which extends $\displaystyle{\widehat F\,\colon\,\widehat S\longrightarrow\M}$.

We claim that $\displaystyle{\breve  Q_{\M}=\left(\breve  S,\breve \pi,\breve  j,\breve  F,\M    \right)}$ is the standard Riemann surface of the element $\left(U,f    \right)$.
To prove this statement, let $Q_{\M}=\left(S,\pi,j,F,\M    \right)$ be any continuation of $\left(U,f    \right)$; set $$S^{-}=\{p\in S\ :\ \pi \hbox{ is regular at } p\}.$$

We have that $Q_{\M}^{-}=\left(S^{-},\pi,j,F,\M    \right)$ is a continuation of $\left(U,f    \right)$, hence there exists a morphism $h\,\colon\,S^{-}\longrightarrow\widehat S$.

Pick $p_0\in S\setminus S^{-}$: we can find a local chart $\left(W,\xi    \right)$ around $p_0$ such that 
$$
\cases{
\displaystyle{
\xi:W\longrightarrow\CI_z}\cr
\displaystyle{
\xi(p_0)=0}\cr
\displaystyle{
\xi(W)=D(0,1)}\cr
\displaystyle{
\pi(p)=[\xi(p)]^n\ n\geq 2\ \hbox{for every }p\in W}.
}
$$

$\bullet$ {\bf At first} let's suppose that for every pair of points $p^{\prime},p^{\prime\prime}$ in $W$ we have:
$$
\pi(p^{\prime})=\pi(p^{\prime\prime})\hbox{ and }
F(p^{\prime})=F(p^{\prime\prime}) \Longrightarrow p^{\prime}=p^{\prime\prime};
$$
in this case the mapping $h\,\colon\,W\setminus \{p_0\}\longrightarrow\widehat S $ is injective, because for every $p^{\prime},p^{\prime\prime}$ in $W$ with $\pi(p^{\prime})=\pi(p^{\prime\prime})$ and two neighbourhoods $W^{\prime},W^{\prime\prime}$ respectively of $\pi(p^{\prime})\hbox{ and }\pi(p^{\prime\prime})$ such that $\pi(W^{\prime})=\pi(W^{\prime\prime})$ and $\pi\vert_{W^{\prime}},\pi\vert_{W^{\prime\prime}}$ are immersions on $A=\pi(W^{\prime})=\pi(W^{\prime\prime})$ the two path elements $\left(A,F\circ\pi\vert_{W^{\prime}}^{-1}    \right),\left(A,F\circ\pi\vert_{W^{\prime\prime}}^{-1}    \right)$ are different and without equivalent points in $S_0$ (see the proof of theorem \ref{ex_max_sup}); therefore $h\vert_{W}$ is an open immersion of $W$ in $\widehat S$.

Set $W^{*}_{k}=\xi\left(D\left(0,\sqrt[n]{\frac{1}{k}}    \right)\setminus \{0\}    \right)$ and $V_k=h\left(W^{*}_{k}   \right)$ for every $k\geq 1$, the sequence $\{V_k\}_{k\geq 1}$ is a removable singularity in $\widehat S$.

Indeed,
\begin{itemize}
\item (SPS1):
there couldn't exist a connected open set contained in
$$\left(\widehat \pi    \right)^{-1}\left(D\left(0,\frac{1}{k}    \right)\setminus \{0\}    \right)$$ and
properly containing $V_k$, since every chain $Z_1,....,Z_u$ of open sets in 
$\left(\widehat \pi    \right)^{-1}\left(D\left(0,\frac{1}{k}    \right)\setminus \{0\}    \right)$ 
such that $\pi\vert_{V_j}\ j=0,....,u$ is an immersion, with $Z_1\subset V_k$ and 
$Z_j\bigcap Z_j\bigcap Z_{j+1}\not=\emptyset\ j=0,....,u-1$, is contained in $V_k$;
\item (SPS2):
this property holds for $\left(V_k,\widehat \pi\vert_{V_k}    \right)$ in the same way as for $\left(W^{*}_k,\pi\vert_{W^{*}_k}    \right)$;
\item (SPS3):
if there existed $\widehat p_{0}$ in $\bigcap_{k\geq 1}\overline{V_k}$, then $\widehat \pi$ would be regular at $\widehat p_0$ and $\pi$ at $p_0$, which is a contradiction.
\item (RMS) in definition \ref{singol2}:
the behaviour of $\widehat F$ on $V_k$ is the same as that of $F$ on $W^{*}_k$.
\end{itemize}

Let $\widehat p_0\in \breve S\setminus \widehat S$ be the removable singularity spotted by the sequence $\{V_k\}_{k\geq 1}$; we get a holomorphic extension of $h$ up to $p_0$ by setting $h(p_0)=\widehat p_0$.
\vskip1truecm
$\bullet\bullet$ {\bf Instead}, if there exist some pairs $\left(p^{\prime},p^{\prime\prime}    \right)$ of distinct points in $W$ with $\pi(p^{\prime})=\pi(p^{\prime\prime} )$ and 
$F(p^{\prime})=F(p^{\prime\prime} )$ we can find again a local chart $\left(W,\xi    \right)$ around $p_0$ such that $\xi(p_0)=0$, $\xi:W\longrightarrow\CI_z$, $\xi(W)=D(0,1)$ and $\pi(p)=[\xi(p)]^n\ n\geq 2$ for every $p\in W$.

Let $y=F(p_0)$ and $({\cal Z},(\Theta^1...\Theta^N))$ a local chart around $y$; maybe shrinking $W$, we may suppose that $F(W)\subset{\cal Z}$; set $E=F\xi^{-1}$ and $z^mG_{\nu}(z)=\Theta^{\nu}\circ E(z)$, with $m\in\ENNE$
and $G_{\nu}(z)\not=0,\ \nu=1\cdots N$, for every $z\in D\left(0,1    \right)$.
Set $m=l.c.m.\left(m_1\cdots m_N    \right)$, $r=G.C.D.\left(n,m \right)$ and $u={\bf e}^{\frac{2\pi {\bf i}}{r}}$.

Let $\breve U$ be a (small enough) neighbourhood of $0$,$$
s=min \left\{k>0\,;\, k\vert r,\ G_{\nu}\left(u^k z    \right)=G_{\nu}(z),\ \nu=1\cdots N, \hbox{for every}\ z\in\breve U\right\} 
$$
and $t=r/s$; we claim that, maybe shrinking $\breve U $, for every pair $\left(z^{\prime}    \right),\left(z^{\prime\prime}    \right)$ in $\breve U$
$$
\left(\left(z^{\prime}    \right)^{n}=\left(z^{\prime\prime}    \right)^{n}\hbox{ and } E(z^{\prime})=E(z^{\prime\prime})    \right)\iff 
\left(\left(z^{\prime}    \right)^{t}=\left(z^{\prime\prime}    \right)^{t}\right).
$$
Indeed, if $\left(z^{\prime}    \right)^{n}=\left(z^{\prime\prime}    \right)^{n}\hbox{ and } E(z^{\prime})=E(z^{\prime\prime})  $, then 
$$
\left[G_{\nu}(z^{\prime})    \right]^n=
\frac{\left[\Psi_{\nu}\circ E(z^{\prime})    \right]^n}{\left(z^{\prime}\right)^{nm_{\nu}}}=
\frac{\left[\Psi\circ E(z^{\prime\prime})    \right]^n}{\left(z^{\prime\prime}\right)^{nm_{\nu}}}=
\left[G_{\nu}(z^{\prime\prime})    \right]^n,
$$
hence, if $\breve U$ is small enough, we have $G_{\nu}(z^{\prime}) =G_{\nu}(z^{\prime\prime}),\ {\nu}=1\cdots N $: this implies $\left(\left(z^{\prime}    \right)^{m_{\nu}}=\left(z^{\prime\prime}    \right)^{m_{\nu}}\right),\ {\nu}=1\cdots N $, thus
$\left(\left(z^{\prime}    \right)^{m}=\left(z^{\prime\prime}    \right)^{m}\right)$.

Set $\vartheta^{\prime}=\arg\left(z^{\prime}    \right)$ and 
$\vartheta^{\prime\prime}=\arg\left(z^{\prime\prime}    \right)$, then, for suitable $a,b\in\ZETA$, 
$$
\cases{
\vert z^{\prime}\vert= \vert z^{\prime\prime}\vert\cr
n\left(\vartheta^{\prime}- \vartheta^{\prime\prime}   \right)=2a\pi\cr
m\left(\vartheta^{\prime}- \vartheta^{\prime\prime}   \right)=2b\pi,
}
$$
hence $r\left(\vartheta^{\prime}- \vartheta^{\prime\prime}   \right)=2e\pi$, , with suitable $e\in\ZETA$ so, $(z^{\prime})^r =(z^{\prime\prime})^r$, that is to say $z^{\prime\prime}= {\bf e}^{\frac{2\pi k{\bf i}}{r}} z^{\prime}$, with $1\leq k \leq r$.

Perhaps shrinking again $\breve U$, the set
$$
I=\left\{k\in\ZETA\,:\, \forall z\in\breve U, {\nu}=1\cdots N,\  G_{\nu}\left(u^kz    \right)=G_{\nu}(z)    \right\}
$$
is a subgroup in $\ZETA$ containing $r$, generated by $s$, hence $r=st$.

Let's shrink $\breve U$ again, so that, for every $k,\ 1\leq k\leq r$ and $r$ not multiple of $s$, $G\left(u^kz^{\prime}    \right)=G(z^{\prime})\Longrightarrow z^{\prime}=0 $: then, for every $z^{\prime},z^{\prime\prime}\in \breve U\setminus \{0\}$ such that $z^{\prime}\not=z^{\prime\prime}$, if $G(z^{\prime}) =G(z^{\prime\prime}) $ and $\left(z^{\prime\prime}    \right)= {\bf e}^{\frac{2\pi k{\bf i}}{r}} z^{\prime}$ then $k=ps$ and
$$
\left(z^{\prime\prime}    \right)^t=\left({\bf e}^{\frac{2\pi k{\bf i}}{r}}\right)^{stp}\left( z^{\prime}    \right)^t=\left( z^{\prime}    \right)^t
.
$$
Vice versa, if $\left(z^{\prime\prime}    \right)^t=\left( z^{\prime}    \right)^t$, then $\left(z^{\prime\prime}    \right)=\left({\bf e}^{\frac{2\pi k{\bf i}}{r}}\right)^{ps}\left( z^{\prime}    \right)$, so that
\begin{eqnarray*}
G_{\nu}\left(z^{\prime}     \right)&=&
G_{\nu}\left(\left({\bf e}^{\frac{2\pi k{\bf i}}{r}}\right)^{s}z^{\prime}     \right)\\
&=&
G_{\nu}\left(\left({\bf e}^{\frac{2\pi k{\bf i}}{r}}\right)^{2s}z^{\prime}     \right)\\
&=&
\cdots=\\
&=& G_{\nu}\left(\left({\bf e}^{\frac{2\pi k{\bf i}}{r}}\right)^{ps}z^{\prime}     \right)\\
&=&
G_{\nu}\left(z^{\prime\prime}     \right)
\ ({\nu}=1\cdots N).
\end{eqnarray*}
Thus $\left(z^{\prime}     \right)^n=\left(z^{\prime\prime}     \right)^n$, hence $\left(z^{\prime}     \right)^m=\left(z^{\prime\prime}     \right)^m$ and $E\left(z^{\prime}     \right)=E\left(z^{\prime\prime}     \right)$.

We may suppose $\breve U=D\left(0,\varepsilon    \right)$, for a suitable real number $\varepsilon$.
Consider $\ell:\CI_z\longrightarrow\CI_w$, defined by setting $w=z^t$ and set $\breve U^{\sharp}=\ell(\breve U)$: then the functions
$$
\cases{\displaystyle{
\pi^{\sharp}:\breve U^{\sharp}\longrightarrow\CI}\cr
\displaystyle{
\pi^{\sharp}(w)=w^{n/t}}
}
$$
and 
$$
\cases{\displaystyle{
F^{\sharp}:\breve U^{\sharp}\longrightarrow\M}\cr
\displaystyle{
F^{\sharp}(w)=F\circ\xi^{-1}(z)}}
$$
are well defined and holomorphic throughout $\breve U^{\sharp}$.

Moreover, $\pi\circ\xi^{-1}=\pi^{\sharp}\circ\ell$ and $F\circ\xi^{-1}=F^{\sharp}\circ\ell$; now let $V\subset W$ be an open set such that $\pi\vert_{V}$ is biholomorphic and $j=(\pi\vert_{V})^{-1}$: then $\displaystyle{\left(W,\pi,j,F,\M\right)}$ and 
$\displaystyle{\left(\breve U^{\sharp},\pi^{\sharp},\ell\circ\xi\circ j,
 F^{\sharp}\right)}$ are analytical continuations of the same path element 
$\displaystyle{\left(\pi(V),F\circ j\right)}$, hence their standard Riemann surfaces are the same up to isomorphism: we shall denote both them by $\breve S$.

Let $h^{\prime}\,:\,\left(\breve U^{\prime}\setminus \{0\}\longrightarrow \widehat S    \right)$ be the natural morphism: since $h^{\prime }\circ\ell= h$, it may be holomorphically extended up to $0$, hence $h$ may be holomorphically extended up to $p_0$.
\QUAN
\begin{corollary}
\rm
Every germ of path admits a standard Riemann surface.\QUAN
\end{corollary}
\subsection{The extended Riemann surface of a path element}
\begin{definition}
\rm Let now 
$$\widehat Q_{\M}=\left(\widehat S,\widehat\pi,\widehat j,\widehat F,\M    \right)\left(\hbox{\rm\ resp. } \breve Q_{\M}=\left(\breve S,\breve\pi,\breve j,\breve F,\M    \right)\right)
$$
be the regular, (resp. standard), maximal analytical continuation of the path element $(U,f)$: the {\bf extended Riemann surface} (or the {\bf extended maximal analytical continuation}) of $(U,f)$ is the quintuple $$\widetilde Q_{\M}=\left(\widetilde S,\widetilde \pi,\widetilde j,\widetilde F,\M    \right),$$ where
\begin{itemize}
\item $\widetilde S$ is the superstructural completion of $\widehat S$, or, which is the same, of $\breve S$ (see definition \ref{supcomp2});
\item $\widetilde\pi:\widetilde S\longrightarrow\PI^1$ is the unique holomorphic extension of $\widehat\pi$ (or $\breve\pi$) to $\widetilde S$;
\item $\widetilde j=id_{\widehat S\longrightarrow\widetilde S}\circ\widehat j$;
\item $\widetilde F=\breve F:\breve S\longrightarrow\M$.
\end{itemize}
\label{ext_rs}
\end{definition}
\subsection{Other features}
\begin{lemma}\rm
Let $\M$ and $\N$ be complex manifolds, $G\in{\cal O}\left(\M,\N    \right)$ and $\left(U,f    \right)$ a path element in $\M$; if $\left(\widetilde R,\pi_{\widetilde R},j_{\widetilde R},F_{\widetilde R},\M    \right)$ is the standard Riemann surface of $\left(U,f    \right)$ in $\M$ and $\left(\widetilde S,\pi_{\widetilde S},j_{\widetilde S},F_{\widetilde S},\N   \right)$ is the standard Riemann surface of $\left(U,G\circ f    \right)$ in $\N$, then there exists a holomorphic function $\displaystyle{h \,:\, \widetilde R\longrightarrow\widetilde S}$ such that $\displaystyle{\pi_{\widetilde S}\circ \psi=\pi_{\widetilde R}}$.
\label{passage}
\end{lemma}
{\bf Proof:} we should merely note that $\left(\widetilde R,\pi_{\widetilde R},j_{\widetilde R},G\circ F_{\widetilde R},\N    \right)$ is a continuation of $\left(U,G\circ f    \right)$ in $\N$, hence the existence of $h$ is postulated by definition \ref{max_sup_1}.
\QUAN
\vfill\eject

\section{Logarithmic singularities}
\subsection{Generalities}

This section deals with 'logarithmic singularities' in analytical continuations, i.e; points resembling $0$ in connection with $z\mapsto \log z$: it will turn out that a complex structure at such 'points' couldn't be made up;
notwithstanding, they may be bept into account by introducing a weaker structure.
\begin{definition}
\rm
{\bf A logarithmic singularity} $q$ of a regular maximal analytical continuation 
$\widehat Q_{\M}=\left(\widehat S,\widehat\pi,\widehat j,\widehat F,\M    \right)$ of some path element is a sequence of decreasing open sets $\{V_k\}_{K\geq K}$ of $\widehat S$ such that there exist a point $z_0\in\PI^1$, depending on $q$ such that:
\begin{itemize}
\item (LS1) 
for every $k\geq K$ $V_k$ is a connected component of $$\displaystyle{\widehat \pi^{-
1}\left(D\left(z_0,\frac{1}{k}\setminus \{z_0\}\right)\right)};$$
\item (LS2) 
for every $k\geq K$ and every (real) nonconstant closed path $\gamma:[0,1]\longrightarrow D(z_0,1/k)\setminus\{z_0\}$, with nonzero winding number around $z_0$, every lifted path $$\beta:[0,1]\longrightarrow \widehat\pi^{-1}\left(D(z_0,1/k)\setminus\{z_0\}\right)$$ with respect to the topological covering $\pi$ is not a closed path, i.e. $\beta(0)\not=\beta(1)$;
\item (LS3) 
$\displaystyle{
\bigcap_{k\geq K}\overline{V_k}=\emptyset}$
\end{itemize}
Of course, even in this case, we consider the sequences $\{V_k\}_{k\geq K}$ and $\{V_k\}_{k\geq H}$ equivalent if $H\geq K$ and say that they spot the same logarithmic singularity in $\widehat Q_{\M}=\left( \widehat S, \widehat\pi, \widehat j, \widehat F,\M    \right)$.
\label{log_sing}
\end{definition}
\begin{lemma}\rm
Let $\widehat Q_{\M}=\left(\widehat S,\widehat\pi,\widehat j,\widehat F,\M    \right)$ be the regular Riemann surface of
some path element: then $$\{logarithmic\ singulartities\} \bigcap \{superstructural\ singularities\} =\emptyset.
$$
\label{vuotoint}
\end{lemma}
{\bf Proof:} let $b$ be a superstructural singularity: then, by SPS2 in definition \ref {sps}, there exist:
\begin{itemize}
\item
a local chart $(U,\phi)$ around $b$ with $\phi(b)=0$;
\item
a local chart $(V,\eta)$ around $0$, with $\eta(0)=0$;
\item
an integer $N$
\end{itemize}

such that $$\eta\circ\widehat\pi\circ\phi^{-1}(z)=z^N .$$
Consider now the nonconstant closed path $$\gamma:[0,1]\longrightarrow D(z_0,1/k)\setminus\{z_0\}$$ defined by setting
$\displaystyle{\gamma(t)=\phi^{-1}\left((1/2k)e^{2\pi iNt}\right)}$: its winding number around $z_0$ is $N$, but the lifted path $\beta:[0,1]\longrightarrow\widehat S $ defined by setting $\beta(t)=\eta\left(\sqrt{N}{1/2k}e^{2\pi it}\right)$ is closed, contradicting (LS2) in definition \ref{log_sing}, hence $b$ is not logarithmic.

On the other hand, suppose $q$ is a logarithmic singularity: by (LS2), $\pi\vert_{V_k}$ couldn't be a n-sheeted covering, contradicting (SPS2) in definition \ref{sps}, hence $q$ is not superstructural.
\QUAN
Consider now the set $B$ of logarithmic singularities of the regular Riemann surface $\widehat Q_{\M}=\left(\widehat S,\widehat\pi,\widehat j,\widehat F,\M    \right)$: set $S^{\sharp}=\widehat S\bigcup B$ as a set and introduce a topology on $S^{\sharp}$: open sets should be the open sets in $\widehat S$ and a fundamental neighbourhood system of
a logarithmic singularity $q_{\{V_k\}_{k\geq K}}\in S^{\sharp}\setminus\widehat S$ should be yielded by the sets $V_k^{\sharp}=V_k\bigcup \{q\}$.
\begin{lemma}\rm
$S^{\sharp}$ admits no complex structure at $q$.
\end{lemma}
{\bf Proof:} were there one, we could find charts $\left({\cal W},\phi   \right)$ around $q$ and $\left({\cal V},\psi   \right)$ around $z_0$ such that 
$$
\psi\circ\pi\circ\phi^{-1}\left(\zeta   \right)
=\zeta^N
$$
for some integer $N>0$.

This fact would imply $\pi\vert_{{\cal W}\setminus\{q\}}$ to be a n-sheeted covering of ${\cal V}\setminus\{z_0\}$, contradicting, as shown in lemma
\ref{vuotoint}, (LS2) in definition \ref{log_sing}.
\QUAN
We say that $q$ is a 'boundary' point because, in the given topology, every point 'near' $q$ does admit a complex structure.
\begin{lemma}\rm
$\widehat\pi$ admits a unique continuous extension $\pi^{\sharp}$ to $S^{\sharp}$.
\label{proest}
\end{lemma}
{\bf Proof:} let $b\in B$ and $\{V_k\}$ be the sequence spotting $b$: define
$$
\pi^{\sharp}(q)=\cases{\widehat\pi(q) & if $q\in V_k$\cr
z_0 & if $q=b,$
} 
$$
where $z_0$ is the common centre of the discs onto which the $V_k's$ are projected.
Now $\pi^{\sharp}$ is continuous at all points in $V_k$; moreover, for every neighbourhood $G$ of $z_0$, 
$$\pi^{\sharp\ -1}(G)\supset \pi^{\sharp\ -1}(z_0) \bigcup\widehat\pi^{-1}(G\setminus\{z_0\}),$$
hence, if we set  $$H=\{b\}\bigcup\widehat\pi^{-1}(G\setminus\{z_0\}), $$ we have that $H$ is a neighbourhood of $b$ in $S^{\sharp}$ such that $\pi^{\sharp}(H)\subset G$, proving continuity at $b$.

Arguing by density, we conclude that this extension is unique.
\QUAN
\begin{definition}\rm
A logarithmic singularity of $\widehat Q_{\M}=\left(\widehat S,\widehat\pi,\widehat j,\widehat F,\M    \right)$  is
\begin{itemize}
\item (RMLS)
a {\bf removable logaritmic singularity} for $\widehat F$ if there exists $\eta\in\M$
such that $\bigcap_{k}\overline{\widehat F(V_k)}=\{\eta\}$;
\item (PLS)
a {\bf polar logarithmic singularity} for $\widehat F$ if there exist:
\begin{itemize}
\item {\TTT an open set } $\Omega\subset\M$;
\item {\TTT complex submanifolds }$\N\subset\Omega$ {\TTT and }$\P\subset\Omega$
\end{itemize}
such that
\begin{itemize}
\item
$dim(\P)\geq 1$;
\item
$\Omega$ {\TTT and }$\N\times\P$ {\TTT are biholomorphic};
\item
 {\TTT for every } $k$, $\widehat F\left(V_k\setminus\{p\}   \right)\subset\Omega$;
\item
$pr_1\circ F:V_k\longrightarrow\N$ 
{\TTT has a removable singularity at } $p$;
$\bigcap_{k=\geq K}\overline{pr_2\circ \widehat F(V_k)}=
\emptyset$.
\end{itemize}
\item (ELS)
an {\bf essential logarithmic singularity} for $\widehat F$ if there exists a n-dimensional complex submanifold ($1\leq n \leq m$) $\N\subset\M$ such that $\bigcap_{k=1}^{\infty}\overline{\widehat F(V_k)}=\N$.
\end{itemize}
\label{log_singol2}
\end{definition}
\begin{lemma}\rm
\rm
For every removable logarithmic singularity $r$ of $\widehat Q_{\M}$, $\widehat F$ admits a unique continuous extension $F^{\sharp}$ to $r$.
\label{funest}
\end{lemma}
{\bf Proof:} let $\{V_k\}$ be any sequence spotting $r$: set $\bigcap_{k}\overline{\widehat F(V_k)}=\{\eta\}$ and define
$$
F^{\sharp}(q)=\cases{\widehat F(q) & if $q\in V_k$\cr
\eta & if $q=r.$
} 
$$
Now $\pi^{\sharp}$ is continuous at all points in $V_k$; moreover, for every neighbourhood $G$ of $\eta$, 
$$F^{\sharp\ -1}(G)\supset \pi^{\sharp\ -1}(\eta) \bigcup\widehat\pi^{-1}(G\setminus\{\eta\}),$$
hence, if we set  $H=\{r\}\bigcup\widehat\pi^{-1}(G\setminus\{\eta\}) $, we have that $H$ is a neighbourhood of $r$ in $S^{\sharp}$ such that $F^{\sharp}(H)\subset G$, proving continuity at $r$.

Arguing by density, we conclude that this extension is unique.
\QUAN
By lemmata \ref{vuotoint}, \ref{proest} and \ref{funest}, we are allowed to add logarithmic singularities
to the {\it extended} Riemann surface $\widetilde Q_{\M}$ of a path element, and this may be done without ambiguity.
\subsection{Riemann surfaces with boundary}
We axiomatize:
\begin{definition}
\rm
A quintuple $\displaystyle Q_{\M}^{\natural}=\left(S^{\natural},\pi^{\natural},j^{\natural},F^{\natural},\M\right)$, is 
an {\bf analytical continuation with boundary},
or {\bf with logarithmic singularities} of the function element $\displaystyle(U,f)$ if there exists an analytical continuation $\displaystyle\widetilde Q_{\M}=
\left(S,\pi,j,F,\M    \right)$ of 
$\displaystyle(U,f)$ such that
$\widetilde Q_{\M}$ and 
$Q_{\M}^{\natural}$ share the following properties:
$\displaystyle S^{\natural}$ is a Riemann surface with boundary such that $\displaystyle int\left( S^{\natural}  \right)=S$,
$\pi$ admits a  unique continuous extension 
$\pi^{\natural}: S^{\natural}
\longrightarrow\PI^1$ to $S^{\natural}$,
$j^{\natural}=id_{S\longrightarrow S^{\natural}}\circ j$ and 
$F$ admits a unique continuous 
$ F^{\natural}$ to $S^{\natural}$.
\label{loganalcont}
\end{definition}
\begin{definition}
\rm let $\widetilde Q_{\M}=\left(\widetilde S,\widetilde \pi,\widetilde  j,\widetilde F,\M    \right)$ be the  extended Riemann surface of the element $(U,f)$: the {\bf extended Riemann surface with boundary}, or the {\bf extended Riemann surface with logarithmic singularities} of $(U,f)$ is the quintuple $\widetilde Q_{\M}^{\sharp}=\left(\widetilde S^{\sharp},\widetilde \pi^{\sharp},\widetilde  j^{\sharp},\widetilde F^{\sharp},\M    \right)$,
where
\begin{itemize}
\item $\widetilde S^{\sharp}=\widetilde S\cup B_R\cup B_P\cup B_E$ is the Riemann surface with boundary associated with $\widetilde S$ and
$B_R\cup B_P\cup B_E$; moreover 
\begin{itemize}
\item
$B_R=
\{\hbox{\rm removable logarithmic singularities of }
\widehat Q_{\M}\}$;
\item
$B_P=\{\hbox{\rm polar logarithmic singularities of }
\widehat Q_{\M}\}$;
\item
$B_E=\{\hbox{\rm essential logarithmic singularities of }\widehat Q_{\M}\}
$
\end{itemize}
and 
\item $\widetilde\pi^{\sharp}:\widetilde S^{\sharp}\longrightarrow\CI$ is the unique continuous extension of $\widetilde\pi$ to $\widetilde S^{\sharp}$;
\item $\widetilde j^{\sharp}=id_{\widetilde S\longrightarrow\widetilde S^{\sharp}}\circ\widetilde j$;
\item $\widetilde F^{\sharp}$ is the unique continuous extension of $\widetilde F$ to $\widetilde S\bigcup B_R$.
\end{itemize}
\label{ext_rsb}
\end{definition}
\chapter{Ordinary differential equations}
{\TTT in the complex domain}
\vfill\eject
\section{Fundamentals}
The theory of ordinary differential equations in the complex domain is quite a classical subject: we start by reporting some statements which we shall use in the following; for further details, the reader is referred to \cite{hille} or \cite{ince}; we go on by reformulating the basic theory in a more geometrical fashion, in view of merging these results with the techniques of analytical continuation.
In the following, the bracket symbol with a subscript ($[f]_p$) will denote the germ at $p$ of the holomorphic mapping $f$ and, given a germ $\u_{v_0}$ at $v_0$, the $\u_{v_0}(v_0)$ will denote the common value at $v_0$ of all the representatives $u$ of $\u_{v_0}$.

Endow now $\CI^n$ with the maximum coordinate norm: 
$$
\vert\vert W \vert\vert=\vert\vert (w^1....w^n)\vert\vert=max(\vert w^j\vert).
$$
We state without proof the classical existence and uniqueness theorem of the theory of ordinary differential equations
in the complex domain:
\begin{theorem}
\TTT
Let the first order N-dimensional Cauchy's problem
\begin{equation}
\spadesuit\quad
\cases{
W^{\prime}=F(W(z),z)\cr
W(z_0)=W_0
}
\label{probcauchy}
\end{equation}
be given, where $(W,z)\longmapsto F(W,z)$ is a mapping of $N+1$ complex variables, defined and holomorphic in a $(N+1)$-cylinder 
$$
\Delta=\{\vert w^1-w_0^1\vert\leq b^1.....\vert w^N-w_0^N\vert\leq b^N,\vert z-z_0\vert\leq a\}
$$
and taking values in $\CI^N$.
Assume that
\begin{itemize}
\item (A) $\max_{\Delta}\vert\vert F(W,z)\vert\vert=M$;
\item (B) $\vert\vert F(U,z)-F(V,z)\vert\vert\leq K\vert\vert U-V\vert\vert$ for every $(U,z)$ and $(V,z)$ in $\Delta$.
\end{itemize}
Define a disc $$D=\{\vert z-z_0\vert \leq r\},$$ where 
$$r<min(a,b^1/M.....b^N/M,1/K):$$
then there exists a unique vector valued holomorphic mapping $$z\mapsto W(z):D\longrightarrow\CI^N$$ 
(such that $graph(W)\subset \Delta$)
satisfying 
\ref{probcauchy} in a neighbourhood of $z_0$.
\QUAN
\label{exunduniq}
\end{theorem}
\begin{remark}
\rm
If $N=1$ condition B is automatically fulfilled, since $\partial F/ \partial z$ is bounded on $\Delta$; the same holds, for any $N$, if $\spadesuit$ is autonomous, i.e. if $F$ does not depend explicitly on $z$.
\label{exundun}
\end{remark}
We now reformulate the theory of local solutions of ordinary differential equations in the complex domain in terms of germs of holomorphic mappings: we confine ourselves to first order systems in normal form, without loss of generality, since every N-th order normal differential equation is in turn equivalent to a system of N first order equations in normal form.

We start by defining the composition operation of a funcion and a germ: let 
$U_{v_0}$ be a germ of a $\CI^P$-valued holomorphic function at $v_0\in\CI$ and $S$ be a holomorphic function of $P$ complex variables, holomorphic at $U_{v_0}(v_0)$ and taking values in $\CI^Q$; let $U$ be any representative of $U_{v_0}$ in a neighbourhood of $v_0$.
\begin{definition}\rm
$S\circ U_{v_0}=[S\circ U]_{v_0}$.
\label{fungerm}
\end{definition}
The above definition is well posed, since, if $V$ is another representative of $U_{v_0}(v_0)$, then $U=V$ in a neighbourhood of $v_0$, hence  $S\circ U=S\circ V$ too.
\begin{definition}
\rm
A {\bf first order differential operator} on $({\cal O}_{v_0},\CI^N)$ is any mapping $D_S:({\cal O}_{v_0},\CI^N)\longrightarrow({\cal O}_{v_0},\CI^N)$ defined 
by setting:
$$
\heartsuit\quad D_S(U_{v_0})=U_{v_0}^{\prime}+S\circ(U_{v_0}\times\id_{v_0}),
$$ 
where 
\begin{itemize}
\item (1) 
$({\cal O}_{v_0},\CI^N)$ is the stalk of germs of holomorphic mappings at $v_0$, taking values in $\CI^N$;
\item (2)
$ U_{v_0}^{\prime}=[U^{\prime}]_{v_0}$, $U$ being any representative of $ U_{v_0}$
\item (3) $S$ is a mapping of $N+1$
 complex variables, taking values in $\CI^N$ and holomorphic at $(U_{v_0}(v_0),v_0)$.
\end{itemize}
\begin{definition}\rm
A complex-analytic (resp a real-analytic) mapping $S$ is {\bf uniformly Lipschitz-like} at 
$(U_0,v_0)$  
provided that there exists $K\in\ERRE$ such that 
$$
\vert\vert S(U,v)-S(V,v)\vert\vert\leq K\vert\vert U-V \vert\vert,
$$
uniformly on $v$ in a neighbourhooud of $(U_0,v_0)$.
\end{definition}

We shall say that $D_S$ is 
\begin{itemize}
\item
{\bf autonomous} provided that $S$ does not depend explicitely on $v$;
\item
{\bf uniformly Lipschitz-like} provided that so is $S$.
\end{itemize}
\label{diffop}
\end{definition}
\begin{definition}
\rm
\begin{itemize}
\item
A {\bf system} of N {\bf first order ordinary differential equations } at $v_0$ is an expression of the form $D(U_{v_0})=0$, where $D$ is a first order differential operator on $({\cal O}_{v_0},\CI^N)$; that system is {\bf autonomous} or {\bf uniformly  Lipschitz-like} provided that so is $D$;
\item
a {\bf solution}, sometimes a {\bf germ solution} of the above system is any element $U_{v_0}\in D^{-1}(0)$;
\item a N-dimensional, first order Cauchy's problem at $v_0$ is a pair consisting in a system of N first order ordinary
differential equations at $v_0$ and an `initial-value` assignment $U_{v_0}(v_0)=U_0\in\CI^N$.
\end{itemize}
\label{sys}
\end{definition}
We state a natural consequence of theorem \ref{exunduniq} and remark \ref{exundun}:
\begin{theorem}
\TTT
Every first order one dimensional Cauchy's problem at a given 
point admits a unique solution; every uniformly  Lipschitz-like N-dimensional Cauchy's problem admits a unique solution. In particular, the same holds about autonomous systems.
\label{gerexuq}
\end{theorem}
The question naturally arises about how far may be pushed the analytical continuation of the germ solution whose existence is asserted by theorem \ref{gerexuq}: in this section we confine ourselves to a quite general statement.
\begin{theorem}
\TTT
Let the N-dimensional Cauchy's problem 
$$
\cases{
U_{v_0}^{\prime}=S\circ (U_{v_0}\times\id_{v_0})\cr
U_{v_0}(v_0)=U_0
}
$$
be given, with $S$ holomorphic at $(U_0,v_0)$; suppose that:
\begin{itemize}
\item
the (unique) germ solution $U_{v_0}$ admits regular  analytical $\quad $ continuation 
$\phi$
along the embedded analytic arc $$\gamma:[a,b]\longrightarrow\CI;$$
\item 
the holomorphic mapping $S$ admits regular analytical continuation 
$\Omega$
along $(\phi\circ\gamma)\times\gamma$; 
\end{itemize}
then, set $v_1=\gamma(b)$, $U_{v_1}=[\phi]_{v_1}$, we have: $$
U_{v_1}^{\prime}=\Omega
\left(U_{v_1}\times\id_{v_1}\right).$$
\label{ancont}
\end{theorem}
{\bf Proof:} it should be merely noted that there exists an open neighbourhood ${\cal V}$ of $v_0$ such that for every $v\in{\cal V}$, $\phi^{\prime}(v)-\Omega(\phi(v),v)=0$: by hypothesis the above expression is analytically continuable up to $v_1$ across $\gamma$ and the result follows.
\QUAN
\begin{proposition}\rm
Let $U$ be a $\ERRE^P$-valued analytic function on $(a,t_0)\bigcup (t_0,b)$, which is continuous at $t_0$, with $U(t_0)=U_0$.
Suppose that $U^{\prime}(t)=F\left(U(t),t    \right)$ in $(a,t_0)$ and $(t_0,b)$, where $F$ is a $\ERRE^P$-valued real analytic mapping 
in $U(a,b)\times (a,b)$, uniformly  Lipschitz-like at $\left(U_0,t_0 \right)$: then $U^{\prime}(t)=F\left(U(t),t    \right)$ in the whole $(a,b)$ and $U$ is analytic at $t_0$.
\label{raccordo}
\end{proposition}
{\bf Proof:} note that, by assumption, $F$ is analytic in a neighbourhood of $\left(U(t_0),t_0    \right)$, hence we have
$$
\lim_{t\to t_0}U^{\prime}=\lim_{t\to t_0}F\left(U(t),t    \right)=F\left(U(t_0),t_0     \right). 
$$
Thus $U$ is derivable in $t_0$ and $U^{\prime}(t_0)=F\left(U(t_0),t_0     \right)$; moreover $U$ equals, in a neighbourhood of $t_0$, the unique solution of the Cauchy's problem
$$
\cases{
U^{\prime}(t)=F\left(U(t),t    \right)\cr
U(t_0)=U_0;
}
$$
this fact proves analyticity at $t_0$.
\QUAN
\begin{proposition}\rm
Let $V$ be a $\ERRE^P$-valued analytic function on $(a,t_0)$ and suppose that $V^{\prime}(t)=F\left(V(t),t    \right)$ in $(a,t_0)$ , where $F$ is a uniformly  Lipschitz-like $\ERRE^P$-valued real analytic mapping in an open neighbourhood of 
$U(a,t_0)\times (a,t_0)$.
Moreover suppose that there exists 
$$
\lim_{t\to t_0^{-}}V(t)=V_0\in\ERRE^P
$$
and that $F$ is analytic at $\left(V_0,t_0    \right)$; then $V$ admits analytical continuation 
up to $t_0$.
\label{continuazione}
\end{proposition}
{\bf Proof:} the Cauchy's problem 
$$
\cases{
V^{\prime}(t)=F\left(V(t),t    \right)\cr
V(t_0)=V_0
}
$$
admits a unique analytic solution $W$ in a neighbourhood $I\left(t_0,r    \right)$ of $t_0$; by proposition \ref{raccordo}, the function
$U$ defined by setting
$$
U(t)=\cases{V(t) & if $t\in(a,t_0)$\cr
W(t) & if $t\in [t_0,t_0+r)$
}
$$
is a solution of the equation 
$U^{\prime}(t)=F\left(U(t),t    \right)$ in the whole $(a,t_0+r)$ and it is analytic at $t_0$, ending the proof.
\QUAN
\begin{theorem}
\TTT
Let $U$ be a $\ERRE^P$-valued analytic
function on $(a,t_0)$, solution of the equation $$U^{\prime}(t)=F\left(U(t),t    \right)$$
in $(a,t_0)$ , where $F$ is a  uniformly  Lipschitz-like $\ERRE^P$-valued real analytic mapping in a neighbourhood $G$ of 
$U(a,t_0)\times (a,t_0)$.
Moreover suppose that there exists a strictly increasing sequence $\{t_h\}\longrightarrow t_0$, such that 
$$
\cases
{
\lim_{h\to\infty}U(t_h)=U_0\quad\heartsuit\cr
\left(U_0,t_0 \right)\in G;
}
$$
then there holds
$$
\lim_{t\to t_0^{-}}U(t)=U_0\in\ERRE^P,
$$
hence then $U$ admits analytical continuation up to $t_0$.
\label{sequenza}
\end{theorem}
{\bf Proof:} now there exist neighbourhoods $I$ of $t_0$ and $J$ of $U_0$ such that $F$ is analytic in $\overline{I\times J}$.
Set
$$
M=\sup_{I\times J}\vert F\left(U,t    \right)\vert.
$$
We claim that for each (small) $\varepsilon>0$ there exists $j$ such that 
$$
\clubsuit\quad\left(t_j<t<t_0    \right)\Longrightarrow \vert U(t)-U_0\vert<\varepsilon.
$$
Let $\varepsilon>0$ be such that $D(U_0,\varepsilon)\subset J$: then, by $\heartsuit$, there exists $j\in\ENNE$ such that $t_j\in I$ and 
$$
\cases{
\vert t_j-t_0 \vert< \varepsilon / 4M\cr
\vert U(t_j)-U_0) \vert <\varepsilon/2.}
$$
Let's show that $$
\triangle\quad\vert U(t)-U(t_j)\vert<\varepsilon/2:
$$
suppose, on the contrary, that 
$$
\left\{
t\in \left(t_j,t_0    \right):\vert U(t_j)-U_0\vert \geq\varepsilon/2 \right\}\not=\emptyset,
$$
and let $\tau=\inf E$: by continuity, $\vert U(t_j-U_0) \vert \geq\varepsilon/2$, hence $\tau>t_j$.

On the other hand, if $t_j<\xi<\tau$, there holds $\vert U(\xi)-U(t_j)\vert<\varepsilon/2$, hence $u(\xi)\in J$ and $\vert U^{\prime}(\xi)\vert=\vert F\left(U(\xi),\xi    \right)\vert\leq M$.
Now
\begin{eqnarray*}
\varepsilon/2 &\leq& \vert U(\tau)-U(t_j)\vert\\
&=&\left\vert\left(
U^1(\tau)-U^1(t_j)...U^P(\tau)-U^P(t_j)
\right)\right\vert\\
&=&\left\vert\left(
U^{1\ \prime}(\xi_1)\left(\tau-t_j    \right)...
U^{P\ \prime}(\xi_P)\left(\tau-t_j    \right)
\right)\right\vert,
\end{eqnarray*}
where the $\xi_n$'s are suitable points in $(t_j,\tau)$.
Hence $$
\varepsilon/2<M\vert t_0-t_j \vert<\varepsilon/4,
$$
which is a contradiction.
By $\triangle$ just proved, if $t_j<t<t_0$, 
$$
\vert U(t)-U_0 \vert\leq \vert U(t)-U(t_j)\vert+\vert U(t_j)-U_0\vert<\varepsilon:
$$
this ends the proof.
\QUAN

The following corollary applies the previous real-variable results to the 
complex domain:
\begin{corollary}\rm
let $\left({\cal V},U\right)$ be a $\CI^N$-valued holomorphic function element, solution of the equation $U^{\prime}(z)=G\left(U(z),z    \right)$ in ${\cal V}$ , where $G$ is a  uniformly  Lipschitz-like $\CI^P$-valued holomorphic mapping in a neighbourhood of 
$U\left({\cal V}\right)\times {\cal V}$.
Pick now $z_0\in bd({\cal V})$ and suppose
that there exists an embedded real analytic regular curve $\gamma:[0,1]\longrightarrow {\cal V}\cup\{z_0\}$ with endpoint $z_0$ and a strictly increasing sequence $\{t_h\}\longrightarrow 1$, such that there exists $\lim_{h\to\infty}U\circ\gamma(t_h)=U_0$ and that $G$ is holomorphic at $\left(U_0,z_0    \right)$;
then $U$ admits analytical continuation up to $z_0$.
\label{complseq}
\end{corollary}
{\bf Proof:} let ${\cal W}\subset\CI^N$ be an open set such that $G$ is holomorphic in ${\cal W}\times{\cal V}$ and consider the mappings
$$
\cases{\Phi_{\Delta}:\CI^N\longrightarrow\ERRE^{2N}\cr
\left(z^1...z^N    \right)\longmapsto\left(\Re(z^1),\Im(z^1)......\Re(z^N),\Im(z^N)    \right)
}
$$

$$
\cases{
H:\Phi_{\Delta}\left({\cal W}\right) \times [0,1]    
\longrightarrow\ERRE^{2N}\cr
H\left(x^1,y^1......x^N,y^N,t    \right)=\Phi_{\Delta}\left[\GAMMA^{\bullet}(t)
G \left(x^1+\hbox{\bf i}y^1...x^N+\hbox{\bf i}y^N,\gamma(t)   \right)\right],
}
$$
and set $
V=\Phi_{\Delta}\circ U\circ\gamma
$; then  
\begin{eqnarray*}
\VII^{\bullet}(t)&=&\frac{d}{dt}\left[\Phi_{\Delta} \circ U\circ\gamma    \right](t)\\
&=&\Phi_{\Delta}\left(\frac{d}{dt}\left(U
\circ\gamma(t)\right)\right)\\
&=&\Phi_{\Delta} \left(U^{\prime}\left(\gamma(t)    \right)\GAMMA^{\bullet}(t)    \right)\\
&=&
\Phi_{\Delta}\left[G\left(U\circ\gamma(t),\gamma(t)    \right)\GAMMA^{\bullet}(t)    \right]\\
&=&H\left(\Phi_{\Delta}\circ U\circ\gamma(t),t    \right)\\
& &\  \\
&=& H\left(\left(V(t),t    \right)\right),
\end{eqnarray*}
hence theorem \ref{sequenza} may be applied since $\{V(t_h)\}\longrightarrow \Phi_{\Delta}(U_0)$ and $H$ is real analytic at $\left(\Phi_{\Delta}(U_0),1    \right)$: then $V$ admits analytical continuation up to $1$, that is to say, $U$ is analytically continuable up to $z_0$.
\QUAN
\vfill\eject

\section{An extension of Painlev\'e's theorem}
\subsection{Introduction}
In this section we pursue the goal of extending 
the well known Painlev\'e's determinateness theorem (see e.g. \cite{hille}, theorem 3.3.1);
let ${\cal N}$ be one of the following objects:
\begin{itemize}
\item 
the empty set $\emptyset$;
\item
a (not necessarily connected) algebraic curve in $\CI^2$  
\end{itemize}
and $(R,p)$ a proper Riemann domain over $\CI^2\setminus{\cal N}$.

Introduce complex-valued holomorphic functions $P$ and $Q$
on $R$, with $Q$ not everywhere vanishing,
so that $P/Q$ is a meromorphic function on $R$; if $X_0\in R$, $(u_0,v_0)=p(X_0)$ and $\eta$ is a local inverse of $p$ in a neighbourhood of $(u_0,v_0)$, consider the following Cauchy's problem at $v_0$:
\begin{equation}
\cases{
\displaystyle
{
\u_{v_0}^{\prime}
=\frac{P\circ\eta}{Q\circ\eta}
(\u_{v_0}\times\id_{v_0})
}
\cr
\displaystyle
{
\u_{v_0}(v_0)=u_0.
}
\cr
\ 
}
\label{cauchya}
\end{equation}
We wonder about the analytical continuability of its (unique) germ solution: this is in general quite a difficult argument, but if $P$ and $Q$ resemble, in a sense to be precised, 'polynomials', something more could be asserted.
We shall state and prove a theorem of determinatess, representing an intermediate step towards the goal of proving analytical continuability theorems.
\subsection[Set of hypotheses for Painlev\'e's theorem]{\label{secpain}Set of hypotheses for theorem \ref{painleve'} } \label{pain2}
\begin{itemize}
\item (H1)  $R$, $p$ and ${\cal N}$ are as in {\bf(2.2.1)};
\item (H2)  the curve $p\left(Q^{-1}(0)\right)$ is an algebraic curve, or, alternatively, the empty set.
\end{itemize}
\subsection
[Notations for Pailev\'e's theorem ]
{Notations for theorem \ref{painleve'} }
\begin{definition}\rm
The complex number $k\in\CI_v$ is a {\bf horizontal flatness point} for the polynomial $\left(u,v    \right)\longmapsto Z\left(u,v    \right)$ (or for the curve $Z\left(u,v    \right)=0$) if $\left(v-k    \right)$ divides
$Z$.
\label{hflat}
\end{definition}
\subsection{Painlev\'e's theorem extended}
\begin{theorem}
\TTT
Let $\left({\cal U},u    \right)$ be a holomorphic function element, re\-pre\-sen\-ta\-tive of the unique germ solution of (\ref{cauchya})
and $A$ the set of all horizontal flatness points of $p\left(Q^{-1}(0) \right)\bigcup{\cal N}$; let $v_1\in\CI\setminus A$ and suppose $\gamma:[0,1]\longrightarrow\CI\setminus A$ is an embedded rectifiable analytic arc connecting $v_0$ and $v_1$; then, if an analytical continuation $\omega$ of $u$ may be got along $\gamma\vert_{[0,1)}$, i.e if there exists a function element $({\cal V},\omega)$ such that
$$
\cases{
\displaystyle{{\cal V}\supset\gamma\left([0,1)   \right)}\cr
\displaystyle{[\omega]_{v_0}=[u]_{v_0}}\cr
\displaystyle{v_1\in bd({\cal V})}\cr
}
$$
and all the hypotheses of the above set are fulfilled, there exists $\displaystyle\lim_{v\to v_1, v\in {\cal V}}\omega(v)$ in $\PI^1$.
\label{painleve'}
\end{theorem}
{\bf Proof:} suppose, on the contrary, that 
$\displaystyle\lim_{v\to v_1, v\in {\cal V}}\omega(v)$ in $\PI^1$ does not exist.

For every $\nu\in\CI_v$, set 
$$
W_{\nu}=
pr_1
\left(
\left(
\left(
p(Q^{-1}(0)
\right)\bigcup{\cal N}
\right)
\bigcap 
(\CI\times\{\nu\})
\right).
$$

By (H2) $W_{v_1}$ is finite or empty; the case it is empty is trivial; otherwise, let's say that $W_{v_1}=\{u_k\}_{k=1....q}$.

For each $k$, set 
$$
\cases{
D_k=D(u_k,\varepsilon)\cr
T_k=bd(D_k),
}
$$
for  $\varepsilon>0$ small enough: then there exists $\varrho>0$ such that 
$$
v\in D\left(v_1,\varrho\right)\Longrightarrow
W_{v}
\subset\bigcup_{k=1}^{q}D_k.
$$ 
This implies that there exists $M>0$ such that 
$$
X\in p^{-1}\left( \bigcup_{k=1}^{q}T_k \times D\left(v_1,\varrho    \right)   \right)\Longrightarrow \vert\frac{P}{Q} (X) \vert\leq M.
$$
Take now, as it is possible,  $R>0$ so large that
$$
X\in p^{-1}\left(bd(D(0,R))\times D\left(v_1,\varrho    \right)   \right)\Longrightarrow \vert\frac{P}{Q} (X) \vert\leq M.
$$
Set now $\Theta_u=\overline{D(O,R)}\setminus\bigcup_{k=1}^{q}D_k$: this is a compact set and, for each $v\in D\left(v_1,\varrho    \right)  $, $p^{-1}\left(\CI_u\times \{v\}    \right)$ is a Riemann surface, hence, by maximum principle,
$$
X\in p^{-1}\left(\Theta_u\times \{v\}  \right)\Longrightarrow \vert\frac{P}{Q} (X) \vert\leq M,
$$
and, by the arbitrariness of $v$ in $D\left(v_1,\varrho    \right)  $,
$$
X\in p^{-1}\left(\Theta_u\times D\left(v_1,\varrho    \right)   \right)\Longrightarrow \vert\frac{P}{Q} (X) \vert\leq M.
$$
Now we have assumed that $\omega$ does not admit limit as $v\longrightarrow v_1$ in ${\cal V}$, hence there exists a sequence $\{v_i\}\subset {\cal V}$ converging to $v_1$ such that $\{\omega(v_i)\}\subset \Theta_u$.
Without loss of generality, we may suppose that $\{v_i\}\subset D\left(v_1,\varrho/2   \right) $; since $p$ is proper, $p^{-1}\left( \Theta_u\times\overline{D\left(v_1,\varrho/2   \right) }   \right)$ is compact, hence we may extract a convergent subsequence $\Omega_k$ from $p^{-1}\{(\omega(v_i),v_i)\}$: let $\Omega$ be its limit and ${\cal H}$ a neighbourhood of $\Omega$ such that $p\vert_{\cal H}$ is a biholomorphic mapping: the Cauchy's problem at $v_1$
$$
\cases{\displaystyle{
\u_{v_1}{}^{\prime}=\frac{P\circ\left(p\vert_{\cal H}    \right)^{-1}}{Q\circ\left( p\vert_{\cal H}      \right)^{-1}}(\u_{v_1}\times\id_{v_1})}\cr
\displaystyle{
\u_{v_1}(v_1)=pr_1\left(p(\Omega)    \right)
}
\cr
\ 
}
$$
admits a unique holomorphic solution $\u_{v_1}$.

By corollary \ref{complseq}, a representative holomorphic mapping  element of $\u_{v_1}$ is connectible with $\left({\cal V},\omega    \right)$, hence with $({\cal U},u)$: this eventually implies that $\lim_{v\to v_1, v\in {\cal V}}\omega(v)=pr_1\left(p(\Omega)    \right)$, which is a contradiction.
\QUAN

The main difference between theorem \ref{painleve'} and theorem 3.3.1 of \cite{hille} is that, in the latter, the proper Riemann domain $\left(R,p    \right)$ is restricted to be of the particular form $\left(S\times\CI    \right)$, for some algebraic Riemann surface $S$.
\vfill\eject
\chapter{Paths} 
\vfill\eject
\section{Introduction}
Let's start with a slight reformulation of the notion of path: to achieve this goal, we adopt the point of view according to which a 'path' or even a 'curve' are analytical continuations of some initial germs, generally yielded by local solutions of systems of differential equations.

Moreover our approach will be inclusive of all types of singularities encountered in the chapter about analytical continuation; on the other side, we agree to use the term 'path' to designate functions defined on one-dimensional manifolds (Riemann surfaces), whereas 'curve' will denote, informally speaking, a collection of local relations among coordinates: in this setting a germ of curve will hence be regarded as a germ of a subset (see \cite{gunros}).

Another aspect of a path, not of negligible importance, will be its velocity field: to define it we shall need a vector field on its domain of definition, to be related with the natural derivation field $d/dz$ on $\CI$: this approach is closely related to the analytical continuation procedure.
\section{Main definitions and theorems}

Let $\M$ be a connected complex manifold: in the continuation, abusing language but following Wells (see e.g. \cite{wells} or \cite{gunros}), we shall name $T\M$ (resp.$T^{*}\M$) its holomorphic tangent (resp. cotangent) bundle and, more generally, ${\cal T}_{r}^{s}\M$ its holomorphic r-covariant and s-contravariant tensor bundle; as usual, $\Pi\colon{\cal T}_{r}^{s}\M\longrightarrow\M$ will appoint the natural projection.
\begin{definition}\rm
A {\bf closed hypersurface} ${\cal F}$ in a complex manifold $\M$ is a closed subset such that there exists 
a maximal atlas $\{U_n\}$ for $\M$ and, for each $n$, a holomorphic function $\Psi_n$, not vanishing everywhere,
such that $$U_n\bigcap{\cal F}=\{X\in U_n : \Psi_n(X)=0\} .$$
\label{clohyp}
\end{definition}

\subsection{Meromorphic sections}
The following definition is adapted from \cite{oneill}, definition 2.4 and lemma 2.5:
\begin{definition}\rm
let ${\cal E}$ be a closed hypersurface in $\M$, $\N$ another connected complex manifold and $F\in{\cal O}(\M,\N)$:
{\bf an ${\cal E}$-meromorphic section} of 
${\cal T}_{r}^{s}\N$ {\bf over} $F$ is a 
holomorphic section $\Lambda$ of ${\cal T}_{r}^{s}\N$ over $F\vert_{\M\setminus
{\cal E}}$ such that
\begin{itemize}
\item
$\pi\circ\Lambda$ admits analytical continuation up to the whole $\M$ (coinciding with $F$);
\item
for every $p\in{\cal E}$ and every coordinate system $\left({\cal U},(z^1...z^n)\right)$ around
$F(p)$, there exists a neighbourhood $U$ of $p$ and $r\cdot s$ pairs of $\CI-$valued 
holomorphic functions $\phi_{i_1...i_r},\
\psi_{l_1...l_s}$, with $\psi_{l_1...l_s}\not=0$ on 
$U\setminus{\cal E}$, such that
$$
\Lambda
\left( dz^{l_1}...dz^{l_s},
\frac
{\partial}
{\partial z^{i_1}}...
\frac{\partial}
{\partial z^{i_r}}
\right)=
\frac{\phi_{i_1...i_r}}
{\psi_{l_1...l_s}}.
$$
\end{itemize}
\label{section}
\end{definition}
\begin{remark}\rm
If $\M=\N$ and $F=id_{\M}$ we are dealing with ${\cal E}$-meromorphic section of ${\cal T}_{r}^{s}\M$ tout-court; if
${\cal E}=\emptyset$ we shall talk about  holomorphic sections of (holomorphic) tensor bundles.
\end{remark}
\subsection{Paths}
\begin{definition}
\rm
A {\bf path} in $\M$ is a quintuple $Q_{\M}=\left(S\bigcup S^{b},\pi,j,F,\M\right)$,
where
\begin{itemize}
\item
$S\bigcup S^{b}$ is a connected Riemann surface with boundary $S^{b}$;
\item
$\pi\colon S\bigcup S^{b}
\longrightarrow\PI^1$ is a branched covering of $S$ over $\pi(S)$, with set of branch points $P$ and $\pi\in {\bf C}^0\left( S\bigcup S^{b},
\PI^1   \right)$;
\item
$F\in{\cal A}_{\Sigma}(S,\M)\bigcap{\bf C}^0\left( S\bigcup S^{b},\M   \right)$, where $\Sigma$ is a discrete set in $S$ and $\Sigma\bigcap S^b=\emptyset$;
\item
$U\subset{\CI}$ is an open set wich admits a holomorphic (hence open) immersion $j\colon U\longrightarrow S\setminus\Sigma$ such that $\pi\circ j=id\vert_{U}$.
\end{itemize}
The path $Q_{\M}$ will be told to be
\begin{itemize}
\item 
$z_0$-starting at $m\in\M$ if $z_0\in U$ and $F\circ j (z_0)=m$;
\item
{\bf without boundary}, or {\bf without logaritmic singularities}, if $S^b=\emptyset$
\end{itemize}
\label{path1}
\end{definition}
In the following, 'starting' will mean
$0$-starting.

\subsection{The velocity field}
We are now turning to define the {\bf velocity field} of a path $Q_{\M}$: it will be defined, except for points in the boundary, as a suitable meromorphic section over $F$ of the holomorphic tangent bundle $T\M$: to achieve this purpose, we need to lift the vector field ${d}/{dz}$ on $\CI$ with respect to $\pi$; of course, in general, contravariant tensor fields couldn't be lifted.

Notwithstanding, we may get through this obstruction by keeping into account that $\CI$ and $S$ are one-dimensional manifolds and allowing the lifted vector field to be meromorphic: these matters are fathomed in next statements: recall that $P$ is the set of branch points of $\pi$.

\begin{lemma}\rm
there exists a unique $P$-meromorphic vector field $\widetilde{{d}/{dz}}$ on $S$ 
such that, for every $r\in S\setminus P$, $\pi_{*}\vert_r\left(\widetilde{{d}/{dz}}\vert_{r}    \right)=\left({d}/{dz}\right)\vert_{\pi(r)}$.
\label{lift}
\end{lemma}
{\bf Proof:} the {\it 1-form} $dz$ and the metric $dz\odot dz$ are covariant tensor fields, hence they may be pulled back with respect to $\pi$, getting a holomorphic 1-form $\omega=\pi^{*}dz$ and a holomorphic metric $\Lambda=\pi^{*}\left(dz\odot dz    \right)$ on $S$: the latter is nondegenerate on $S\setminus P$, hence it establishes an isomorphism between its holomorphic cotangent  and  tangent bundles.

Call $V$ the holomorphic vector field corresponding to $\omega$ in the above isomorphism: we claim that $V=\widetilde{{d}/{dz}}$ on $S\setminus P$.
To show this fact, we explicitely compute the components of $V$ with respect to a maximal atlas ${\cal B}=\left\{\left(U_{\nu},\zeta_{\nu}    \right)\right\}$ for $S\setminus P$: let
$$
\cases{\omega_{(\nu)\ 1}=\omega\left(\partial/\partial \zeta_{(\nu)}    \right)\cr
g_{(\nu)\ 11}=\Lambda\left(\partial/\partial \zeta_{(\nu)},\partial/\partial \zeta_{(\nu)}    \right);
}
$$
then, set $\displaystyle{V_{(\nu)}^1=\frac{\omega_{(\nu)\ 1}}{g_{(\nu)\ 11}}}$ the collection $\left\{\left( {\cal U_{\nu}}, V_{(\nu)}^1  \right)\right\}$ of open sets and holomorphic functions is such that, on overlapping local charts $\left(U_a,\zeta_a    \right)$ and $\left(U_b,\zeta_b   \right)$, we have:
\begin{eqnarray*}
V_{(a)}^1&=&\frac{\omega_{(a)\ 1}}{g_{(a)\ 11}}\\
&=&\frac{\omega_{(b)\ 1}\frac{\displaystyle d\zeta_{(b)}}{\displaystyle d\zeta_{(a)}}}{g_{(b)\ 11}\left(\frac{\displaystyle d\zeta_{(b)}}{\displaystyle d\zeta_{(a)}}    \right)^2}\\
&=& \frac{\omega_{(a)\ 1}}{g_{(a)\ 11}}\frac{\displaystyle d\zeta_{(a)}}{\displaystyle d\zeta_{(b)}}\\
&=&V_{(b)}^1
\frac{d\zeta_{(a)}}{d\zeta_{(b)}},
\end{eqnarray*}
that is to say, that collection defines a holomorphic, simply contravariant tensor field, i.e. a holomorphic vector field.

Now for every $r\in S\setminus P$, 
$$
\pi_{*}\vert_r\left(\widetilde{{d}/{dz}}\vert_{r}    \right)=\left({d}/{dz}\right)\vert_{\pi(r)}\iff
dz\vert_{\pi(r)}\left(\pi_{*}\vert_r \widetilde{{d}/{dz}}\vert_{r}    \right)=1.
$$
We prove the right side to be true:
\begin{eqnarray*}
dz\vert_{\pi(r)}
\left( 
\pi_{*}\vert_r
\widetilde{{d}/{dz}}\vert_{r}
\right)=\pi^*dz\vert_r
\left(
\widetilde{{d}/{dz}}\vert_{r}
\right)=\\
=\pi^*dz\vert_r
\left(
\left.
\frac
{1}
{dz\vert_{\pi(r)}(\pi_*\partial/\partial\zeta\vert_r)} 
\frac{\partial}{\partial\zeta}
\right\vert_r
\right)
=\\
=\frac{\pi^*dz\vert_r (\partial/\partial\zeta\vert_r) }
{dz\vert_{\pi(r)}(\pi_*\partial/\partial\zeta\vert_r) 
}
=1,
\end{eqnarray*}
hence the left one is true too,
proving the asserted. 

Let's prove now that $\widetilde{d/dz}$ may be extended to a meromorphic vector field on $S$: let indeed $p\in P$: then
we could find local charts $\left(U,\psi    \right)$ around $p$ ($\psi : U\longrightarrow\CI_u$)  and $\left(V,\phi   
\right)$ around $\pi(p)$ ($\phi : V\longrightarrow\CI_v$) , and an integer $N>0$, not depending on $\left(U,\psi    \right)$ and $\left(V,\phi    \right)$  such that $$\phi\circ\pi\circ\psi^{-1}(u)=u^N.$$
Now we have
\begin{eqnarray*}
(\psi^{-1\ *}\pi^{*}\phi^{*}(dw)(d/du))(u)=dw(\phi_{*}\pi_{*}\psi_{*}^{-1}(d/du)\vert_u)\\
=dw((\phi\pi\psi^{-1})^{\prime}(d/dw))=Nu^{N-1};
\end{eqnarray*}
 but $\phi$ and $\psi$ are charts, hence $\pi^{*}dz $ itself is vanishing of order $N-1$ at $p$; but, as we have already proved, $\pi_{*}\vert_r\left(\widetilde{{d}/{dz}}\vert_{r}    \right)=\left({d}/{dz}\right)\vert_{\pi(r)}$ on $U\setminus \{p\}$ and, consequently,
$$(\pi^{*}dz)(\widetilde{d/dz})=dz(\pi_{*}\widetilde)=dz(d/dz)=1$$ on $U\setminus \{p\}$, hence, by Riemann's removable singularity theorem, on the whole $U$.
Now, in local coordinates, 
$$
\cases{\displaystyle(\pi^{*}dz)=\alpha d\phi\cr
\displaystyle\widetilde{d/dz}=y\,\,\frac{\partial}{\partial\phi},\cr
\cr
}
$$
where $\alpha$ is a holomorphic function on $U$, vanishing of order $N-1$ at $p$ and $y$ is a holomorphic function on $U\setminus \{p\}$.
By the argument above, $y\alpha=1$, hence $y$ has a pole of order $N-1$ at $p$: this holds for each isolated point in $P$, hence the meromorphic behaviour of $\widetilde{d/dz}$ is proved for the special atlas ${\cal B}={\cal A}\bigcup\{(U,\phi)_p\}_{p\in P}$, hence for every atlas on $S$.
The uniqueness of $\widetilde {{d}/{dz}}$ as a meromorphic vector field follows from its uniqueness 
as a {\it holomorphic} vector field in a neighbourhood of a regular point of $\pi$ in $S$.
\QUAN

\begin{definition}\rm
A {\bf finite-velocity  point} of a path $Q_{\M}=\left(S,\pi,j,F,\M\right)$ 
is a point $r\in S$ such that $\widetilde {{d}/{dz}}$ is holomorphic at $r$.
\label{path3}
\end{definition}

We are ready to define the velocity field (this will be done in definition \ref{vf}): let at first be $r$ a finite-velocity  point of $Q_{\M}$; since $\widetilde{{d}/{dz}}$ is holomorphic at $r$, we could define the holomorphic velocity at $r$ as $V_r=F_{*}\left((\widetilde{{d}/{dz}})\vert_r\right)$.
\begin{lemma}\rm
The mapping 
\begin{eqnarray}
& & V\left(Q_{\M}\right)\colon S\setminus P \longrightarrow T\M\nonumber\\
& & r\longmapsto \left(F,F_{*}\vert_r
\left(\widetilde{\frac{d}{dz}}\vert_r\right)\right)\nonumber
\end{eqnarray}
may be extended to a P-meromorphic section  of $T\M$ over $F$.
\label{velo_def}
\end{lemma}
{\bf Proof:} trivially $\Pi\circ V\vert_{R\setminus P}=F\vert_{R\setminus P}$.
Let's show the meromorphic behaviour of $V$: if $p\in P$ there is a neighbourhood $U$ of $p$ such that, for every local chart $\zeta\colon U\longrightarrow\CI_w$ there exist holomorphic functions $f,g\in{\cal H}\left(\zeta(U)\right)$ such that 
$$\widetilde{\frac{d}{dz}}\vert_{\zeta_{-1}(U)}=\zeta^{-1}_{*}\left(\frac{f}{g}(w)\frac{d}{dw}\vert_w\right);$$
moreover, for every local chart $\Psi=\left(u^1...u^m,du^1...du^m\right)$ in $T\M$ we obtain
\begin{eqnarray}
& &\Psi\circ V \circ \zeta^{-1}(w)=\nonumber\\
&=& \Psi\circ \left(F\circ \zeta^{-1}(w),F_{*}\vert_{\zeta^{-1}(w)}\left(\widetilde{\frac{d}{dz}}\vert_{\zeta^{-1}(w)} \right)\right) \nonumber\\
&=& \Psi\circ \left(
F\circ \zeta^{-1}(w),F_{*}\vert_{\zeta^{-1}(w)}\zeta^{-1}_{*}\left(\frac{f}{g}(w)\frac{d}{dw}\vert_w\right)\right) \nonumber\\
&=& \Psi\circ \left( F\circ \zeta^{-1}(w),\frac{f}{g}(w)\frac{d}{dw}(F\circ\zeta^{-1})(w)\right) \nonumber\\
&=& \left(u^1\circ F\zeta^{-1}(w)...u^m\circ  F\zeta^{-1}(w),\right.\nonumber\\ 
& &
\left.\frac{f}{g}(w)\frac{d}{dw}\left(u^1\circ F\circ\zeta^{-1}\right)(w)...\frac{f}{g}(w)\frac{d}{dw}\left(u^m\circ
F\circ\zeta^{-1}\right)(w))\right)  \nonumber
\end{eqnarray}
\QUAN
\begin{definition}\rm
The {\bf velocity field } of a path $Q_{\M}=(S,\pi,j,F,\M)$ is the meromorphic mapping
\begin{eqnarray}
& & V\left(Q_{\M}\right)\colon S\setminus P \longrightarrow T\M\nonumber\\
& & r\longmapsto \left(F,F_{*}\vert_r
\left(\widetilde{\frac{d}{dz}}\vert_r\right)\right)\nonumber
\end{eqnarray}
\label{vf}
\end{definition}

\subsection{Completeness}
We eventually yield a definition of completeness:
\begin{definition}
\rm
\label{completezza}
A $\M$-valued path (with or without logarithmic singularities) $\displaystyle\left(S,\pi,j,F,\M    \right)$ is {\TTT complete} provided that $\PI^1\setminus\displaystyle\pi\left(S    \right)$ is a finite set in the Riemann sphere.
\label{complete}
\end{definition}
\vfill\eject

\chapter{Complex-Riemannian metric structures}
\vfill\eject 
\section{Introduction and main definitions}
The intuitive geometry of the real euclidean space $\ERRE^3$ can be easily brought back to its natural inner product, which allows basic geometrical operations, like measuring the length of a tangent vector, or angles between tangent vectors: Riemannian real geometry generalizes all this to 'curved' spaces, which is based on the concept of positive definite bilinear forms: weakening definiteness to nondegeneracy leads us in the realm of Lorentz geometry, originating from the problems posed by Einstein's general
relativity theory.

A less intuitive idea is that of starting from the basic geometry of $\CI^3$, meant as a 'complexification' of the usual real euclidean space to get formal extension of the geometric properties of real 'curved' manifolds.
Introducing this complex environment could allow us to hope to get able to handle some types of metrical singularities which naturally arise in dealing with real manifolds with indefinite metrics; it is soon seen that the nondegeneracy hypothesis itself should be dropped, since it would cut down the purport of this investigations, when combined with analyticity assumptions.

Wide treatises about semi-Riemannian or Lorentz  geometry are 
\cite{oneill} and \cite{beemerlich}; a different approach to 
holomorphic geometry could be found in \cite{manin}. Finally, we owe \cite{lebrun} for the definition of nondegenerate holomorphic metric.

Let now $\M$ be a complex manifold, ${\cal D}$ and ${\cal E}$ closed hypersurfaces in $\M$.

\begin{definition}\rm
{\TTT A  holomorphic (resp.${\cal E}$-meromorphic) metric on} $\M$ is a holomorphic (resp.${\cal E}$-meromorphic)  section 
$\Lambda\colon\M\longrightarrow {\cal T}_{0}^{2}\M$ which is symmetric, that is to say, for every $m\in\M$ and every pair
of holomorphic tangent vectors $V_m$ and $W_m$ at $m$, there holds 
$$\Lambda(m)\left( V_m, W_m  \right)=\Lambda(m)\left( W_m, V_m  \right) ;
$$
\begin{itemize}
\item
the {\bf rank} of $\Lambda$ at $p\in\M$ is the rank of the bilinear form $\Lambda(p)$;
\item
$\Lambda$ is {\bf nondegenerate} at $p$ if $rk(\Lambda(p))=dim(\M)$, {\bf degenerate} otherwise;
\item
if ${\cal D}$ is a hypersurface in $\M$ and $\Lambda$ is degenerate only on $\cal D$, we shall say that $\Lambda$ is ${\cal D}$-degenerate.
\end{itemize}
\label{rank_met}\label{hol_met}
\end{definition}
In the following we shall consider only metric which degenerate only on closed hypersurfaces.
\begin{remark}\rm
We emphasize that holomorphic and meromorphic metrics {\rm are not} pseudo-Riemannian ones on the underlying real manifolds (see \cite{lebrun}).

\end{remark}

\begin{definition}\rm
We say that $p$ is a {\bf metrically ordinary point} in $\M$ if $\Lambda$ is holomorphic and nondegenerate at $p$; a {\bf metrically extraordinary point} otherwise.
\label{reg_point}
\end{definition}
\begin{definition}\rm
Two meromorphic metrics $\Lambda$ and $H$ on $\M$ are 
\begin{itemize}
\item[(a)]
{\bf holo-conformal} if for every $p\in\M$ there exists a neighbourhood $U$ of $p$ and $F\in{\cal O}(U)$, 
never zero, such that $\Lambda\vert_U=F\cdot H\vert_U$
\item[(b)]
{\bf mero-conformal} if if for every $p\in\M$ there exists a neighbourhood $U$ of $p$ and $F\in{\cal M}(U)$, 
not everywhere vanishing, such that $\Lambda\vert_U=F\cdot H\vert_U$
\end{itemize}
\label{conf}
\end{definition}
\begin{lemma}\rm
Every meromorphic metric $H$ is locally mero-conformal to a holomorphic one, that is, for every $p\in\M$ there exists a neighbourhood $U$ of $p$ such that (\ref{conf}) (b) holds on $U$. 

\end{lemma}
{\bf Proof:} it is enough to prove this statement for meromorphic metrics on open sets in $\CI^N$: let 
$$\Lambda=\sum_{i\leq k}\frac{E_{ik}}{F_{ik}}\left( dx^i\otimes dx^k+dx^k\otimes dx^i   \right)$$
and $\psi=l.c.m.(F_{ik})$: then $\psi\cdot\Lambda$ is a holomorphic metric on $U$, meroconformal to $\Lambda$ by construction.
\QUAN
\begin{definition}\rm
\begin{itemize}
\item
A {\bf holomorphic Riemannian manifold} is a complex manifold endowed with a holomorphic metric ;
\item
a {\bf nondegenerate holomorphic Riemannian manifold} is a complex manifold endowed with a nondegenerate holomorphic metric ;
\item
a {\bf meromorphic  Riemannian manifold} is a complex manifold endowed with a 
meromorphic metric ;
\end{itemize}
\label{riemann}
\end{definition}
Thus, strictly speaking, all the above objects are pairs consisting in complex manifolds and  metrics, but we shall often understand metrics and denote them by teo only underlying complex manifolds.
\vskip0.5truecm
{\bf Short examples}
(We postpone more articulate examples after introducing geodesics
).
\vskip0.5truecm
\begin{itemize}
\item[(a)] 'The complex-euclidean space': endow $\CI^m$ with metric 
$$\Lambda=\sum_{i=1}^m dz^i\otimes dz^i :$$ such a metric structure is holomorphic, everywhere nondegenerate; note that $\Lambda\left(Z,Z    \right)=\sum_{i=1}^m Z^{i\ 2}$ is not a norm in the usual sense: for example, if $m=2$, then the (complex) vector subspaces generated by $\left(1,i    \right)$ or by $(1,-i)$ entirely consist of null vectors.
\item[(b)]
'Warped products': let $\left(\M,\Lambda    \right)$ and $\left(\N,H    \right)$ be meromorphic Riemannian manifolds, $F\in{\cal M}(\M)$;
$$\cases{
\varrho\colon \M\times\N\longrightarrow\M\cr
\sigma\colon \M\times\N\longrightarrow\N}
$$
be the canonical projection of their Cartesian product: then $\varrho^{*}\left(\Lambda    \right)+\left(F\circ    \right)\sigma^{*}\left(H    \right)$ is a meromorphic metric on $\M\times\N$.
Consider e.g. $\CI^2_{(u,v)}$ endowed with $\Lambda\left(u,v    \right)=du\otimes du+F(u)dv\otimes dv$: if $F$ is holomorphic then $\Lambda$ is a holomorphic metric, degenerate only on the locus $F(u)=0$; if $F$ is $P$-meromorphic, then $\Lambda$ is $P\times\CI$-meromorphic.
\item[(c)]
Let $H$ a noncostant doubly periodic function on $\CI$: then $H(z)dz\otimes dz$ defines a meromorphic metric on $\CI$, 
which may be clearly pushed down to the quotient torus.
\item[(d)]
Vice versa, a holomorphic metric on the torus ${\CI}/{\left(a,b    \right)}$ defines a doubly periodic holomorphic function on $\CI$ with periods $a,b$, hence a constant function.
This proves that the only holomorphic metrics on complex one-dimensional tori may be identified with complex constants.
\item[(e)]
Symmetric product of $1$-forms: let $\omega_1$ and $\omega_2$ be holomorphic (resp. meromorphic)
$1$-forms: then $\omega_1\odot\omega_2=\frac{1}{2}\left(\omega_1\otimes\omega_2+\omega_2\otimes\omega_1    \right)$ is a holomorphic
(resp. meromorphic) metric.
Symmetric product is commutative and distributive with respect to sum, but in general, not associative.

\item[(f)]
It is a consequence of general theory that $\CI\PI^N$ does not admit any holomorphic metric; on the other hand, it does
carry  meromorphic ones, which may be defined e.g. in the following way: start from the metric  $\Lambda(z)=\sum_{i,j}R_{ij}(z^1\cdots z^N)\,dz^i
\odot dz^j$ defined on $\CI^N $, where the $R_{ij}$ are rational functions; then 
extend $\Lambda$ to the other copies of $\CI^N $ which $\CI\PI^N$ consists of by pulling it back along the coordinate-changing mappings.
These metrics are  meromorphic by construction.
\end{itemize}
\vfill\eject

\section{The meromorphic Levi-Civita connection}
We begin this section by introducing the holomorphic Levi-Civita connection induced on a holomorphic nondegenerate Riemannian manifold by its metric structure: this is done in a quite similar way to that pursued in (real) differential geometry, apart from a slight difference, which naturally arises: the action of the Levi-Civita connection is defined at first on 'local' vector fields, producing local ones as well, then it is globalized as a collection of local operators.

Let now $\left(\M,\Lambda    \right)$ be a nondegenerate Riemannian holomorphic manifold, ${\cal A}$ a maximal atlas for $\M$, ${\cal U}\in{\cal A}$ a domain of a local chart.
Let also 
\begin{itemize}
\item
${\cal X}\left({\cal U}    \right)$ be the Lie algebra of holomorphic vector fields on ${\cal U}$;
\item
${\cal O}\left({\cal U}    \right)$ the ring of holomorphic functions on ${\cal U}$.
\end{itemize}

\begin{definition}\rm
A {\bf connection on ${\cal U}$} is a mapping 
$$D \colon{\cal X}\left({\cal U}    \right)\times{\cal X}\left({\cal U}    \right)\longrightarrow {\cal X}\left({\cal U}    \right)$$ such that:
\begin{itemize}
\item[(D1)]
$D_V W$ is ${\cal H}\left({\cal U}    \right)$-linear in $V$;
\item[(D2)]
$D_V W$ is $\CI$-linear in $W$;
\item[(D3)]
$D_V\left(fW    \right)=\left(Vf\right)W+fD_VW$ for every $f\in{\cal H}\left({\cal U}    \right)$;
\end{itemize}
\label{met_conn_germ}
\end{definition}

$D_V W$ is called the {\bf covariant derivative } of $W$ with respect to $V$
 in the connection $D$. By axiom (D1), $D_V W$ has tensor character in $V$, while axiom (D3) tells us that it is not a tensor in $W$.

Next step is showing that there is a unique connection characterized by two further properties, (D4) and (D5) below, namely being anti-Leibnitz like with respect to the Lie bracket operation and Leibnitz like wuth respect to the metric.
In the following we use the alternative notation $\left\langle V,W   \right\rangle$ instead of $\Lambda\left(V,W    \right)$.
\begin{lemma}\rm
Let ${\cal U}$ be an open set belonging to a maximal atlas ${\cal A}$ for the nondegenerate holomorphic Riemannian manifold $\M$. If $V\in{\cal X}\left({\cal U}    \right)$, let $V^*$ be the holomorphic one-form on ${\cal U}$ such that $V^*(X)=\left\langle V,X   \right\rangle$ for every $X\in{\cal X}\left({\cal U}    \right)$: then the mapping $V\mapsto V^*$ is a ${\cal O}$-linear isomorphism from ${\cal X}\left({\cal U}    \right)$ to ${\cal X}^*\left({\cal U}    \right)$.
\label{vefo}
\end{lemma}
{\bf Proof:} since $V^*$ is ${\cal O}$-linear, it is in fact a one-form, and $V\mapsto V^*$ is ${\cal O}$-linear too.
We claim:
\begin{itemize}
\item[(a)] if $\left\langle V,X   \right\rangle=\left\langle W,X   \right\rangle$ for every $X\in{\cal X}\left({\cal U}    \right)$ then $V=W$;
\item[(b)] given any one-form $\omega\in{\cal X}^*\left({\cal U}    \right)$ there is a uique vector field $V\in{\cal X}\left({\cal U}    \right)$ such that $\omega(X)=\left\langle V,X   \right\rangle$ for every $X\in{\cal X}\left({\cal U}    \right)$.
\end{itemize}
Let $U=V-W$: then (a) is proved if we show that 
$$
\left(\left\langle U_p,X_p   \right\rangle=0\  
\hbox{\rm for every } X\in{\cal X}\left({\cal U}    \right)
\hbox{\rm and }
p\in{\cal U}    \right)\Longrightarrow U=0.
$$
This simply follows from the nondegeneracy of the metric tensor.

To prove (b), let $\left(z^1...z^N    \right)$ be local coordinates on ${\cal U}$.

Then $\omega=\sum_{i=1}^N\omega_i dz^i$; let $\{g_{ij}\}$ be the representative matrix of $\Lambda\vert_{\cal U}$ in $\left(z^1...z^N    \right)$: by nondegeneracy, it admits a holomorphic inverse matrix $\{g^{ij}\}$: set now $$V=\sum_{j=1}^N\left(\sum_{i=1}^N g^{ij}\omega_i   \right)\frac{\partial}{\partial z^j}; 
$$ we have
\begin{eqnarray*}
\left\langle V,X   \right\rangle &=& \left\langle \sum_{j=1}^N\left(\sum_{i=1}^N g^{ij}\omega_i   \right)\frac{\partial}{\partial z^j},\sum_{k=1}^N X^k \frac{\partial}{\partial z^k}  \right\rangle\\
&=& \sum_{ijk}g^{ij}\omega_i X^k g_{jk}\\
&=& \sum_{ik}\delta_k^i X^k \omega^i \\
&=&\sum_k X^k\omega^k\\
&=&\omega\left(X    \right)
\end{eqnarray*}
\QUAN
Corresponding vector fields and one-forms in the above isomorphism are told to be be metrically equivalent.
\begin{theorem}
\TTT
Let ${\cal U}$ be an open set belonging to a maximal atlas ${\cal A}$ for the nondegenerate holomorphic Riemannian manifold $\M$. There exists a unique connection $D$ on ${\cal U}$, called the Levi-Civita connection, such that:
\begin{itemize}
\item[(D4)]
$\left[V,W\right]=D_V W-D_W V$;
\item[(D5)]
$X\left\langle V,W   \right\rangle=\left\langle D_X V,W   \right\rangle + \left\langle V,D_X  W   \right\rangle $ for every $X,V,W\in {\cal X}\left({\cal U}    \right)$.
\end{itemize}
Moreover $D$ is characterized by the 'Koszul's formula':
\begin{eqnarray}
2\left\langle D_V W,X   \right\rangle &=& V\left\langle W,X   \right\rangle+
W\left\langle X,V   \right\rangle - X\left\langle V,W   \right\rangle\nonumber\\
&-&\left\langle V, [W,X]   \right\rangle
+\left\langle W, [X,V]   \right\rangle
+\left\langle X, [V,W]   \right\rangle,\nonumber
\end{eqnarray}
for every $X,V,W\in {\cal X}\left({\cal U}    \right)$.
\label{levicivita}
\end{theorem}
{\bf Proof:} suppose that such a connection exists: then 
\begin{eqnarray*}
& & V\left\langle W,X   \right\rangle+
W\left\langle X,V   \right\rangle
- X\left\langle V,W   \right\rangle\\
& & \qquad-\left\langle V, [W,X]   \right\rangle
+\left\langle W, [X,V]   \right\rangle
+\left\langle X, [V,W]   \right\rangle= \\
&=&
\left\langle D_V W,X   \right\rangle
+\left\langle W,D_V X   \right\rangle
+ \left\langle D_W X,V   \right\rangle
\\
& &\quad
+\left\langle X,D_W V   \right\rangle
-  \left\langle D_X V,W   \right\rangle
 -\left\langle V,D_X W   \right\rangle\\
& &
\qquad
-\left\langle V,D_W X - D_X W   \right\rangle
+\left\langle W,D_X V - D_V X   \right\rangle
+\left\langle X, D_V W-  D_W V   \right\rangle
\\
&=&
2\left\langle  D_V W,X  \right\rangle,
\end{eqnarray*}
proving uniqueness. 

To prove existence, define 
\begin{eqnarray}
F\left(V,W,X       \right) &=& V\left\langle W,X   \right\rangle+
W\left\langle X,V   \right\rangle - X\left\langle V,W   \right\rangle\nonumber\\
&-&\left\langle V, [W,X]   \right\rangle
+\left\langle W, [X,V]   \right\rangle
+\left\langle X, [V,W]   \right\rangle:\nonumber
\end{eqnarray}
for fixed $V,W\in{\cal X}\left({\cal U}    \right)$, the function
$X\mapsto F\left(V,W,X       \right) $ is ${\cal O}\left({\cal U}    \right)$-linear, indeed
$$
F\left(V,W,X+Y      \right)=F\left(V,W,X       \right)+F\left(V,W,Y       \right)
$$
trivially, and, if $\phi\in{\cal O}\left({\cal U}    \right)$
\begin{eqnarray*}
F\left(V,W,\phi X       \right) &=& V\left\langle W,\phi X   \right\rangle+
W\left\langle \phi X,V   \right\rangle -\phi  X\left\langle V,W   \right\rangle\\
&-&\left\langle V, [W,\phi X]   \right\rangle
+\left\langle W, [\phi X,V]   \right\rangle
+\left\langle \phi X, [V,W]   \right\rangle\\
&=& \left(V\phi     \right)\left\langle W,X   \right\rangle+\phi V\left\langle W,X+   \right\rangle +\left(W\phi     \right)\left\langle X,V   \right\rangle\\
&+& \phi W\left\langle X,V   \right\rangle-\phi X\left\langle V,W   \right\rangle\\
&-& \phi \left\langle V,[W,X]   \right\rangle-\left(W\phi     \right)\left\langle V,X   \right\rangle+\phi\left\langle W,[X,W]   \right\rangle\\
&-&\left(V\phi    \right)\left\langle W,X   \right\rangle+\phi\left\langle X,[V,W]   \right\rangle\\
&=&\phi F\left(V,W,X       \right).
\end{eqnarray*}
Thus the map $X\mapsto F\left(V,W,X       \right) $ is in fact a holomorphic one-form: hence, by lemma \ref{vefo}, there exists a unique holomorphic vector field on ${\cal U}$, which we denote by $D_V W$, such that $2\left\langle DV W,X   \right\rangle=F\left(V,W,X    \right)$ for every $X\in{\cal X}\left({\cal U}    \right)$.
Then the Koszul's formula holds and we are able to deduce \hbox{(D1)$\div$ (D5):}
\begin{eqnarray*}
\bullet\quad \hbox{(D1):}\quad 2\left\langle D_{\phi V}W,X   \right\rangle &=&
\phi V\left\langle W,X   \right\rangle+W\left\langle X,\phi
W    \right\rangle-X\left\langle \phi
V,W   \right\rangle\\
&-&\phi
\left\langle V,[W,X]   \right\rangle+\left\langle W,[X,\phi
V]   \right\rangle+\left\langle X,[\phi
V,W]   \right\rangle\\
&=&\phi
V\left\langle W,X   \right\rangle+(W\phi
)\left\langle X,V   \right\rangle+\phi
W\left\langle X,V   \right\rangle\\
&-&(X\phi
)\left\langle V,W   \right\rangle-\phi
X\left\langle V,W   \right\rangle-\phi
\left\langle V,[W,X]   \right\rangle\\
&+&\phi
\left\langle W,[X,V]   \right\rangle+(X\phi
)\left\langle W,V   \right\rangle+\phi
\left\langle X,[V,W]   \right\rangle\\
&-&(W\phi )
\left\langle  X,V  \right\rangle =2\phi
\left\langle D_V W,X   \right\rangle,
\end{eqnarray*}
where $\phi
\in {\cal O}\left({\cal U}    \right)$ and $V,W,X\in{\cal O}\left({\cal U}    \right)$;
finally $\left\langle D_{V+U} W,X   \right\rangle=\left\langle D_{V} W,X   \right\rangle+\left\langle D_{U} W,X   \right\rangle$ as a consequence of the bilinearity of $\left\langle \ \    \right\rangle$ and $\left[\ \     \right]$ with respect to holomorphic vector fields.

$\bullet\quad $(D2) is trivial;

\begin{eqnarray*}
\bullet\quad \hbox{(D3):}\quad 2\left\langle D_V(fW),X   \right\rangle&=&V\left\langle fW,X   \right\rangle+fW\left\langle X,V   \right\rangle-X\left\langle V,fW   \right\rangle\\
&-&\left\langle V,[fW,X]   \right\rangle+\left\langle fW,[X,V]   \right\rangle+\left\langle X,[V,fW]   \right\rangle\\
&=&(Vf)\left\langle W,X   \right\rangle+fV\left\langle W,X   \right\rangle+fW\left\langle X,V   \right\rangle\\
&-&(Xf)\left\langle V,W   \right\rangle-fX\left\langle V,W   \right\rangle-f\left\langle V,(W,X)   \right\rangle\\
&+&Xf\left\langle V,W   \right\rangle+f\left\langle W,[X,V]   \right\rangle+f\left\langle X,[V,W]   \right\rangle\\
&+&(Vf)\left\langle X,W   \right\rangle=2\left\langle VfW+fD_V W,X   \right\rangle
\end{eqnarray*}
for $f\in {\cal O}\left({\cal U}    \right)$ and $V,W,X\in{\cal O}\left({\cal U}    \right)$;
thus $D_V(fW)+(Vf)W+fD_V W$.

\begin{eqnarray*}
\bullet\quad \hbox{(D4):} & &\quad 2\left\langle D_V W-D_W V,X    \right\rangle =
V\left\langle W,X   \right\rangle+
W\left\langle X,V   \right\rangle - X\left\langle V,W   \right\rangle\\
&-&\left\langle V, [W,X]   \right\rangle
+\left\langle W, [X,V]   \right\rangle
+\left\langle X, [V,W]   \right\rangle\\
&-&W\left\langle V,X   \right\rangle-
V\left\langle X,W   \right\rangle + X\left\langle W,V   \right\rangle\\
&+&\left\langle W, [V,X]   \right\rangle
-\left\langle V, [X,W]   \right\rangle
-\left\langle X, [W,V]   \right\rangle\\
&=& \left\langle X,[V,W]   \right\rangle-\left\langle X,[W,V]   \right\rangle\\
&=& 2\left\langle [V,W],X   \right\rangle;
\end{eqnarray*}
\begin{eqnarray*}
\bullet\quad \hbox{(D5):}
& &\quad 2\left(\left\langle D_X V,W   \right\rangle+\left\langle V,D_X W   \right\rangle    \right) =\\
&=& X\left\langle V,W   \right\rangle+V\left\langle W,X   \right\rangle -W\left\langle X,V   \right\rangle\\
&-&\left\langle X,[V,W]   \right\rangle+\left\langle V,[W,X]   \right\rangle+\left\langle W,[X,V]   \right\rangle\\
&-& V\left\langle X,W   \right\rangle+X\left\langle W,V   \right\rangle+W\left\langle V,X   \right\rangle\\
&+&\left\langle V,[X,W]   \right\rangle-\left\langle X,[W,V]   \right\rangle=\left\langle W,[V,X]   \right\rangle\\
&=&2X\left\langle V,W   \right\rangle
\end{eqnarray*}
\QUAN

If we have to emphasize the open set ${\cal U}$ in theorem \ref{levicivita} we shall write $D\left[ {\cal U}   \right]$ instead of $D$:  if $\ {\cal U}_1,{\cal U}_2 \subset \M$ in a maximal atlas
${\cal A}$ for $\M$ are overlapping open sets, then the open set ${\cal U}_{1}\bigcap{\cal U}_{2} $ is in ${\cal A}$   too and 
$$
D\left[{\cal U}_{1}    \right]\vert_{{\cal X}\left({\cal U}_{1}\bigcap{\cal U}_{2}    \right)}=
D\left[{\cal U}_{1}    \right]\vert_{{\cal X}\left({\cal U}_{1}\bigcap{\cal U}_{2}    \right)}
,$$
hence we can collect all local definitions of Levi-Civita connections: 
\begin{definition}\rm
the {\bf Levi-Civita connection} (or {\bf metric connection}) $D$ of $\left(\M,\Lambda    \right)$
is the collection consisting of all the metric connections $\{D\left[{\cal U}_i    \right]\}_{i\in I}$ as ${\cal U}_i $ runs over any maximal atlas ${\cal A}=\left(\{{\cal U}_i\}    \right)_{i\in I}$ on $\M$. 
\label{met_conn}
\end{definition}

So far we have studied nondegenerate holomorphic Riemannian manifolds: this situation is quite similar to real Riemannian geometry.

Things are different, instead, if we allow metrics to have meromorphic behaviour, or even simply to lower somwhere in their ranks.
These metric 'singularities' will be generally supposed to lie in closed hypersurfaces; Levi Civita  connections may still be defined, but, as one could expect, they will turn out to be themselves 'meromorphic'.

Let now $\left(\N,\Lambda    \right)$ be a meromorphic Riemannian manifold admitting closed hypersurfaces ${\cal D}$ and ${\cal E}$ such that 
\begin{itemize}
\item
$\Lambda\vert_{\N\setminus{\cal E}}$ is holomorphic;
\item
$\Lambda\vert_{(\N\setminus{\cal E})\setminus{\cal D}}$ is nondegenerate;
\end{itemize}
Since $\N\setminus{\cal E}$ is connected, we have that $(\N\setminus{\cal E})\setminus{\cal D},\Lambda\vert_{(\N\setminus{\cal E})\setminus{\cal D}}$ is a nondegenerate holomorphic Riemannian manifold admitting, as such, a canonical holomorphic Levi-Civita connection $D$.

Now, if $p\in{\cal D}\bigcup{\cal E}$ and $V,W$ are holomorphic vector fields in a neighbourhood ${\cal V} $ of $p$, it will result that we are able to define the vector field $D_V W$ on ${\cal V}\setminus\left({\cal D}\bigcup{\cal E}    \right)$, and this will be a meromorphic vector field.

Let's state all this more precisely:
\begin{definition}\rm
the {\bf Christoffel symbols } of a coordinate system $$Z=(z^1\cdots z^m)$$ on an open set ${\cal U}\subset\N$ are those complex valued functions, defined on
${\cal U}\setminus\left({\cal D}\bigcup{\cal E}    \right)$ by setting $$\Gamma_{ij}^k =dz^k\left(  D_{\frac{\partial}{\partial z^i}}\left({\frac{\partial}{\partial z^j}}   \right)  \right).$$
\label{christoffel}
\end{definition}
Now the representative matrix $(g_{ij})$ of $\Lambda$ with respect to the coordinate system $Z$ 
is holomorphic in ${\cal U}$, with nonvanishing determinant function on 
${\cal U}\setminus\left({\cal D}\bigcup{\cal E}    \right)$; as such it admits a inverse 
matrix ${g^{ij}}$, whose coefficients hence result in being ${\cal D}\bigcup{\cal E} 
$-meromorphic functions.
\begin{lemma}\rm
\begin{itemize}
\item[(a)]
$D_{\frac{\partial}{\partial z^i}}\left(\sum_{j=1}^m W^j {\frac{\partial}{\partial z^j} }   \right)=
\sum_{k=1}^m\left(\frac{\partial W^k}{\partial z^i}+\sum_{j=1}^m \Gamma_{ij}^k W^j   \right)\frac{\partial}{\partial z^k}$ as meromorphic vector fields;
\item[(b)]
$2\Gamma_{ij}^k =\sum_{m=1}^N g^{km}\left(-g_{ij,m}+g_{im,j}+g_{jm,i}\right)=2\Gamma_{ij}^k $ as meromorphic functions;
\end{itemize}
\label{christoffel_2}
\end{lemma}
{\bf Proof:} at first note that the operation of associating Christoffel symbols 
to a coordinate system is compatible with restrictions, in the sense that the Christoffel symbols of the restriction of $Z$ to a smaller open set are its Christoffel symbols restricted to  that  set:  now, if $$
p\in{\cal U}\bigcap\{n\in\N:\Lambda\hbox{ is holomorphic and nondegenerate at $n$}\}
$$ 
and ${\cal V}_p\subset {\cal U}$ is a neighbourhood of $p$, contained in ${\cal U}$, we have that $\Lambda$ is holomorphic and nondegenerate in ${\cal V}_p$: hence
(a):  by Koszul's formula we have 
$$
2\sum_{a=1}^N \Gamma_{ij}^a g_{am}=2\left\langle D_{\frac{\partial}{\partial z^i}} \frac{\partial}{\partial z^j},\frac{\partial}{\partial z^m}  \right\rangle=
\frac{\partial}{\partial z^i}g_{jm}+\frac{\partial}{\partial z^j}g_{im}+\frac{\partial}{\partial z^m}g_{ij};
$$
multiplying both side by $g^{mk}$ and summing over $m$ yields the desired result; (b) follows immediately from (D3) of definition \ref{met_conn_germ}.
Now the fact that (a) and (b) hold in fact on ${\cal U}$ follows by analytical continuation: note that this result does not depend on the choice of $p$.
\QUAN
\begin{proposition}\rm
For every pair $V,W$ of holomorphic vector field on the open set ${\cal U}$ (belonging to a maximal atlas) in the meromorphic Riemannian manifold $\left(\N,\Lambda     \right)$, $D_V W$ is a well defined vector field, holomorphic on ${\cal U}\bigcap\{n\in\N:\Lambda\hbox{ is holomorphic and nondegenerate at $n$}\}$ and may be extended to a meromorphic vector field on ${\cal U}$.
\end{proposition}
{\bf Proof:} there exist holomorphic functions $\{V^i\}$, $\{W^j\}$ and a coordinate system $Z=\left(z^1.....z^N    \right)$
on ${\cal U}$ such that 
$$
\cases{
\displaystyle{V=\sum_{i=1}^N V^i\frac{\partial}{\partial z^i}}\cr
\displaystyle{W=\sum_{j=1}^N W^i\frac{\partial}{\partial z^j}}
.}
$$

By lemma \ref{christoffel_2}(a),
\begin{eqnarray*}
D_V W=\sum_{i=1}^N V^i D_{\frac{\partial}{\partial z^i}}\left(\sum_{j=1}^N W^j   \frac{\partial}{\partial z^i} \right)=\\
=\sum_{k=1}^N \left(\sum_{i,j=1}^N V^i\left( \frac{\partial W^k}{\partial z^i}+\Gamma^k_{ij}W^j   \right)    \right)\frac{\partial}{\partial z^k}:
\end{eqnarray*}
this is a vector field whose components are meromorphic functions.
\QUAN
Summing up, we yield:
\begin{definition}
\rm
Given a ${\cal D}$-degenerate and ${\cal E}$-meromorphic Riemannian manifold $\left(\N,\Lambda    \right)$, with ${\cal
D} $ and ${\cal E}$ closed hypersurfaces in $\N$, the {\bf Levi-Civita metric connection} (or {\bf meromorphic metric
connection}) of $\N$ is the collection consisting of the metric connections $\{D\left[{\cal U}_i\setminus({\cal D}\bigcup
{\cal E}  ) \right]\}_{i\in I}$ as ${\cal U}\}_i $ runs over any maximal atlas ${\cal B}=\left(\{{\cal U}\}_i    \right)_{i\in I}$ on $\N$. 

\end{definition}
\vfill\eject

\section{Meromorphic parallel translation}
We turn now to study vector fields on paths: an obvious example is the velocity field, defined in definition \ref{vf}: just as in semi-Riemannian geometry, there is a natural way of defining the rate of change $X^{\prime}$ of a meromorphic vector field $X$ on a path.
We study at first paths without boundary, with values in a nondegenerate holomorphic Riemannian manifold $\M$:
let
\begin{itemize}
\item
$Q_{\M}=\left(S,\pi,j,\gamma,\M\right)$ be a path in $\M$, without boundary;
\item
$P$ be the set of branch points of $\pi$;
\item
$r\in S\setminus P$ be a finite-velocity point of $Q_{\M}$;
\item
${\cal V}\subset S\setminus P$ be a neighbourhood of $r$ such that $\gamma\left({\cal V}\right)$ is contained in a local chart in $\M$;
\item
${\cal H\left(V    \right)}$ be the ring of holomorphic functions on ${\cal V}$;
\item
${\cal X}_{\gamma}\left({\cal V}    \right)$ the Lie algebra of holomorphic vector fields over $\gamma$ on ${\cal V}$.
\end{itemize}

\begin{proposition}\rm
There exists a unique mapping $$\nabla_{\gamma^{\prime}}\colon{\cal X}_{\gamma}\left({\cal V}    \right) \longrightarrow {\cal X}_{\gamma}\left({\cal V}    \right) ,$$
called {\bf induced covariant derivative} on $Q_{\M}$ in ${\cal V}$, (or on $\gamma$ in ${\cal V}$) such that:
\begin{eqnarray}
\hbox{(a)}\ & & \nabla_{\gamma^{\prime}}\left(aZ_1+bZ_2    \right)=a\nabla_{\gamma^{\prime}}
Z_1+b\nabla_{\gamma^{\prime}} Z_2,\quad a,b\in\CI\nonumber;\\
\hbox{(b)}\  & & \nabla_{\gamma^{\prime}} \left(hZ\right)=\left(\widetilde{\frac{d}{dz}}h    \right)Z
+h\nabla_{\gamma^{\prime}} Z,\quad h\in{\cal H\left(V    \right)}; \nonumber\\
\hbox{(c)}\ & & \nabla_{\gamma^{\prime}} \left(V\circ\gamma \right)(r)=
D_{\gamma_{*}\vert_r(\widetilde{\frac{d}{dz}}\vert_r)}\ r\in{\cal V},\nonumber
\end{eqnarray}
where $V$
is a holomorphic vector field in a neighbourhood of $\gamma(r)$.
Moreover,
$$\widetilde{\frac{d}{dz}}\left\langle X,Y\right\rangle=\left\langle \nabla_{\gamma^{\prime}}
 X,Y   \right\rangle+\left\langle X,\nabla_{\gamma^{\prime}} Y \right\rangle
\quad X,Y\in {\cal X}_{\gamma}({\cal V}).$$
\label{traspar}
\end{proposition}
{\bf Proof:} let's prove at first uniqueness: suppose existence to hold and let $\left(z^1....z^N    \right)$ be a coordinate system containing $\gamma({\cal V})$; we have, at $\gamma(r)$, 
$$
Z(r)=\sum_{i=1}^N \left(Z\,z^i(r)    \right)\frac{\partial}{\partial z^i}.
$$
Set $Z(z^i)=Z^i$: by (a) and (b),
$$
\nabla_{\gamma^{\prime}}Z(r)=\sum_{i=1}^N \left(\widetilde{\frac{d}{dz}}Z^i   \frac{\partial}{\partial z^i}\vert_{\gamma(r)}+Z^i \nabla_{\gamma^{\prime}}\left( \frac{\partial}{\partial z^i}   \right)\vert_{\gamma(r)}
 \right).
$$
By (c) we deduce
\begin{equation}
\label{diamante}
\nabla_{\gamma^{\prime}}Z=\sum_{i=1}^N  
\left(
\widehat{\frac{d}{dz}}Z^i 
\frac{\partial}{\partial z^i}   
+Z^i D_{\gamma_*}\widetilde{\frac{d}{dz}}
\left( 
\frac{\partial}{\partial z^i}
\right)
\right).
\end{equation}
Thus $\nabla_{\gamma^{\prime}}$ is completely determined by the Levi-Civita connection on $\M$.
Let's prove now existence: 
define $\nabla_{\gamma^{\prime}}$ by (\ref{diamante}): then, using local coordinates $\left(z^1...z^N   \right)$ at $\gamma(r)$ we have:
$$
\nabla_{\gamma^{\prime}}Z=\sum_{k=1}^m\left( \widetilde{\frac{d}{dz}}Z^k+\sum_{i,j=1}^m\Gamma_{ij}^k
\widetilde{\frac{d}{dz}}\left(\gamma^i\right)Z^j   \right)\frac{\partial}{\partial z^k},
$$
hence (a) and (b) are trivial; (c) and (d), instead, require a little bit more of attention.
Within this proof, we shall understand summation with respect to repeated indices, but only with regards to tensor components or Christoffel symbols: we shall instead explicitely write down the summation sign $\Sigma$ when acting on vector fields.

Let's prove (c):
\begin{eqnarray*}
& &\left(V\circ\gamma   \right)(r)=\\
&=&\sum_{k=1}^N
\left(
\widetilde{\frac{d}{dz}}
\left(V^k\circ\gamma 
\right)(r)+
\Gamma_{ij}^k(\gamma(r))
\widetilde{\frac{d}{dz}}
\left(\gamma^i
\right)(r)V^j(\gamma(r))
\right)
\frac{\partial}{\partial z^k }
=\\
&=&
\sum_{k=1}^N 
\left(\frac{\partial V^k}{\partial z^i}
\left(\gamma(r)
\right)
\widetilde{\frac{d}{dz}}\left(\gamma^i\right)(r)+
\Gamma_{ij}^k(\gamma(r))
\widetilde{\frac{d}{dz}}
\left(\gamma^i
\right)(r)V^j(\gamma(r))
\right)\frac{\partial}{\partial z^k }\\
&=&
\sum_{k=1}^N 
\left[
\widetilde{\frac{d}{dz}}
\left(\gamma^i
\right)(r)
\left(
\frac{\partial V^k}{\partial z^i}
\left(\gamma(r)
\right)
+
\Gamma_{ij}^k(\gamma(r))
V^j(\gamma(r))
\right)
\frac{\partial}{\partial z^k }
\right]
\\
&=&
\sum_{k=1}^N 
D_{\sum_{k=1}^N \widetilde{\frac{d}{dz}}
\left(\gamma^i
\right)(r)
 \frac{\partial}{\partial z^i }}
\left(
\frac{\partial V^k}{\partial z^i}
\left(\gamma(r)
\right)
+
\Gamma_{ij}^k(\gamma(r))
V^j(\gamma(r))
\right)\\
&=&
D_{\gamma_{*}\vert_r(\widetilde{\frac{d}{dz}}\vert_r)}.
\end{eqnarray*}
While proving (d), we set $$g_{kl}=\Lambda\left\langle
\frac{\partial}{\partial z^k }
\frac{\partial}{\partial z^l }
\right\rangle\circ\gamma.$$
Now
\begin{eqnarray*}
& & \left\langle \nabla_{\gamma^{\prime}}(X),Y
    \right\rangle +\left\langle X,\nabla_{\gamma^{\prime}}(Y)\right\rangle=\\
&=&
\left\langle \sum_{k=1}^N\left(\widetilde{\frac{d}{dz}}
\left(X^k   \right)+\Gamma_{ij}^k \widetilde{\frac{d}{dz}}
\left(\gamma^i   \right)X^j  \right)\frac{\partial}{\partial z^k },\sum_{l=1}^N
Y^l\frac{\partial}{\partial z^l }
\right\rangle\\
& & \qquad + \left\langle \sum_{l=1}^N
X^l\frac{\partial}{\partial z^l } \sum_{k=1}^N\left(\widetilde{\frac{d}{dz}}
\left(Y^k   \right)+\Gamma_{ij}^k \widetilde{\frac{d}{dz}}
\left(\gamma^i   \right)Y^j  \right)\frac{\partial}{\partial z^k }   \right\rangle\\
&=& g_{kl}\left(\widetilde{\frac{d}{dz}}
\left(X^k   \right)+\Gamma_{ij}^k \widetilde{\frac{d}{dz}}
\left(\gamma^i   \right)X^j  \right)Y^l+ 
g_{kl}\left(\widetilde{\frac{d}{dz}}
\left(Y^k   \right)+\Gamma_{ij}^k \widetilde{\frac{d}{dz}}
\left(\gamma^i   \right)Y^j  \right)X^l\\
&=& g_{kj}\left(\widetilde{\frac{d}{dz}}
\left(X^k   \right)+\Gamma_{il}^k \widetilde{\frac{d}{dz}}
\left(\gamma^i   \right)X^l  \right)Y^j+ 
g_{kl}\left(\widetilde{\frac{d}{dz}}
\left(Y^k   \right)+\Gamma_{ij}^k \widetilde{\frac{d}{dz}}
\left(\gamma^i   \right)Y^j  \right)X^l\\
&=& g_{kl}\widetilde{\frac{d}{dz}}
\left(X^k   \right)Y^l+g_{kl}\widetilde{\frac{d}{dz}}
\left(Y^k   \right)X^l+\widetilde{\frac{d}{dz}}
\left(\gamma^i   \right)X^lY^j\left(g_{kj}\Gamma^k_{il}+g_{kl}\Gamma^k_{ij}\right)\\
&=& g_{kl}\widetilde{\frac{d}{dz}}
\left(X^k   \right)Y^l+g_{kl}\widetilde{\frac{d}{dz}}
\left(Y^k   \right)X^l+\\
& & \qquad +\frac{1}{2} \widetilde{\frac{d}{dz}}\left(\gamma^i   \right)X^lY^j
\left(g_{kj}g^{km}\left(-\frac{\partial g_{il}}{\partial z^m }+\frac{\partial g_{im}}{\partial z^l }+\frac{\partial g_{lm}}{\partial z^i }
   \right)\right.+\\ 
& & \qquad\qquad +
\left. g_{kl}g^{km}\left(-\frac{\partial g_{ij}}{\partial z^m }+\frac{\partial g_{im}}{\partial z^j}+\frac{\partial g_{jm}}{\partial z^i }
   \right)
  \right)\\
&=& g_{kl}\widetilde{\frac{d}{dz}}
\left(X^k   \right)Y^l+g_{kl}\widetilde{\frac{d}{dz}}
\left(Y^k   \right)X^l+\\
& & \qquad + \frac{1}{2}\widetilde{\frac{d}{dz}}\left(\gamma^i   \right)X^lY^j
\left(\delta_j^m\left(-\frac{\partial g_{il}}{\partial z^m }+\frac{\partial g_{im}}{\partial z^l }+\frac{\partial g_{lm}}{\partial z^i }
   \right)\right.+\\
& & \qquad\qquad + 
\left.\delta_l^m\left(-\frac{\partial g_{ij}}{\partial z^m }+\frac{\partial g_{im}}{\partial z^j}+\frac{\partial g_{jm}}{\partial z^i }
   \right)
  \right)\\
&=& g_{kl}\widetilde{\frac{d}{dz}}
\left(X^k   \right)Y^l+g_{kl}\widetilde{\frac{d}{dz}}
\left(Y^k   \right)X^l+\\
& & \qquad + \frac{1}{2}\widetilde{\frac{d}{dz}}\left(\gamma^i   \right)X^lY^j
\left(-\frac{\partial g_{il}}{\partial z^j}+\frac{\partial g_{ij}}{\partial z^l }+\frac{\partial g_{lj}}{\partial z^i }
-\frac{\partial g_{ij}}{\partial z^l }+\frac{\partial g_{il}}{\partial z^j}+\frac{\partial g_{jl}}{\partial z^i }
  \right)\\
&=& g_{kl}\widetilde{\frac{d}{dz}}
\left(X^k   \right)Y^l+g_{kl}\widetilde{\frac{d}{dz}}
\left(Y^k   \right)X^l+\widetilde{\frac{d}{dz}}\left(\gamma^i   \right)X^lY^j \frac{\partial g_{lj}}{\partial z^i }\\
&=& g_{kl}\widetilde{\frac{d}{dz}}
\left(X^k   \right)Y^l+g_{kl}\widetilde{\frac{d}{dz}}
\left(Y^k   \right)X^l+\left(\widetilde{\frac{d}{dz}}
\left(g_{kl}   \right)   \right)\\
&=& \widetilde{\frac{d}{dz}}
\left\langle X,Y  \right\rangle
\end{eqnarray*}
\QUAN
Now let ${\cal R}=\{{\cal V}_k\}_{k\in K}$ be a maximal atlas for $S\setminus P$; we may assume that, for every $k$,
maybe shrinking ${\cal V}_k   $, $\gamma\left( {\cal V}_k   \right)$ is contained in some local chart ${\cal U}_i$ in the
already introduced atlas ${\cal A}$ for $\M$.

By proposition \ref{traspar}, if ${\cal V}_1  $ and ${\cal V}_2  $ are overlapping open sets in ${\cal R}$, ${\cal V}_1 \bigcap{\cal V}_2\in {\cal R} $ too, and 
$$
\nabla_{\gamma^{\prime}}
\left[ {\cal V}_1   \right]
\vert_{{\cal V}_1 \bigcap{\cal V}_2}
=\nabla_{\gamma^{\prime}}
\left[ {\cal V}_2   \right]
\vert_{{\cal V}_1 \bigcap{\cal V}_2}.
$$

Now let's complete ${\cal R}$ to an atlas ${\cal S}$ for $S$: keeping into account that the local coordinate expression of the induced covariant derivative is
$$
\nabla_{\gamma^{\prime}}Z=\sum_{k=1}^m\left( \widetilde{\frac{d}{dz}}Z^k+\sum_{i,j=1}^m\Gamma_{ij}^k
\widetilde{\frac{d}{dz}}\left(u^i\circ\gamma\right)Z^j   \right)\frac{\partial}{\partial u^k}.
$$
and arguing as about the meromorphic Levi-Civita connection, we are able to show that pairs of holomorphic vector fields on $\gamma$ are transormed into $P$-meromorphic vector fields on $\gamma$.

\begin{definition}\rm
The $P$-meromorphic {\bf induced covariant derivative}, or the $P$-meromorphic parallel translation on a path $Q_{\M}=\left(S,\pi,j,\gamma,\M    \right)$ with set of branch points $P$ and taking values in a nondegenerate Riemannian manifold $\M$ is the collection consisting of the induced covariant derivatives $\nabla_{\gamma^{\prime}}\left[ {\cal V}_k \setminus P  \right] $ as ${\cal V}_k $ runs over a maximal atlas ${\cal S}=\left(\{{\cal V}_k\}    \right)_{k\in K}$ on $S$. 
\label{indcovdev}
\end{definition}

Let's turn now to dealing with meromorphic parallel translations induced on a path $Q_{\N}=\left(T,\varrho,j,\delta,\N    \right)$, without boundary, in a meromorphic Riemannian manifold
$(\N,\Lambda)$ admitting closed hypersurfaces ${\cal D}$ and ${\cal E}$ such that 
\begin{itemize}
\item
$\Lambda\vert_{\N\setminus{\cal E}}$ is holomorphic;
\item
$\Lambda\vert_{(\N\setminus{\cal E})\setminus{\cal D}}$ is nondegenerate.
\end{itemize}
We set ${\cal F}={\cal D}\bigcup{\cal E}$ and restrict our attention to paths $z_0$-starting at metrically ordinary points, supposing, without loss in generality, that $z_0=0$.
\begin{lemma}\rm
Set $\M=\N\setminus{\cal F}$, $S=\delta^{-1}(\M)$: then $T\setminus S$ is discrete, hence $S$ is a connected Riemann surface.
\end{lemma}
{\bf Proof:} suppose that there exists a subset ${\cal V}\subset T\setminus S$ admitting an accumulation point $t\in{\cal V}$ and consider a countable atlas for ${\cal B}=\{U_n\}_{n\in\ENNE}$ for $\N$ such that, for every $n$, there exists $\Psi_n\in{\cal O}\left( \{U_n\}   \right)$ such that 
$$
U_n\bigcap{\cal F}=\{X\in U_n : \Psi_n=0\}
$$
(If $U_n\bigcap{\cal F}=\emptyset
$,  pick a never vanishing holomorphic function). 

Set $\delta^{-1}(U_n)=T_n\subset T$ and suppose, without loss of generality, that $\delta(t)\in U_0$.

We have $\Psi_0\circ\delta\vert_{{\cal V}\cap T_0}=0$ and $t\in{\cal V}\cap T_0$ is an accumulation point of ${\cal V}\cap T_0$ , hence $\Psi_0\circ\delta\vert_{ T_0}=0$ and $\delta(T_0)\subset {\cal F}$.

Suppose now that $T_N\not=\emptyset$ for some $N$: we claim that this implies $\delta(T_N)\subset{\cal F}$: to prove the asserted, pick two points $\tau_0\in T_0$ and $\tau_n\in T_n$ and two neighbourhoods $T^{\prime}_0$, $T^{\prime}_N$ of $\tau_0 $ and $\tau_n$ in $T_0$ and $T_n$ respectively, such that $\varrho\vert_{T^{\prime}_0}$ and $\varrho\vert_{T^{\prime}_N}$ are biholomorphic functions.
Now the function elements $\left(\varrho(T^{\prime}_0),\delta\circ\left( \varrho\vert_{T^{\prime}_0}   \right)^{-1}    \right)$ and $\left(\varrho(T^{\prime}_N),\delta\circ\left( \varrho\vert_{T^{\prime}_N}   \right)^{-1}    \right)$ are connectible, hence there exists a finite chain $\{W_{\nu}\}_{\nu=0...L}$ such that $W_0=\varrho(T^{\prime}_0)$, $W_L=\varrho(T^{\prime}_N)$, $W_{\nu}\bigcap W_{\nu+1}\not=0$ for every $\nu$.

Without loss of generality, we may suppose that each $W_{\nu}$ admits a holomorphic, hence open, immersion $j_{\nu}\longrightarrow T$, hence, setting
$$
\cases{
S_0=T_0\cr
S_{\lambda}=j_{\lambda}(W_{\lambda}) \hbox{\rm \ for }
\lambda=1...L\cr
S_{L+1}=T_N,
}
$$
yields a finite chain of open subsets $\{S_{\lambda}\}_{\lambda=0...M}$ of $T$ connecting $T_0$ and $T_N$.

Let's prove, by induction, that, for every $\lambda$, $\delta(S_{\lambda})\subset{\cal F}$.

$\bullet$ {\bf At first} recall that
$\delta(S_0)\subset U_0\bigcap{\cal F}$ as already proved; suppose now that $\delta(S_{k-1})\subset {\cal F}$. 
We have $S_{k-1}\bigcap S_k\not=\emptyset$, hence $\delta(S_{k-1})\bigcap\delta(S_{k})\not=\emptyset$.

For every $m$ set $$\Sigma_{km}=\delta(S_{k-1})\bigcap\delta(S_{k})\bigcap U_m:$$ if $\Sigma_{km}\not=\emptyset$, then $$\Psi_m\circ\delta\vert_{\delta^{-1}(\Sigma_{km})\bigcap S_{k-1}\bigcap S_k }\equiv 0;$$ but $\delta^{-1}(\Sigma_{km})\bigcap S_{k-1}\bigcap S_k$ is open in $\delta^{-1}\left(\delta(S_k)\bigcap U_m    \right)\bigcap S_k$, thus
$$\Psi_m\circ\delta\vert_{\delta^{-1}\left(\delta(S_k)\bigcap U_m    \right)\bigcap S_k}\equiv 0,$$ that is to say $\delta(S_k)\bigcap U_m\subset{\cal F}$.

$\bullet$ {\bf On the other hand},
if $\Sigma_{km}=\emptyset$, but $\delta(S_k)\bigcup U_m\not=\emptyset$ we claim that $\delta(S_k)\bigcap U_m\subset{\cal F}$ as well: proving this requires a further induction: pick a $U_M$ such that $\Sigma_{kM}\not=\emptyset$ and a finite chain of open sets
${\cal B}^{\prime}=\{U^{\prime}_{\mu}\}_{\mu=0...J}\subset {\cal B}$ (with $U_{\mu}^{\prime}\bigcap\delta(S_k)\not=\emptyset$ for each $\mu$) connecting $U_M$ and $U_m$.
Since $\Sigma_{kM}\not=\emptyset$, 
$$\delta(S_k)\bigcap U_0^{\prime}=\delta(S_k)\bigcap U_M\subset {\cal F};$$ suppose by induction that 
$$\delta(S_k)\bigcap U_{l-1}^{\prime}\subset {\cal F} :$$
then
$$\Psi_l\circ\delta\vert_{\delta^{-1}\left(\delta(S_k)\cap U_{l-1}^{\prime}\cap U_{l}^{\prime}   \right)\cap S_k}\equiv 0,$$
hence
$$
\Psi_l\circ\delta\vert_{\delta^{-1}\left(\delta(S_k)\cap U_{l}^{\prime}   \right)\cap S_k}\equiv 0,$$
i.e. $\delta(S_k)\bigcap U_{l}^{\prime}  \subset {\cal F}$: this ends the induction and eventually implies 
$$\delta(S_k)\bigcap U_{m}= \delta(S_k)\bigcap U_{J}^{\prime}\subset {\cal F}.$$

Summing up, $$\delta(S_k)=\bigcup_m \left( \delta(S_k)\bigcap U_m   \right)\subset {\cal F},$$ for each $k$;  hence $$\delta(T_N)=\delta(S_M)\subset{\cal F}$$  and eventually $$\delta(T)=\delta\left(\bigcup_{N\in\ENNE} T_N    \right)\subset {\cal F},$$ hence $\delta$ couldn't start at a point in $\N\setminus{\cal F}$.
\QUAN

In the following considerations, there will still hold all notations introduced in
preceding lemma: given a path $Q_{\N}^{\sharp}=\left(T\bigcup T^{\flat},\varrho^{\sharp}
,j,\delta,\N    \right)$, with possibly nonempty boundary $T^{\flat}$ set $Q_{\N}=\left(T,\varrho,j,\delta,\N    \right)$, set $\pi=\varrho\vert_S$, where $\varrho=\varrho^{\sharp}\vert_{T}$, we call 
$Q_{\N}$ the interior path of $Q_{\N}^{\sharp}$.
Set now $\pi=\varrho\vert_S$,
$\gamma=\delta\vert_S$ and note that, since $Q_{\N}$ is starting from a metrically ordinary point $m$, $j$ may be
supposed to take values in fact in $S$; since the preceding lemma shows that $S$ is a connected Riemann surface,
$Q_{\M}=\left(S,\pi,j,\delta\vert_S\M    \right)$ is in fact a path in $\M$, which we call the {\bf depolarization} of $Q_{\N}$. But $\M$ is a nondegenerate holomorphic Riemannian manifold, hence if $P$ is the set of branch points of $\pi$, there is a $P$-meromorphic induced parallel translation on $Q_{\M}$, got following definition \ref{indcovdev} and its substratum.
Finally, we introduce a maximal atlas ${\cal T}$ for $T$ and yield the following:
\begin{definition}\rm
Let $\left(\N,\Lambda    \right)$ be a ${\cal E}$- meromorphic and ${\cal D}$-degenerate Riemannian manifold,
$\M=\N\setminus\left({\cal D}\bigcup{\cal E}    \right)$,
$Q_{\N}^{\sharp}=\left(T\bigcup T^{\flat},\varrho^{\sharp},j,\delta,\N    
\right)$ a path  with possibly nonempty boundary $T^{\flat}$
and $Q_{\N}$ its interior path: 
the {\bf $\left(P\bigcup\delta^{-1}\left( {\cal D}\bigcup{\cal E}   \right)    \right)$-meromorphic induced covariant derivative} on $Q_{\N}$
 is the collection consisting of all induced covariant derivatives $\nabla_{\gamma^{\prime}}\left[ {\cal V}_k \bigcap S
 \right] $ as ${\cal V}_k$ runs over a maximal atlas ${\cal T}=\left(\{{\cal V}_k\}    \right)_{k\in K}$ for  $T$ and
 $Q_{\M}=\left(S,\pi,j,\delta\vert_S\M    \right)$ is the depolarization of $Q_{\N}$.  
\label{indcovdev2}
\end{definition}
\vfill\eject
\section{Geodesics}
\begin{definition}\rm
A meromorphic (in particular, holomorphic) vector field $Z$ on a path without boundary $Q_{\M}=\left(S,\pi,j,\gamma,\M\right)$
is {\bf parallel} provided that $\nabla Z=0$ (as a {\bf meromorphic} field on $Q_{\M}$).
\label{parall}
\end{definition}
\begin{definition}\rm
The {\bf acceleration} $\aleph\left( Q_{\M}\right)$ of $Q_{\M}$ is the meromorphic field $\nabla \left(V\left(Q_{\M}\right)\right)$ on $Q_{\M}$ yielded by the induced covariant
derivative of its velocity field; 
\label{acceleration}
the {\bf speed} of a path is the 'amplitude' function of its velocity field:
 $S\left(Q_{\M}    \right)(r)=
 \left\langle \gamma_{*}\vert_r\left(\widetilde{\frac{d}{dz}}\right), \gamma_{*}\vert_r\left(\widetilde{\frac{d}{dz}}   \right)   \right\rangle$. This is a meromorphic function.
A path is {\bf null} provided that its speed is zero everywhere. \label{speed}
\end{definition}
\begin{definition}\rm
A {\bf geodesic} in a meromorphic (in particular, holomorphic) Riemannian manifold is a path whose interior path's velocity field is parallel, or, equivalently, one of zero acceleration (see definition \ref{acceleration}).
A geodesic is {\bf null} provided that so is as a path.
\label{geodesic}
\end{definition}
The local equations of elements of geodesics $\left(U,\beta   \right)$

$$
\cases{ \BETA^{\bullet\bullet}{}^{k}+\sum_{i,j=1}^N \Gamma_{ij}^k(\beta)\BETA^{\bullet}{}^{i}\BETA^{\bullet}{}^{j}=0\cr
 \ \cr
k=1.....N\cr 
}
$$
are a system of $N$  second-order ordinary differential equations in the complex domain, with meromorphic coefficients, in turn equivalent to an autonomous system of $2N$ first-order equations, hence, as a consequence of the general theory (see theorem \ref{gerexuq}) we have the following
\begin{theorem}
\TTT
For every metrically ordinary point $p\in\M$, every holomorphic tangent vector $V_p\in T_p\M$ and every $z_0\in\CI$, there exists a unique germ 
\hbox{\boldmath{}$\beta_{z_0}$\unboldmath}
of geodesic such that 
$$
\cases{
\hbox{\boldmath{}$\beta_{z_0}$\unboldmath $(z_0)=p$} \cr
\hbox{\boldmath{}$\beta_{z_0\ *}$\unboldmath$(d/dz)\vert_{z_0}=V_p$;}
}
$$
moreover any analytical continuation of \hbox{\boldmath{}$\beta_{z_0}$\unboldmath}
is a geodesic (see theorem \ref{ancont}).
\end{theorem}
\QUAN
\begin{definition}\rm
Let now $Q_{\N}=\left(S,\pi,j,\gamma,\N    \right)$ be a nonconstant path: another path $P_{\N}=\left(R,\varrho,i,\delta,\N    \right)$ is a {\bf reparametrization} of $Q_{\N}$, or, shortly, $\delta$ is a reparametrization of $\gamma$ if there exists a holomorphic function $\phi:R\longrightarrow S$ such that $\delta=\gamma\circ\phi$.
\end{definition}
$$
\bullet (\delta=\gamma\circ\phi)
$$
$$
\ 
$$
$$
\def\normalbaselines{\baselineskip20pt\lineskip3pt \lineskiplimit3pt}
\def\mapright#1{\smash{\mathop{\longrightarrow}\limits^{#1}}}
\def\mapdown#1 {\Big\downarrow\rlap{$\vcenter{\hbox{$\scriptstyle #1$}}$}}
\matrix{
R   &\mapright{\phi}    &S         &\mapright{\gamma}  &\N\cr
\mapdown{\rho}    &     &\mapdown{\pi}       &                      &  \cr
\CI_{\tau}          &     &\CI_t     &                       &  \cr
}
$$

\begin{proposition}\rm
Suppose that $Q_{\N}$ is a geodesic: then $P_{\N}$ is a geodesic if and only if there exist two complex contstants $a,b$, with $a\not=0$ and, for every regular point $r$ of $\varrho$ and $s$ of $\pi$ such that $s=\phi(r)$, there exist two neighbourhoods $U_r$
of $r$ and $V_s$ of $s$ such that $\varrho\vert_{U_r}$ and $\pi\vert_{V_s}$ are biholomorphic and $\pi\circ\phi\circ\left(\varrho\vert_{U_r}    \right)^{-1}(z)=a\,z+b$.
\end{proposition}
{\bf Proof:} the 'if' implication is trivial; on the other hand, suppose that $\gamma\circ\phi $ is a geodesic and $r$ is a regular point for the branched covering $\varrho$: then there exist neighbourhoods $U_r$
of $r$ and $V_s$ of $s=\phi(r)$ such that $\varrho\vert_{U_r}$ and $\pi\vert_{V_s}$ are biholomorphic and 
$\gamma\circ\phi\circ\varrho\vert_{U_r}^{-1}$
is a geodesic.
 
Pick now a chart $\left({\cal W},\Theta    \right)$ around $\gamma(s)$ and set 
$$
\cases
{\gamma^{k}=\Theta^k\circ\gamma\circ\pi\vert_{V_s}^{-1}\cr 
\psi=\pi\circ\phi\circ\varrho\vert_{U_r}^{-1}; 
}
$$
denote derivation with respect 
to $\tau$ ba an apex and with respect to $t$ by a dot;
then, for each $k$:
\begin{eqnarray*}
0 &=& \left(\gamma^k\left(\psi(\tau)    \right)    \right)^{\prime\prime}+\Gamma_{ij}^k 
            \left(\gamma^i\left( \psi(\tau)     \right)    \right)^{\prime}
            \left(\gamma^j\left( \psi(\tau)     \right)    \right)^{\prime}\\
&=& \GAMMA^{\bullet\bullet}{}^{k}\left(\psi(\tau)    \right) \left(\psi^{\prime} (\tau)   \right)^2+\GAMMA^{\bullet}{}^{k}\left(\psi(\tau)    \right) \psi^{\prime\prime}(\tau)+\\
&\ &\qquad\Gamma_{ij}^k\left(\gamma\circ\psi(\tau)    \right)\GAMMA^{\bullet}{}^{i}\left(\psi(\tau)    \right) 
\GAMMA^{\bullet}{}^{j}\left(\psi(\tau)    \right)\left(\psi^{\prime(\tau)}    \right)^2\\
&=& \GAMMA^{\bullet}{}^{k}\left(\psi(\tau)    \right) \psi^{\prime\prime}(\tau).
\end{eqnarray*}

Since $\gamma$ is a nonconstant path, then, for some $k$, $\GAMMA^{\bullet}{}^{k}\not=0$, hence it should be $\psi^{\prime\prime}=0$, i.e $\pi\circ\phi\circ\varrho\vert_{U_r}^{-1}(z)=a\,z+b$, for some complex constants $a,b$, with $a\not=0$.
The uniformity of $a$ and $b$ with respect to $r$ follows now from the fact that, if we pick another regular point $r^{\clubsuit}$ of $\varrho$, the holomorphic function element $\pi\circ\phi\circ\varrho\vert_{U_r^{\clubsuit}}^{-1}$ is connectible with $\pi\circ\phi\circ\varrho\vert_{U_r}^{-1}$.
\QUAN
\begin{definition}\rm
A {\bf pregeodesic} is a path admitting a reparametrization as a geodesic.
\label{pregeo}
\end{definition}
\subsection{A local reparametrization theorem}
Somewhat surprisingly, the natural parameter of any nonconstant geodesic could be cut off from its local equations: in other terms, we are able to express, at least locally, the behaviour of such geodesic in terms of a coordinate function.
Of course this is not in general  a geodesic parametrization, but this may be although restored by imposing a constant-speed like condition.
We state all this more precisely:
\begin{theorem}
\TTT
For  every nonconstant geodesic function element $\quad$ $\left(z_0,{\cal U},\beta    \right)$, 
$z_0$-starting at the metrically ordinary point $m\in\N$ in the meromorphic Riemannian 
manifold $\N$ and every local chart $$\left(W,\Theta    \right)=\left(W,\Theta^1...\Theta^N    \right)$$ around $m$ (with $\Theta:W\longrightarrow\CI^N_{v^1...v^N}$ and $\Theta(m)=0$), there exist:
\begin{itemize}
\item
an integer $l$, $1\leq l\leq N$, which may be supposed to equal $N$;
\item
a (one-dimensional) neighbourhood ${\cal V}$ of $0$ in $\CI_{v^N}$;
\end{itemize}
such that:
\begin{itemize}
\item
$\beta^N\left(=\Theta^N\circ\beta    \right):{\cal U}\longrightarrow{\cal V}$ is biholomorphic, with inverse function $\eta:{\cal V}\longrightarrow{\cal U}$;
\item
$\beta\circ\eta$ is pregeodesic (see definition \ref{pregeo}).
\item
$\left(\beta^1\circ\eta...\beta^{N-1}\circ\eta    \right)$
satisfies the following $N-1$ dimensional second order Cauchy's problem, an apex denoting derivation with respect to $v_N$ and $$\left(\gamma^1...\gamma^{N-1}    \right)=\left(\beta^1\circ\eta...\beta^N\circ\eta    \right) :$$
\end{itemize}
$$
\cases{
\displaystyle{
\gamma^{k\ \prime\prime}-\gamma^{k\ \prime}
\left(
\sum_{i,j=1}^{N-1}\Gamma_{ij}^N (\gamma,v^N) 
\gamma^{i\ \prime}\gamma^{j\ \prime} +2 \Gamma_{iN}^N (\gamma,v^N)        \gamma^{i\ \prime} +\Gamma_{NN}^N(\gamma,v^N)   
 \right) }
\cr
\displaystyle{
+\sum_{i,j=1}^{N-1}\Gamma_{ij}^k (\gamma,v^N) \gamma^{i\ \prime}\gamma^{j\ \prime} +2 \Gamma_{iN}^k (\gamma,v^N) \gamma^{i\ \prime} +\Gamma_{NN}^k(\gamma,v^N) =0
}
\cr 
\ 
\cr
\displaystyle{
\gamma(0)=0
}
\cr
\displaystyle{
\gamma^{k\ \prime}(0)=\left(\BETA^{\bullet}/\VI^{\bullet}{}^{N} \right) (z_0)
}
\cr
\displaystyle{
k=1...N-1
}
\cr
\
}
$$   
\label{repam}
\end{theorem}
{\bf Proof:} since $\beta$ a nonconstant geodesic, one at least among the $\BETA^{\bullet}{}^{J}(v_0)$ is nonzero: maybe reordering the $N$ coordinate functions, we may suppose, without loss of generality, that $\BETA^{\bullet}{}^{N}(v_0)\not=0$: then $\beta^N$ is biholomorphic at $v_0$ and, as such, admits a holomorphic inverse function $\eta:{\cal V}\longrightarrow{\cal U}$ for suitable neighbourhoods ${\cal V}$ of $0$ and ${\cal U}$ of $z_0$.
The fact that $\beta\circ\eta$ is pregeodesic simply follows from being $\eta$  biholomorphic.
Now, if $k\not=N$ we have, by the chain rule $\gamma^{ k\ \prime}=\BETA^{\bullet}{}^{k}\eta^{\prime}$; set $w=v^N\circ\beta$: then $\BETA^{\bullet}{}^{k}=\gamma^{ k\ \prime}$ and $\BETA^{\bullet\bullet}{}^{k}=\gamma^{ k\ \prime\prime}\WI^{\bullet}{}^{2}+\gamma^{ k\ \prime}\WI^{\bullet\bullet}$.
On the other hand, 
$$
\spadesuit\quad
\cases
{\BETA^{\bullet\bullet}{}^{k}+\sum_{i,j=1}^N \Gamma_{ij}^k(\beta)\BETA^{\bullet}{}^{i}\BETA^{\bullet}
{}^{j} 
\cr
k=1.....N,
}
$$
in particular 
\begin{eqnarray*}
\WI^{\bullet\bullet}
&=&-\sum_{i,j=1}^N \Gamma_{ij}^k(\beta)\BETA^{\bullet}{}^{i}\BETA^{\bullet}{}^{j}\\
&=&-\left(\sum_{i,j=1}^{N-1}\Gamma_{ij}^N (\gamma,v)\gamma^{i\ \prime}\gamma^{j\ \prime} +2 \Gamma_{iN}^N (\gamma,v) \gamma^{i\ \prime} +\Gamma_{NN}^N(\gamma,v)    \right)\WI^{\bullet}{}^2, 
\end{eqnarray*}
hence, for $k\not=N$,
$$
\BETA^{\bullet\bullet}{}^{k}=\WI^{\bullet}{}^2
\left(
\gamma^{k\ \prime\prime}-\gamma^{k\ \prime}
\left(
\sum_{i,j=1}^{N-1}\Gamma_{ij}^N (\gamma,v)  \gamma^{i\ \prime} \gamma^{j\ \prime} +2 \Gamma_{iN}^N (\gamma,v)\gamma^{i\ \prime} +\Gamma_{NN}^N(\gamma,v)    \right)    \right);
$$
on the other hand, 
$$
\sum_{i,j=1}^N \Gamma_{ij}^k(\beta)\BETA^{\bullet}{}^{i}\BETA^{\bullet}{}^{j}
=\WI^{\bullet}{}^2
\left(\sum_{i,j=1}^{N-1}\Gamma_{ij}^N (\gamma,v)\gamma^{i\ \prime}\gamma^{j\ \prime} +2 \Gamma_{iN}^N (\gamma,v)   \gamma^{i\ \prime} +\Gamma_{NN}^N(\gamma,v)    \right),
$$
hence substituting in $\spadesuit$ ends the proof.
\QUAN
\subsection{Geodesic completeness}
We eventually yield a definition of geodesic completeness, which will be an extension of the classical one:
\begin{definition}
\rm
A meromorphic Riemannian manifold is {\TTT geodesically complete} provided that the Riemann surface, with logarithmic singularities, of each geodesic
germ starting at a metrically ordinary point is complete as a²path with boundary.
\end{definition}
\vfill\eject
\chapter{Completeness theorems}
\section{Warped products}

In this section we shall be concerned with warped products of Riemann surfaces, each one endowed with some meromorphic metric: in this framework we shall prove a geodesic completeness criterion.

Let now $\,{\cal U}_i$,
$\left(i=1....N   \right)$, $N\geq 2$ be either a copy of the unit ball in the complex plane, or the complex plane itself, whose coordinate function we shall call $u^i$.

Moreover, let each $\,{\cal U}_i\,$ be endowed with a (not everywhere vanishing) meromorphic metric, which we denote by $b_1(u^1)\,du^1\odot du^1$ on $\,{\cal U}_1\,$, or by $f_i(u^i)\,du^i\odot du^i$ if $i\geq 2$, where
both $b_1$ and the $f_i$'s are nonzero meromorphic functions.

Consider now the meromorphic Riemannian manifold

$$
{\cal U}={\cal U}_1\times_{a_2(u^1)}
{\cal U}_2\times_{a_3(u^1)}
{\cal U}_3\times
........
\times_{a_N(u^1)}{\cal U}_N,
$$
where the $a_k$'s ($k\geq 2$) are
nonzero meromorphic warping functions defined on ${\cal U}_1$, i.e. depending solely on $u^1$.

We could write down the meromorphic metric $\Lambda$ of ${\cal U}$ in the form 

$$
\Lambda\left(u^1.....u^N\right)=
b_1(u^1)\,du^i\odot du^i+
\sum_{i=2}^N a_i(u^i)f_i(u^i)
\,du^i\odot du^i.
$$

In other words, the representative matrix of $\Lambda $, with respect to the canonical coordinates of ${\cal U}$, inherited from $\CI^N$, is of the form
$$
(g_{ik})=\left(
\matrix
{
b_1(u^1) &\       &\       &\       &\       &\       \cr
\ &a_2(u^1)f_2(u^2)&\     &\       &\       &\       \cr
\ &\     &a_3(u^1)f_3(u^3)&\       &\       &\       \cr
\ &\     &\             &\       &\       &\ddots &\       \cr
\ &\   &\      &\   &\    &\  & a_N(u^1)f_N(u^N)\cr
}
\right),
$$
where the blanks should be filled in with zeroes.
\begin{lemma}\rm
\label{connessione}
The meromorphic Levi-Civita connection induced on ${\cal U}$ by ${\Lambda}$ is characterized by admitting the following Christoffel symbols $\Gamma_{ij}^k$:
$$
\hbox{\TTT  if\ } k=1
\cases{
\displaystyle
2\Gamma_{11}^1=\frac{b_1^{\prime}(u^1)}
{b_1(u^1)};
\cr
\displaystyle
\Gamma_{ij}^1=0\ \hbox{if\  } i\not=j;\cr
\displaystyle
2\Gamma_{ii}^1=-\frac{a_i^{\prime}(u^1)f_i(u^i)}
{b_1(u^1)}\ \ \hbox{if\  }1\leq i\leq N;\cr
}
$$
$$
\hbox{\TTT \ if\ }2\leq k\leq N\,
\cases{
\displaystyle
2\Gamma_{kk}^k=\frac{f_k^{\prime}(u^k)}
{f(u^k)} \cr
\displaystyle
2\Gamma_{ik}^k =\frac{a_k^{\prime}(u^1)}
{a_k(u^1)}\ \hbox{if } i=1\cr
\displaystyle
\Gamma_{ij}^k=0\ \hbox{otherwise.}\cr
}
$$
\end{lemma}
{\bf Proof:} the inverse matrix of the representative one of $\Lambda$ is
$$
(g^{ik})=\left(
\matrix
{\displaystyle
\frac{1}{b_1(u^1)} &\       &\       &\       &\       &\       \cr\displaystyle
\ &\frac{\displaystyle
1}{\displaystyle
a_2(u^1)f_2(u^3)}&\     &\       &\       &\       \cr\displaystyle
\ &\     &\frac{\displaystyle
1}{\displaystyle
a_3(u^1)f_3(u^3)}&\       &\       &\       \cr\displaystyle
\ &\     &\             &\       &\       &\ddots &\       
\cr\displaystyle
\ &\   &\      &\   &\    &\  &\frac{\displaystyle
1}{\displaystyle
b_1(u^N)f_N(u^N)}\cr
}
\right),
$$
where the blanks should be filled in with zeroes;
therefore there holds 
$$
2\Gamma_{ij}^k=g^{kk}\left(
-\frac{\partial g_{ij}}{\partial u^k}
+\frac{\partial g_{ik}}{\partial u^j}
+\frac{\partial g_{jk}}{\partial u^i}
\right),
$$
where $k$ is meant to be fixed.

Now, if $k=1$, then
$$
2\Gamma_{ij}^1=g^{kk}\left(
-\frac{\partial g_{ij}}{\partial u^1}
+\frac{\partial g_{i1}}{\partial u^j}
+\frac{\partial g_{j1}}{\partial u^i}
\right);
$$
but
\begin{itemize}
\item $\displaystyle\frac{\partial g_{ij}}{\partial u^1}
\not=0$ only if $i=j$;
\item
$\displaystyle\frac{\partial g_{i1}}{\partial u^j}\not=0$ only if
$i=j=1$;
\item
$\displaystyle\frac{\partial g_{j1}}{\partial u^i}\not=0$  only if
$i=j=1$, 
\end{itemize}
hence 
$$
2\Gamma_{11}^1=g^{11}\frac{\partial g_{11}}{\partial u^1}=\frac{b_1^{\prime}(u^1)}
{b_1(u^1)};
$$
on the other hand, if $2\leq i\leq N$,
$$
2\Gamma_{ii}^1=
-g^{11}\frac{\partial g_{ii}}{\partial u^1}=-\frac{a_i^{\prime}(u^1)f_i(u^i)}
{b_1(u^1)}.
$$

Instead, if $2\leq k\leq N$,
\begin{itemize}
\item $\displaystyle\frac{\partial g_{ij}}{\partial u^k}
\not=0$ only if $i=j=k$;
\item
$\displaystyle\frac{\partial g_{ik}}{\partial u^j}\not=0$ only if
$i=j=k$ or if $i=k\not=j$ and $j=1$;
\item
$\displaystyle\frac{\partial g_{jk}}{\partial u^i}\not=0$  only if
$i=j=k$ or if $j=k\not=i$ and $i=1$, 
\end{itemize}
hence 
$$
2\Gamma_{kk}^k=g^{kk}\frac{\partial g_{kk}}{\partial u^k}=\frac{f_k^{\prime}(u^k)}
{f_k(u^k)};
$$
moreover
$$
2\Gamma_{ik}^k=
2\Gamma_{ki}^k=
g^{kk}\frac{\partial g_{kk}}{\partial u^1}=
\frac{a_k^{\prime}(u^1)}{a_k(u^1)}\,\delta_{k1}:
$$
this ends the proof.
\QUAN
\begin{lemma}\label{equazionisecondoordine}
\rm
Each element of geodesic of $\left({\cal U},\Lambda   \right)$ satisfies the following system of $N$ ordinary differential equations in the complex domain: 
\end{lemma}
\begin{equation}
\label{equazioninormali}
\cases{
\displaystyle
\U^{\bullet\bullet}{}^1(z)+\frac{b_1^{\prime}(u^1(z))}{2b_1(u^1(z))}\left(\U^{\bullet}{}^1(z)   \right)^2\cr
\displaystyle
\qquad\quad -\sum_{l=2}^N \frac{a_l^{\prime}(u^1(z))f_l(u^l(z))}{2b_1(u^1(z))}\left(\U^{\bullet}{}^l (z)  \right)^2=0\cr
\displaystyle
\ \cr
\displaystyle
\U^{\bullet\bullet}{}^k(z)+
\frac{f_k^{\prime}(u^k(z))}{2f_k(u^k(z))}\left(\U^{\bullet}{}^k(z)\right)^2\cr
\displaystyle
\qquad\quad
+\frac{a_k^{\prime}(u^1(z))}{a_k(u^1(z))}\left(\U^{\bullet}{}^1 (z)  \right)\left(\U^{\bullet}{}^k (z)  \right)
=0,\ k=2...N,
}
\end{equation}
provided that it starts at a metrically ordinary point.

{\bf Proof:} this is an immediate consequence of lemma \ref{connessione}.
\QUAN
\begin{lemma}\rm
\label{integraleprimo1}
The system of equations (\ref{equazioninormali}) of every element of geodesic 
$$\displaystyle
z\longmapsto\left(u^1(z)...u^N(z)    \right)
$$
of $\left({\cal U},\Lambda   \right)$ 
such that
\begin{itemize}
\item the initial values 
$$
\left(
\displaystyle
u^1(z_0).....u^N(z_0),
\U^{\bullet}{}^1(z_0).....\U^{\bullet}{}^N(z_0) 
  \right)
$$
of $\gamma$
yield a metrically ordinary point of $\left({\cal U},\Lambda   \right)$;
\item $u^1$ is not a constant function;
\end{itemize}
admits the following first integral:
\begin{equation}
\cases{
\displaystyle
\left(
\U^{\bullet}{}^1(z) 
\right)^2
\left( b_1\left(u^1(z)   \right)   \right)
=A_1-\sum_{l=2}^N\frac{A_l}
{a_l\left(u^1(z)   \right)}
\quad\spadesuit\cr
\displaystyle
\left(\U^{\bullet}{}^k(z)   \right)^2
f_k\left(u^k(z)   \right)
\left[a_k\left(u^1(z)   \right)\right]^2
=A_k\quad k=2...N\ \clubsuit.
}
\end{equation}
Here the $A_k$'s are suitable complex constants.
\end{lemma}
{\bf Proof:} let's prove at first that the set of equations $\,\clubsuit\,$, corresponding to $k=2...N$, holds.

If $u^k$ is a constant function, then $\U^{\bullet}{}^k\equiv 0$ and the k-th equation in $\,{\clubsuit}\,$
holds, with $A_k=0$.

If, instead, $u^k$ is not a constant function, then we could divide the k-th equation in (\ref{equazioninormali}) by $u^k$, this division being lead within the ring of meromorphic functions in a neigbhourhood of $z_0$.

We get
$$
2\frac{\U^{\bullet\bullet}{}^k(z)}{\U^{\bullet}{}^k(z)}
+\frac{f^{\prime}_k\left(u^k(z)   \right)}
          {f_k(u^k(z))}\U^{\bullet}{}^k(z)+
2\frac{a^{\prime}_k\left(u^1(z)   \right)}
          {a_k(u^1(z))}\U^{\bullet}{}^1(z)=0.
$$
Therefore, integrating once,
$$
\left(\U^{\bullet}{}^k(z)   \right)^2
f_k\left(u^k(z)   \right)
\left[a_k\left(u^1(z)   \right)\right]^2
=A_k,
$$
where we have set
$$
A_k=
\left(\U^{\bullet}{}^k(z_0)   \right)^2
f_k\left(u^k(z_0)   \right)
\left[a_k\left(u^1(z_0)   \right)\right]^2.
$$
Note that $A_k$ is a well defined complex number, since 
$$U\left(z_0   \right)=
\left(u^1(z_0)...u^N(z_0)   \right)$$ is a metrically ordinary point.

Let's eventually prove $\,\spadesuit\,$: we could multiply the first equation of (\ref{equazioninormali}) by $2b_1\left(u^1(z)   \right)\U^{\bullet}{}^1(z)$, since this last function is not everywhere vanishing.

We get
\begin{eqnarray*}
2b_1\left(u^1(z)   \right)\U^{\bullet}{}^1(z)
\U^{\bullet\bullet}{}^1(z)+{b_1^{\prime}(u^1(z))}
\left(\U^{\bullet}{}^1(z)   \right)^3\\
\qquad\quad -\sum_{l=2}^N {a_l^{\prime}(u^1(z))f_l(u^l(z))}
\left(\U^{\bullet}{}^l (z)  \right)^2\U^{\bullet}{}^1(z)=0;
\end{eqnarray*}
by $\clubsuit$, already proved,
$$
\left(\U^{\bullet}{}^l (z)  \right)^2=
\frac
{A_l}
{f_l(u^l(z))\left[a_l(u^1(z))\right
]^2},
$$
hence
\begin{eqnarray*}
& &2b_1\left(u^1(z)   \right)\U^{\bullet}{}^1(z)
\U^{\bullet\bullet}{}^1(z)+{b_1^{\prime}(u^1(z))}
\left(\U^{\bullet}{}^1(z)   \right)^3\\
& &\quad\quad -\sum_{l=2}^N 
A_l
\frac
{a_l^{\prime}(u^1(z))}
{\left[a_l(u^1(z))\right]^2}
\U^{\bullet}{}^1(z)=0;
\end{eqnarray*}
integrating once we get
$$
b_1\left(u^1(z)   \right)\left(\U^{\bullet}{}^1
(z)   \right)^2+\sum_{l=2}^N
\frac
{A_l}
{a_l\left(u^1(z)   \right)}
=K,
$$
where
$$
K=
b_1\left(u^1(z_0)   \right)\left(\U^{\bullet}{}^1
(z_0)   \right)^2+\sum_{l=2}^N
\frac
{A_l}
{a_l\left(u^1(z_0)   \right)};
$$
dividing by $b_1\left(u^1(z)   \right)$, keeping into account that $b_1\left(u^1(z_0)   \right)\not=0$ (due to the metrical ordinariness of the initial point of the geodesic) and eventually 
setting
$\displaystyle A_1=K/b_1\left(u^1(z_0)   \right)$
ends the proof.
\QUAN
\begin{lemma}\rm
\label{integraleprimo2}
Every element of geodesic 
$\displaystyle
z\longmapsto\left(u^1(z)...u^N(z)    \right)
$
of $\left({\cal U},\Lambda   \right)$ 
such that
\begin{itemize}
\item the initial values 
$$
\left(
\displaystyle
u^1(z_0).....u^N(z_0),
\U^{\bullet}{}^1(z_0).....\U^{\bullet}{}^N(z_0) 
  \right)
$$
of $\gamma$
yield a metrically ordinary point of $\left({\cal U},\Lambda   \right)$;
\item $u^1$  is a constant function;
\end{itemize}
admits the following first integral:
\begin{equation}
\cases{
\displaystyle
{u}^1(z)
=A_1\quad\diamondsuit\cr
\displaystyle
\left(\U^{\bullet}{}^k(z)   \right)^2
f_k\left(u^k(z)   \right)
=A_k\quad\heartsuit\  k=2...N.
}
\end{equation}
Here the $A_k$'s are suitable complex constants.
\end{lemma}
{\bf Proof:} $\diamondsuit\,$ holds by hypothesis: let's prove now that the set of equations $\,\heartsuit\,$, corresponding to $k=2...N$, holds.

If $u^k$ is a constant function, then $\U^{\bullet}{}^k\equiv 0$ and the k-th equation in $\,{\diamondsuit}\,$
holds, with $A_k=0$.

If, instead, $u^k$ is not a constant function, then we could divide the k-th equation in (\ref{equazioninormali}) by $u^k$, this division being lead within the ring of meromorphic functions in a neigbhourhood of $z_0$.

By keeping into account that $\U^{\bullet}{}^1(z)\equiv 0$ we get:
$$
2\frac{\U^{\bullet\bullet}{}^k(z)}{\U^{\bullet}{}^k(z)}
+\frac{f^{\prime}_k\left(u^k(z)   \right)}
          {f_k(u^k(z))}\U^{\bullet}{}^k(z)=0
$$
Therefore, integrating once,
$$
\left(\U^{\bullet}{}^k(z)   \right)^2
f_k\left(u^k(z)   \right)
=A_k,
$$
where we have set
$$
A_k=
\left(\U^{\bullet}{}^k(z_0)   \right)^2
f_k\left(u^k(z_0)   \right).
$$
Note that $A_k$ is a well defined complex number, since 
$$U\left(z_0   \right)=
\left(u^1(z_0)...u^N(z_0)   \right)$$ is a metrically ordinary point: this fact ends the proof.
\QUAN
\begin{remark}\rm
\label{radiciquadrate}
In the following we shall be concerned with 'extracting square roots' of nonvanishing elements, or germs, of holomorphic functions at some points in the complex plane: more precisely, let $\left(U,\Psi   \right)$ be a never vanishing holomorphic function element: then there exist two holomorphic function elements
$\left(U,\Xi_1   \right)$ and $\left(U,\Xi_2   \right)$
such that $\Xi_1^2=\Psi$ and $\Xi_2^2=\Psi$ on $U$: the Riemann surfaces of 
$\left(U,\Xi_1   \right)$ and $\left(U,\Xi_2   \right)$ are isomorphic, since either
\begin{itemize}
\item
the Riemann surface $\left(R,p,i,\widetilde{U}   \right)$ of 
$\left(U,\Psi   \right)$ is such that $\widetilde{U}$ is never vanishing, nor has it got any poles; then the Riemann surfaces of $\left(U,\Xi_1   \right)$, $\left(U,\Xi_2   \right)$  and $\left(U,\Psi   \right)$ are all isomorphic;
\item
the Riemann surface $\left(R,p,i,\widetilde{U}   \right)$ of 
$\left(U,\Psi   \right)$ is such that there exists some point $p\in R$ such that $\widetilde{U}(p)=0$ or such that $\widetilde{U}$ has a pole in $p$: then the function elements 
$\left(U,\Xi_1   \right)$ and $\left(U,\Xi_2   \right)$ are connectible, hence their Riemann surfaces are again isomorphic.
\end{itemize}
The same reasoning could be applied without changes to the Riemann surfaces of 
the holomorphic function elements
$\left(U,\int\Xi_1   \right)$ and $\left(U,\int\Xi_2   \right)$. 
\end{remark}
\begin{definition}
\rm
\label{coercive}
A meromorphic warped product
$$
\displaystyle
{\cal U}={\cal U}_1\times_{a_2(u^1)}
{\cal U}_2\times_{a_3(u^1)}
{\cal U}_3\times
........
\times_{a_N(u^1)}{\cal U}_N 
$$
of complex planes or one-dimensional unit balls with metric
$$
\Lambda\left(u^1.....u^N\right)=
b_1(u^1)\,du^i\odot du^i+
\sum_{i=2}^N a_i(u^i)f_i(u^i)
\,du^i\odot du^i,
$$
where $b_1$, the $a_k$'s and the $f_k$'s are 
nonzero meromorphic functions
is {\bf coercive} provided that, for every metrically ordinary
 point  $\displaystyle X_0=\left(x_0^1...x_0^N \right)$ and
\begin{itemize}
\item
for every n-tuple $\displaystyle \left(A_1...A_N   \right)\in\CI^N$ such that 
$$
\cases
{
\displaystyle
b_1(x_0^1)\not=0\cr
A_1-\sum_{l=2}^N\frac
{\displaystyle
A_l}
{\displaystyle
a_l(x_0^1)}\not=0
}
$$
and for each one of the two holomorphic function germs 
$\displaystyle{\alph_1}$
and
$\displaystyle{\alph_2}$
such that 
$$
\left({\alph_i}\right)^2=
\left[
\frac{1}{b_1}
\left(
A_1-\sum_{l=2}^N\frac
{A_l}
{a_l}
\right)
\right]
_{x^1_0}
\quad i=1,2,
$$
the Riemann surface $\displaystyle\left(S_1,\pi_1,j_1,\Phi_1,{\cal U}   \right)$
of both the holomorphic function germs
(see remark \ref{radiciquadrate})
\begin{equation}
\left[
\int_{x_0}^{u^1}
{
\frac{d\,\eta}
{\alph_i(\eta)}
}
\right]_{x_0^1}
\quad i=1,2;\label{ciuno}
\end{equation}
is such that $\PI^1\setminus\Phi_1(S_1)$ is a finite set;
\item
for each $k$, $2\leq k\leq N$
and for each one of the two holomorphic function germs 
$\displaystyle{\phi_{k1}}$
and
$\displaystyle{\phi_{k2}}$
such that 
$$
\left({\phi_{ki}}   \right)^2=
\left[
f_k
\right]_{x_0^1},
\quad i=1,2
$$
the Riemann surface $\displaystyle\left(S_k,\pi_k,j_k,\Phi_k,{\cal U}   \right)$
of both the holomorphic function germs
(see remark \ref{radiciquadrate})
\begin{equation}
\displaystyle\left[\int_{x_0^1}^{u^k}
\phi_{ki}(\eta)
\,d\eta\right]_{x_0^1}\ \ i=1,2\label{cidue}
\end{equation}
is such that $\PI^1\setminus\Phi_k(S_k)$ is a finite set.
\end{itemize}
\end{definition}
\begin{remark}
\rm
{\bf Definition \ref{coercive} may be checked
for just one metrically ordinary point $X_0$: this is proved in lemma \ref{isom}; moreover, we may assume,without loss of generality 
$X_0=0$:} 
were not, we could carry it into $0$ by
applying an automorphism of 
${\cal U}$, that is to say a direct product of automorphisms of the unit ball or of the complex plane, according to the nature of each ${\cal U}_i$. 
Then a simple pullback procedure would yield back the initial situation: {\bf in the following we shall understand this choice}.
\end{remark}
In the following lemma we shall use the 'square root' symbol in the meaning of definition \ref{coercive}, or remark \ref{radiciquadrate}: in other words, given a holomorphic function germ, which is not vanishing at some point, it should denote any one of the two holomorphic function germs yielding it back when squared.
\begin{lemma}
\label{isom}
\rm
For every metrically ordinary point $\displaystyle\left(\xi^1...\xi^N   \right)$ of ${\cal U}$ and every n-tuple $\displaystyle\left(A_1...A_N   \right)\in\CI^N$
such that
$$
\cases
{
b_1(x_0^1)\not=0\cr
A_1-\sum_{l=2}^N\frac
{
A_l}
{
a_l(x_0^1)}\not=0\cr
\ b_1(\xi^1)\not=0\cr
A_1-\sum_{l=2}^N\frac
{
A_l}
{
a_l(\xi^1)}\not=0,
}
$$
set $\Psi(\eta)=\displaystyle{ A_1-\sum_{l=2}^N\frac
{A_l}
{a_l(\eta)}}$,
the Riemann surfaces of the holomorphic function germs
$\displaystyle
\left[
\int_{\xi_1}^{u^1}
\sqrt{\displaystyle b_1(\eta)/\Psi(\eta)}
\,d\eta
\quad
\right]_{\xi_1}
$
and
$\displaystyle
\left[
\int_{0}^{u^1}
\sqrt{\displaystyle b_1(\eta)/\Psi(\eta)}
\,d\eta
\quad
\right]_{0}
$
are isomorphic: moreover so are, for each $k$,
those of
$\displaystyle
\left[
\int_{\xi_k}^{u^k}
\sqrt{
f_k(\eta)
\,d\eta
}
\right]_{\xi_k}
$
and
$\displaystyle
\left[
\int_{0}^{u^k}
\sqrt{
f_k(\eta)
\,d\eta
}
\right]_{0}
$.
\end{lemma}
{\bf Proof:}
the holomorphic function germs
$\displaystyle
\left[
\int_{\xi_1}^{u^1}
\sqrt{\displaystyle b_1(\eta)/\Psi(\eta)}
\,d\eta
\quad
\right]_{\xi_1}
$
and
$\displaystyle
\left[
\int_{0}^{u^1}
\sqrt{\displaystyle b_1(\eta)/\Psi(\eta)}
\,d\eta
\quad
\right]_{0}
$
are connectible.

Moreover so are, for each $k$,
$\displaystyle
\left[
\int_{\xi_k}^{u^k}
\sqrt{
f_k(\eta)
\,d\eta
}
\right]_{\xi_k}
$
and
$\displaystyle
\left[
\int_{0}^{u^k}
\sqrt{
f_k(\eta)
\,d\eta
}
\right]_{0}
$.
\QUAN
\vfill\eject
\section{More analytical continuation}
\begin{lemma}
\label{inverse}
\rm
Let $\left({\cal U},f   \right)$ and $\left({\cal V},g   \right)$ be two 
holomorphic function elements (or two holomorphic function germs), 
each one inverse of the other; let 
$\displaystyle\left(   R,\pi,j,F,\CI\right)$ and 
$\displaystyle\left(S,\rho,\ell,G,\CI   \right)$
be their respective standard Riemann surfaces: then $$F(R)=\rho(S).$$
\end{lemma}
{\bf Proof:} 

a) $F(R)\subset\rho(S)$: let $\xi\in R$ and $F(\xi)=\eta$; there exist:
\begin{itemize}
\item
an open neighbourhood ${\cal U}_1$ of $\xi$;
\item
open subsets ${\cal U}_2\subset\pi\left({\cal U}_1   \right)$ and ${\cal V}_2\subset F\left({\cal U}_1   \right)$;
\item
a biholomorphic function $g_2:{\cal V}_2\longrightarrow{\cal U}_2$, with inverse function $f_2:{\cal U}_2\longrightarrow{\cal V}_2$
\end{itemize}
{\TTT such that}:
\begin{itemize}
\item
$\left({\cal U}_2,f_2   \right)$ and $\left({\cal U},f  \right)$ are connectible;
\item
$\left({\cal V}_2,g_2   \right)$ and $\left({\cal V},g  \right)$ are connectible.
\end{itemize}

By construction there hence exist two holomorphic immersions 
$$\widetilde j:{\cal U}_2\longrightarrow R\hbox{ and }
\widetilde\ell:{\cal V}_2\longrightarrow S$$ 
such that $\displaystyle\pi\circ\widetilde j=\id$ and $\displaystyle\rho\circ\widetilde\ell=\id$.

Let ${\cal V}_1=F(U)_1$ and 
$$
\Sigma=\{\left(x,y   \right)\in{\cal U}_1\times{\cal V}_2: F(x)=y\};
$$
moreover let $J:{\cal V}_2\longrightarrow\Sigma$
be defined by setting $\displaystyle J(v)=\left(\widetilde j\circ  g_2(v),v \right)$.

Then $\displaystyle \left(\Sigma,pr_2,J,\pi\circ pr_1   \right)$ is an analytical continuation of $\left({\cal V}_2,g_2   \right)$; indeed $$\displaystyle \pi\circ pr_1\circ J=\pi\circ\widetilde j\circ g_2=g_2.$$ But $\left({\cal V_2},g_2   \right)$ is connectible with $\left({\cal V},g     \right)$, hence
$\displaystyle \left(\Sigma,pr_2,J,\pi\circ pr_1   \right)$ is an analytical continuation of $\left({\cal V},g \right)$.

There eventually exists a holomorphic function $h:\Sigma\longrightarrow S$ such that $\rho\circ h=pr_2$: hence 
$$
\eta=pr_2\left(\xi,\eta   \right)=\rho\circ h
\left(\xi,\eta   \right)\in\rho\left(S   \right).
$$
\vskip1truecm
b) $\rho(S)\subset F(R)$: let $s\in S$: there is a neighbourhood  $V$ of $s$ in $S$ such that $V\setminus\{s\}$ consists entirely of regular points both of $\rho$ and $G$, not excluding that $s$ itself be regular for $\rho$ or $G$ or both.

This fact means that for each $s^{\prime}\in V\setminus\{s\}$ there exists a holomorphic function element $\left({\rho(s^{\prime})},{\cal V}^{\prime},\widetilde g_{s^{\prime}}   \right)$ connectible with $\left({\cal V},g   \right)$ and, besides, a holomorphic immersion $\widetilde\ell:{\cal V}^{\prime}\longrightarrow V$.

By a) already proved, $G(s)\in\pi(R)$, hence there exist:
\begin{itemize}
\item
$p\in R$ such that $\pi(p)=G(s)$;
\item
a neighbourhood $W$ of $p$ in $R$ such that $\displaystyle \pi^{-1}\left(\widetilde g\left({\cal V}^{\prime}   \right)   \right)\bigcap W\not=\emptyset$.
\end{itemize}

Set
$$\displaystyle W^{\prime}=\pi^{-1}\left(\widetilde g\left({\cal V}^{\prime}   \right)   \right)\bigcap W:$$
we may suppose, without loss of generality, that $\pi$ is invertible on $W^{\prime}$: hence there exists a (open) holomorphic immersion $\displaystyle\widetilde j:\widetilde g\left({\cal V}^{\prime}   \right)\longrightarrow W$.

Therefore, for each $\displaystyle\zeta\in\widetilde j\left(\widetilde g\left({\cal V}^{\prime}   \right)   \right)$, there exists $\displaystyle \eta\in\widetilde\ell\left({\cal V}^{\prime}   \right)$ such that $$\displaystyle F(\zeta)=F\left(\widetilde j \circ \widetilde g\circ \rho(\eta)  \right).$$

Now, by definition of analytical continuation there holds $$\displaystyle F\circ\widetilde j\circ\widetilde g=\id,$$ hence we have $\displaystyle F(\zeta)=\rho(\eta)$.

Consider now the holomorphic function 
$\displaystyle\Xi:W\times V\longrightarrow\CI$ 
defined by setting $$\displaystyle \Xi\left(w,v   \right)=F(w)-\rho(v):$$ we have
$$
\left.\Xi\right\vert_{\widetilde j\left(\widetilde g\left({\cal V}^{\prime}   \right)   \right)\times\widetilde\ell({\cal V}^{\prime})}
\equiv 0.
$$

But $\displaystyle {\widetilde j\left(\widetilde g\left({\cal V}^{\prime}   \right)   \right)\times\widetilde\ell({\cal V}^{\prime})}$ is an open set in $W\times V$, hence $\Xi\equiv 0$ on $W\times V$, which in turn implies $\displaystyle F(p)=\rho(s)$.

Therefore we have proved that for each $s\in S$ there exists $p\in R$ such that $ F(p)=\rho(s)$: this eventually implies that $\rho(S)\subset F(R)$.
\QUAN

\begin{lemma}
\label{quasiinverse}
\rm
Let 
$\f$, $\g$, $\h$
be three holomorphic function germs such that
$\f\circ \g=\h$.
Let $\displaystyle\left(   R,\pi,j,F,\CI\right)$ 
be the Riemann surface of $\f$, $\displaystyle\left(S,\rho,\ell,G,CI   \right)$
the one of $\g$ and 
$\displaystyle\left(T,\sigma,m,H,CI  \right)$
the Riemann surface with logarithmic singularities (see definition \ref{loganalcont}) of $\h$: then $$F(R)\setminus\left(\PI^1\setminus\left(
\sigma(T)
   \right)   \right)\subset\rho(S).$$
\end{lemma}
{\bf Proof:} let $\xi\in R$ such that $\displaystyle\eta=F(\xi)\not\in\PI^1
\setminus\left(
\sigma(T)
   \right)  $; there exist:
\begin{itemize}
\item
an open neighbourhood ${\cal U}_1$ of $\xi$;
\item
open subsets ${\cal U}_2\subset\pi\left({\cal U}_1   \right)$, ${\cal V}_2\subset F\left({\cal U}_1   \right)$ and  ${\cal W}_2\subset \sigma\left(T  \right)$ ;
\item
biholomorphic functions 
$\displaystyle f_2:{\cal U}_2\longrightarrow{\cal W}_2$,
$\displaystyle g_2:{\cal V}_2\longrightarrow{\cal U}_2$ and
$\displaystyle h_2:{\cal V}_2\longrightarrow{\cal W}_2$
\end{itemize}
{\TTT such that}:
\begin{itemize}
\item
$\left({\cal U}_2,f_2   \right)$ and $\f$ are connectible;
\item
$\left({\cal V}_2,g_2   \right)$ and $\g$ are connectible;
\item
$\left({\cal V}_2,h_2   \right)$ and $\h$ are connectible;
\item
$f_2\circ g_2=h_2$.
\end{itemize}

By construction there hence exist three holomorphic immersions 
$$
\cases
{
\widetilde j:{\cal U}_2\longrightarrow R\cr
\widetilde\ell:{\cal V}_2\longrightarrow S\cr 
\widetilde m:{\cal W}_2\longrightarrow T
}
$$ 
such that $\displaystyle\pi\circ\widetilde j=\id$, $\displaystyle\rho\circ\widetilde\ell=\id$ and
$\displaystyle\sigma\circ\widetilde m=\id$.

Let ${\cal V}_1=F\left(U_1\right)$, ${\cal W}_1$ be the connected component of $\displaystyle \sigma^{-1}\left(F\left({\cal U}_1   \right)   \right)$ in $T$ and
$$
\Sigma=\{\left(x,y   \right)\in{\cal U}_1\times{\cal W}_1: F(x)=H(y)\};
$$
moreover let $J:{\cal V}_2\longrightarrow\Sigma$
be defined by setting $\displaystyle J(v)=\left(\widetilde j\circ  g_2(v),\widetilde m(v) \right)$.

Then $\displaystyle \left(\Sigma,pr_2,J,\pi\circ pr_1   \right)$ is an analytical continuation, with logarithmic singularities of $\left({\cal V}_2,g_2   \right)$; indeed 
$$\displaystyle \pi\circ pr_1\circ J=\pi\circ\widetilde j\circ g_2=g_2.$$

But $\left({\cal V_2},g_2   \right)$ is connectible with $\left({\cal V},g     \right)$, hence
$\displaystyle \left(\Sigma,pr_2,J,\pi\circ pr_1   \right)$ is an analytical continuation of $\left({\cal V},g \right)$.

There eventually exists a continuous function $h:\Sigma\longrightarrow S$ (in fact holomorphic on $int(\Sigma)$)
such that $\rho\circ h=pr_2$ : hence 
$$
\eta=pr_2\left(\xi,\eta   \right)=\rho\circ h
\left(\xi,\eta   \right)\in\rho\left(S   \right).
$$
\QUAN
\vfill\eject
\section{The main theorem}
\begin{theorem}
\TTT
\label{teoremaprincipale}
A meromorphic warped product
$$
\displaystyle
{\cal U}={\cal U}_1\times_{a_2(u^1)}
{\cal U}_2\times_{a_3(u^1)}
{\cal U}_3\times
........
\times_{a_N(u^1)}{\cal U}_N 
$$
of complex planes or one-dimensional unit balls with metric
$$
\Lambda\left(u^1.....u^N\right)=
b_1(u^1)\,du^1\odot du^1+
\sum_{i=2}^N a_i(u^1)f_i(u^i)
\,du^i\odot du^i,
$$
is geodesically complete if and only if it is coercive.
\end{theorem}
{\bf Proof:} 

a) suppose that $\displaystyle{\cal U}$ is {\TTT coercive} and that $U$, defined by
$$
z\longmapsto\left(u^1(z)...u^N(z)    \right),
$$
is an element of geodesic, defined in a neighbourhood of $0$ in the complex plane and such that $\left(u^1(0)...u^N(0)    \right)$ is a metrically ordinary point; moreover, let 
$$
\left(\U^{\bullet}{}^1(0)...\U^{\bullet}{}^N(0)    \right)
$$
be the initial velocity of $U$.

Suppose at first that $\displaystyle z\mapsto u^1(z)$ is a constant function (hence $\U^{\bullet}{}^1(0)=0\,$): then, by lemma \ref{integraleprimo2}, the equations of $U$ are
\begin{equation}
\label{riportointpr2}
\cases{
\displaystyle
{u}^1(z)
=A_1\quad\cr
\displaystyle
\left(\U^{\bullet}{}^k(z)   \right)^2
f_k\left(u^k(z)   \right)
=A_k\quad\  k=2...N,
}
\end{equation}
where the $A_k$'s are suitable complex constants.

Now the Riemann surface of the holomorphic function element $\displaystyle z\mapsto {u}^1(z)$ is trivially isomorphic to $\displaystyle \left(\PI^1,\id,\id,A_1 \right)$; if  $\displaystyle A_k=0$ the Riemann surface of $\displaystyle z\mapsto {u}^k(z)$ is isomorphic to $\displaystyle \left(\PI^1,\id,\id,A\right)$ for some complex constant $A$; if $\displaystyle A_k\not=0$ we could rewrite the k-th equation of (\ref{riportointpr2}) in the form:
\begin{equation}
\label{riscritta}
\frac{1}{B_k}\int_{u^k(0)}^{u^{k}(z)} \phi(\eta)\,d\eta
=\,z,
\end{equation}
where $\phi_k^2=f_k$ and $B_k^2=A_k$, the choice of $\phi_k$ and $B_k$ being made in such a way that 
$$\displaystyle \U^{\bullet}{}^k(0)=\frac{B_k}{\phi_k(0)}.$$
By hypothesis, the Riemann surface $\displaystyle\left(S_k,\pi_k,j_k,\Phi_k   \right)$
of the holomorphic function germ $\displaystyle\left[\int_0^{u^k}\phi_k\,d\eta\right]_{0}$
is such that $\displaystyle\PI^1\setminus\Phi_1(S_1)$ is a finite set; by lemma \ref{isom} the Riemann surface of the holomorphic function germ $\displaystyle\left[\int_{u^k(0)}^{u^k}\phi_k\,d\eta\right]_{u^k(0)}$
is isomorphic to $\displaystyle\left(S_k,\pi_k,j_k,\Phi_k   \right)$; but, by (\ref{riscritta}), the germs
$\displaystyle \u^k_{z=0}$ and
$\displaystyle\left[\int_{u^k(0)}^{u^k}\phi_k\,d\eta\right]_{u^k(0)}$
are each one inverse of the other; hence, by lemma \ref{inverse} the Riemann surface of $\displaystyle \u^k_{z=0}$ is complete; this eventually implies that the Riemann surface of the element
$$
z\longmapsto\left(u^1(z)...u^N(z)    \right)
$$
is complete too: this fact ends the proof of a) in the case that $u^1$ is a constant function.

On the other side, suppose that $u^1$ is not a constant function: then, by lemma \ref{integraleprimo1}, the equations of $U$ are
\begin{equation}
\label{riportointpr1}
\cases{
\displaystyle
\left(\U^{\bullet}{}^1(z)   \right)^2
\left( b_1\left(u^1(z)   \right)   \right)
=A_1-\sum_{l=2}^N\frac{A_l}
{a_l\left(u^1(z)   \right)}
\quad\spadesuit\cr
\displaystyle
\left(\U^{\bullet}{}^k(z)   \right)^2
f_k\left(u^k(z)   \right)
\left[a_k\left(u^1(z)   \right)\right]^2
=A_k\quad k=2...N\ \clubsuit.}
\end{equation}
for suitable complex constants $A_1...A_N.$

Consider now the germ $\displaystyle z\mapsto u^1(z)$ in $z=0$:
rewrite the first equation of (\ref{riportointpr1}) in the form:
\begin{equation}
\label{riscritta2}
\int_{u^1(0)}^{u^{1}(z)}
\frac{\displaystyle d\eta}
 {\displaystyle\alph(\eta)_{u^1(0)}}
=\,z,
\end{equation}
where 
$$
\left(
\alph(\eta)_{u^1(0)}
\right)^2=
\displaystyle\frac
{A_1-\sum_{l=2}^N
A_l
/
a_l(\eta)}
{b_1(\eta)}
$$ 
in a neighbourhood of $z=0$,
the choice of the square root  $\alph_k$ being 
made in such a way that
$
\displaystyle \alph_{u^1(0)}\left(u^1(0)\right)=1/\U^{\bullet}{}^1(0)
$.

Denote now  by $\alph_{u=0}$ the holomorphic function germ 
defined by setting
$$
\left(
\alph_{0}
\right)^2=
\left[
\frac
{1}
{\displaystyle b_1}
\left(
{\displaystyle
A_1-\sum_{l=2}^N\displaystyle\frac
{A_l}
{a_l}}
\right)
\right]_{0}
,
$$
the choice of the 'square root' $\alph_0$ being arbitrary.

By hypothesis, the Riemann surface $\displaystyle\left(S_1,\pi_1,j_1,\Phi_1   \right)$
of the holomorphic function germ $\displaystyle\left[\int_0^{u^1}1/\alph_0\right]_0$
is such that $\displaystyle\PI^1\setminus\Phi_1(S_1)$ is a finite set.

By lemma \ref{isom} the Riemann surfaces of $\displaystyle\left[\int_0^{u^1}1/\alph_0\right]_0$ and of $\displaystyle\left[\int_{u^1_0}^{u^1}1/\alph_0\right]_{u^1_0}$
are both isomorphic to $\displaystyle\left(S_1,\pi_1,j_1,\Phi_1   \right)$; but, by (\ref{riscritta}), the germs
$\displaystyle \u^1_{z=0}$ and
$\displaystyle[\int_0^{u^1}1/\alph_0]_{u^1(0)}$ are each one inverse of the other; hence, by lemma \ref{inverse} the Riemann surface of $\displaystyle \u^1_{z=0}$ is complete.

Let now $\displaystyle 2\leq k\leq N$: if $\displaystyle A_k=0$ the Riemann surface of $\displaystyle z\mapsto {u}^k(z)$ is isomorphic to $\displaystyle \left(\PI^1,\id,\id,A\right)$ for some complex constant $A$; if $\displaystyle A_k\not=0$ we could rewrite the k-th equation of (\ref{riportointpr1}) in the form:
\begin{equation}
\label{riscritta3}
\int_{u^k(0)}^{u^{k}(z)} \phi(\eta)\,d\eta
=\,
\int_0^z\frac{B_k\,dz}{a_k\left(u^1(z)   \right)},
\end{equation}
where $\phi_k^2=f_k$ and $B_k^2=A_k$, the choice of $\phi_k$ and $B_k$ being made in such a way that 
$$\displaystyle \U^{\bullet}{}^k(0)\,\phi\left(
u^k(0)   \right)\,
a_k\left(u^1(z)   \right)={B_k}.$$

Denote now  by $\displaystyle[\varphi_k]_{u^k=0}$ the holomorphic function germ 
defined by setting
$\displaystyle[\varphi_k]_{u^k=0}^2=
\left[
f_k
\right]_{u^k=0}
$,
the choice of the "square root" $\displaystyle[\varphi_k]_{u^k=0}$ being arbitrary.

By hypothesis, the Riemann surface $\displaystyle\left(S_k,\pi_k,j_k,\Phi_k   \right)$
of the holomorphic function germ $\displaystyle\left[\int_0^{u^k}\varphi_k\right]_{0}$
is such that $\displaystyle\PI^1\setminus\Phi_1(S_1)$ is a finite set; moreover, by lemma \ref{isom} the Riemann surfaces of the holomorphic function germ $\displaystyle\left[\int_{u^k(0)}^{u^k}\phi_k\,d\eta\right]_{u^k(0)}$
is isomorphic to $\displaystyle\left(S_k,\pi_k,j_k,\Phi_k   \right)$; but, by (\ref{riscritta3}) the germs
\begin{itemize}
\item
$\displaystyle \left[z\longrightarrow\u^k\right]_{z=0}$;
\item
$\displaystyle\left[\int_{u^k(0)}^{u^k}\phi_k\,d\eta\right]_{u^k(0)}$;
\item
$\displaystyle
\left[z\longrightarrow\int_0^z\frac{B_k}{a_k\left(u^1(\zeta)   \right)}\,d\zeta\right]_{z=0}$
\end{itemize}
satisfy, in the above order, the hypotheses of lemma \ref{quasiinverse};
moreover, the Riemann surface with logarithmic singularities of $\displaystyle\left[\int_{u^k(0)}^{u^k}\phi_k\,d\eta\right]_{u^k(0)}$ is complete, since the one of 
$\displaystyle\left[\phi_k\right]_{u^k(0)}$ is complete without logarithmic singularities.

Therefore the Riemann surface with logarithmic singularities of $\displaystyle \u^k_{z=0}$ is complete; this eventually implies that the Riemann surface with logarithmic singularities of the element
$$
z\longmapsto\left(u^1(z)...u^N(z)    \right),
$$
is complete too: this fact ends the proof of a).
\vskip1truecm
Vice versa, suppose that 
$$
\displaystyle
{\cal U}={\cal U}_1\times_{a_2(u^1)}
{\cal U}_2\times_{a_3(u^1)}
{\cal U}_3\times
........
\times_{a_N(u^1)}{\cal U}_N 
$$
is not {\TTT coercive}: then
\begin{itemize}
\item
{\bf either}
there exists a complex
n-tuple $\displaystyle \left(A_1...A_N   \right)\in\CI^N$ such that 
$$
\cases
{
\displaystyle
b_1(x_0^1)\not=0\cr
A_1-\sum_{l=2}^N\frac
{\displaystyle
A_l}
{\displaystyle
a_l(x_0^1)}\not=0
}
$$
and for each one of the two holomorphic function germs 
$\displaystyle{\alph_1}$
and
$\displaystyle{\alph_2}$
such that 
$$
\left({\alph_i}\right)^2=
\left[
\frac{1}{b_1}
\left(
A_1-\sum_{l=2}^N\frac
{A_l}
{a_l}
\right)
\right]
_{0}
\quad i=1,2,
$$
the Riemann surface $\displaystyle\left(S_1,\pi_1,j_1,\Phi_1   \right)$
of both the holomorphic function germs
(see remark \ref{radiciquadrate})
$$
\left[
\int_{x_0}^{u^1}
{
\frac{d\,\eta}
{\alph_i(\eta)}
}
\right]_{x_0^1}
\quad i=1,2;
$$ 
is such that $\PI^1\setminus\Phi_1(S_1)$ is
an infinite set;
\item
{\bf or}
there exists $k$, $2\leq k\leq N$
such that, for each one of the two holomorphic function germs 
$\displaystyle\left[\phi_{k1}\right]_0$
and
$\displaystyle\left[\phi_{k2}\right]_0$
such that 
$$
\displaystyle\left[\phi_{ki}\right]_0=
\left[
f_k
\right]_0,
\quad i=(1,2)
$$
the Riemann surface $\displaystyle\left(S_k,\pi_k,j_k,\Phi_k   \right)$
of both the holomorphic function germs
(see remark \ref{radiciquadrate})
$\displaystyle\left[\int_0^{u^k}
\phi_ki(\eta)
\,d\eta\right]_{0}\ \ i=1,2$
is such that $\PI^1\setminus\Phi_1(S_1)$ is an infinite set.
\end{itemize}

In the first case the geodesic element
$$
z\longmapsto U(z)=\left(u^1(z)...u^N(z)    \right)
$$
starting from $0$ with velocity 
$\displaystyle\left(L_1...L_N    \right)$, such that
$$\displaystyle 
L_1^2=
\frac
{1}
{b_1(0)}
\left(
A_1-\sum_{l=2}^N\frac
{A_l}
{a_l(0)}
\right)
,\ L_k^2=\frac{A_k}{f_k(0)a_k(0)},\ 
k=2...N
,$$
satisfies, among the other ones, the equation
$$
\int_0^{u^1(z)}
\frac
{\displaystyle d\eta}
{\displaystyle \alph_i(\eta)}
=z,
$$
where $i=1$ or $i=2$;
by lemma \ref{inverse}, this fact implies that $\left[
\displaystyle z\longmapsto
u^1(z)
\right]_{0}  $
has got an incomplete Riemann surface, hence  the same holds about $\displaystyle z\longmapsto
U(z)$ too.

Consider now the second case: first construct a geodesic element 
$$
z\longmapsto U(z)=\left(0...u^k(z)...0    \right),
$$
with all components which have to be constant functions except $\displaystyle u^k, k\geq 2$ (this element is easily seen to exist).

Now recall lemma \ref{integraleprimo2} to conclude that 
$
\displaystyle
z\longmapsto u^k(z)
$
satisfies, in a neighbourhood of $z=0$ the equation 
$$
\displaystyle
\frac
{1}
{C_k}
{\int_0^{u^k(z)} \phi_{ki}(\eta)\,d\eta}=z,
$$
for a suitable complex constant $A_k$; therefore its Riemann surface is incomplete by lemma \ref{inverse}; this fact ends the proof.
\QUAN
\begin{definition}
\label{directbiholom}
\rm
Let ${\cal U}$ and ${\cal V}$ be meromorphic warped products of complex planes and unit balls;  ${\cal U}$ and ${\cal V}$ are {\bf directly biholomorphic} provided that they are biholomorphic under a direct product of biholomorphic functions between each  ${\cal U}_i$ and each ${\cal V}_i$.
\end{definition}
\begin{remark}
\label{biholom}
\rm
Definition \ref{coercive} is invariant under direct biholomorphism (see definition \ref{directbiholom}
) : in other words, if ${\cal U}$ and ${\cal V}$ are directly biholomorphic, then 
${\cal U}$ is coercive if and only ${\cal V}$ is too: this is a simple consequence of 'changing variable' in integrals \ref{ciuno} and \ref{cidue}.
\end{remark}
Therefore, we could yield the following
\begin{definition}
\label{equicoercive}
\rm
An equivalence class $\left[{\cal U}    \right]$ of meromorphic warped products of complex planes and unit balls, consisting of mutually directly (see definition \ref{directbiholom}
) biholomorphic elements is {\bf coercive}
provided that any one of its representatives is coercive. 
\end{definition}

Our goal is now to extend definitions \ref{coercive} and \ref{equicoercive} to warped products containg some $\PI^1$'s among their factors.

Keeping into account remark \ref{biholom}, this could be readiliy pursued: indeed, consider a meromorphic warped product 
$$
\displaystyle
{\cal U}={\cal U}_1\times_{a_2(u^1)}
{\cal U}_2\times_{a_3(u^1)}
{\cal U}_3\times
........
\times_{a_N(u^1)}{\cal U}_N 
$$
of Riemann spheres, complex planes or one-dimensional unit balls with metric
$$
\Lambda\left(u^1.....u^N\right)=
b_1(u^1)\,du^i\odot du^i+
\sum_{i=2}^N a_i(u^i)f_i(u^i)
\,du^i\odot du^i.
$$

Let $\displaystyle L\subset \{1...N\}$ be the set of indices such that $\displaystyle {\cal U}_l\simeq \PI^1$ for each $l\in L$.

\begin{definition}
\rm
Let $\displaystyle Y=\left(y^1...y^N \right)\in{\cal U}$: then $\left(Y,L    \right)$ is a {\bf principal multipole} of ${\cal U}$ provided that 
$$
\cases
{
b_1(y^1)=\infty\cr
f_l(y^l)=\infty & for each $l\in L\setminus \{1\}$.\cr
}
$$
\label{multipole}
\end{definition}
\begin{definition}
\rm
A meromorphic warped product 
$$
\displaystyle
{\cal U}={\cal U}_1\times_{a_2(u^1)}
{\cal U}_2\times_{a_3(u^1)}
{\cal U}_3\times
........
\times_{a_N(u^1)}{\cal U}_N 
$$
of Riemann spheres, complex planes or one-dimensional unit balls with metric is {\bf partially projective} if some one of its factors is biholomorphic to the Riemann sphere $\PI^1$.
\end{definition}
\begin{definition}
\rm
\label{multicoer}
A partially projective warped product $\displaystyle {\cal U}=
\prod_{i=1}^N{\cal U}_i$ is {\bf coercive in opposition to the principal multipole } $\left(Y,L    \right)$ if, set 
$$
{\cal W}_i=
\cases
{
{\cal U}_i & if $i\not\in L$\cr
{\cal U}_i\setminus \{y^i\} & if $i\in L$,
}
$$
then $\displaystyle \prod_{i=1}^N {\cal W}_i$ is coercive in the sense of definition \ref{equicoercive}, that is to say, belongs to a coercive equivalence class with respect to direct biholomorphicity.
\end{definition}
\vfill\eject
\section{Warped product of Riemann surfaces}

Consider now a warped product of Riemann surfaces
$$
\displaystyle
{\cal S}={\cal S}_1\times_{a_2}
{\cal S}_2\times_{a_3}
{\cal S}_3\times
........
\times_{a_N}{\cal S}_N,
$$
where each ${\cal S}_i$ is endowed with meromorphic metric $\lambda_i$: ${\cal S}$'s metric $\Lambda$ is defined by setting
$$
\Lambda=\lambda_1+\sum_{k=2}^N a_k\lambda_k,
$$
where each $a_k$ is a meromorphic function on ${\cal S}_i$.

\begin{theorem}
\label{rivestimento}
\TTT
${\cal S}$ admits universal covering $\displaystyle \Psi : {\cal U}\longrightarrow {\cal S}$, where ${\cal U}$ is a direct product 
of Riemann spheres, complex planes or one-dimensional unit balls: this universal covering is unique up to direct biholomorphisms.
\end{theorem}
{\bf Proof}: this is a simple consequence of Riemann's uniformization theorem.
\QUAN
Now ${\cal U}$ could be endowed with the pull-back meromorphic metric $\displaystyle \Psi^*\Lambda$, hence ${\cal U}$ itself results in a meromorphic warped product.
\begin{definition}
\rm
\label{rieco}
The manifold ${\cal S}$ is
\begin{itemize}
\item {totally unelliptic} provided that none of the ${\cal S}_i$ is elliptic;
\item {$L$-elliptic} provided that there exists a nonempty set of indices $L$ such that ${\cal S}_l$ is elliptic if and only if $l\in L$.
\end{itemize}
\end{definition}
\begin{definition}
\rm
Let ${\cal S}$ be a $L$-elliptic warped product, with universal covering $\Psi : {\cal U}\longrightarrow{\cal S}$: then $\left(Z,L    \right)$ is a principal multipole for ${\cal S}$ provided that $Z\in {\cal S}$ and each $Y\in\Psi^{-1}\left(Z    \right)$ is a principal multipole for ${\cal U}$.
\end{definition}
\begin{definition}
\rm
\begin{itemize}
\item
A totally unelliptic warped product of Riemann sur\-faces is {\bf coercive} provided that its universal covering is coercive in the sense of definition \ref{equicoercive};
\item
A $L$-elliptic warped product of Riemann surfaces is {\bf coercive in opposition to the principal multipole} $\displaystyle\left(Z,L    \right)$ provided that its universal covering ${\cal U}$ is coercive in opposition to each principal multipole $\displaystyle\left(Y,L    \right)$ as $Y$ runs over $\Psi^{-1}(Z)$.
\end{itemize}
\end{definition}
\begin{theorem}
\label{unell}
\TTT
A totally unelliptic warped product of Riemann surfaces ${\cal S}$ is geodesically complete if and only if it is coercive.
\end{theorem}
{\bf Proof:} let $\Psi : {\cal U}\longrightarrow
{\cal S}$ be the universal covering of ${\cal S}$:
by definition \ref{rieco} ${\cal U}$ is coercive, hence geodesically complete by theorem \ref{teoremaprincipale}.

Let now $\gam$ be a germ of geodesic in ${\cal S}$, starting at a metrically ordinary point: since $\Psi$ is a local isometry, there exists a germ $\bet$ of geodesic in ${\cal U}$, starting at a metrically ordinary point, such that $\displaystyle\gam=\Psi\circ\bet$.

By definition of completeness, the Riemann surface with logarithmic singularities $\displaystyle\left(\Sigma,\pi,j,B,{\cal U}    \right)$ of $\bet$ is such that $\PI^1\setminus \pi\left(\Sigma    \right)$ is a finite set; moreover, $\displaystyle\left(\Sigma,\pi,j,\Psi\circ B,{\cal S}    \right)$ is an analytical continuation, with logarithmic singularities, of $\gam$.

This proves that, if $\displaystyle\left(\widetilde \Sigma,\widetilde \pi,\widetilde j,G,{\cal S}    \right)$ is the Riemann surface with logarithmic singularities of $\gam$, then $\PI^1\setminus \widetilde \pi\left(\widetilde \Sigma    \right)$ is a finite set too, hence ${\cal S}$ is geodesically complete.

On the other side, if ${\cal S}$ admits an incomplete germ of geodesic $\gam$, starting at a metrically ordinary point, then there exists  an incomplete germ of geodesic $\bet$ in ${\cal U}$, starting at a metrically ordinary point, such that 
$\displaystyle\gam=\Psi\circ\bet$; this means by theorem \ref{teoremaprincipale},
that ${\cal U}$ is not coercive; eventually, by definition \ref{rieco}, ${\cal S}$ is not coercive: this fact ends the proof.
\QUAN
\begin{theorem}
\TTT
A $L$-elliptic warped product of Riemann surfaces ${\cal S}$ is geodesically complete if and only if) it is coercive in opposition to some principal  multipole.
\end{theorem}
{\bf Proof:} suppose that ${\cal S}$ is coercive in opposition to some principal  multipole $\displaystyle\left(Z,L    \right)$: then, by theorem \ref{unell},  ${\cal S}$ is coercive in opposition to $\displaystyle\left(Z,L    \right)$ if and only if ${\cal S}\setminus Z$ is geodesically complete;
since $Z$ is not metrically ordinary, ${\cal S}$ is geodesically complete.

On the other hand, suppose that ${\cal S}$ admits an incomplete geodesic
$\displaystyle\left(\Sigma,\pi,j,\gamma,{\cal S}    \right)$: let $\left(Z,L    \right)$ be a principal multipole of ${\cal S}$ wich is known to exist; set $R=\gamma^{-1}\left({\cal S}\setminus Z    \right)\subset\Sigma$.

Now $\displaystyle\left(R,\pi\vert_R,j,\gamma\vert_R,{\cal S}\setminus Z    \right)$ is an incomplete geodesic of ${\cal S}\setminus Z$: this fact implies that ${\cal S}\setminus Z$ is not geodesically complete, hence it is not coercive, that is to say, ${\cal S}$ is not coercive in opposition to $\left(Z,L    \right)$.

The arbitrariness of $Z$ allows us to conclude the proof.
\QUAN
\vfill\eject
\section{Examples}
We shall now show a wide class of warped products sharing all characteristics defining coercivity: they will hence result in being geodesically complete.

We recall, without proof, the following results from the theory of meromorphic functions (see \cite{nevanlinna} or \cite{hayman}):
\begin{theorem}
\label{mero}
\TTT
\begin{itemize}
\item
A function meromorphic in the complex plane takes all $\PI^1$'s values but at most two ones;
\item
a function meromorphic in the unit disc, whose characteristic function $T$ is such that the ratio
\begin{equation}
T(r)/log(1-r)\label{carat}
\end{equation}
is unlimited as $r\longrightarrow 1$, takes all $\PI^1$'s values but at most two ones.
\end{itemize}
\end{theorem}

In the following we shall need some technicalities from integral calculus, hence we state:
\begin{proposition}\rm
\label{integrale}
Set $\Delta=b^2-4ac$, there holds
$$
\left[
\int \frac{\displaystyle d\,\eta}
{\sqrt{a\eta^2+b\eta+c}}
\right]_0 
=
\cases{
\cases{
\displaystyle
\left[
\frac{1}{\sqrt{a}}\log\left(\eta+\frac{b}{2a}+
\sqrt{\eta^2+\frac{b}{a}\eta+\frac{c}{a}}    \right)+cost
\right]_0
\cr
\hbox{\sl the same branch of $\sqrt{}$, any branch of log}\cr
\hbox{\sl if $a\not=0$ and $\Delta\not=0$}
}
\ \cr
\cases
{
\displaystyle
\left[
\frac{1}{\sqrt{a}}\log\left(\eta+\frac{b}{2a}+
\right)+cost
\right]_0
\cr
\hbox{\sl any branch of log}\cr
\hbox{\sl if $a\not=0$ and $\Delta=0$}\cr
} 
\ \cr
\cases
{
\left[
\frac{\displaystyle 2}{\displaystyle b}\sqrt{b\eta+c}+cost 
\right]_0
\cr
\hbox{\sl the same branch of $\sqrt{}$}
\cr
\hbox{\sl if $a=0$ and $b\not=0$}
} 
\ \cr
\cases
{
\left[ 
\eta/\sqrt{c}+cost\right]_0\cr
\hbox{\sl the same branch of $\sqrt{}$}
\cr
\hbox{\sl if $a=b=0$.}
}
}
$$
\end{proposition}
\QUAN
Let now $h,f_2...f_N$ be mermorphic functions on $\CI$ and $P_2...P_N$ polynomials of degree at most two.

Consider on $\CI^N$ the mermorphic metric 
$$
\Lambda\left(u^1...u^N\right)=
\left[h^{\prime}(u^1)\right]^2 du^1\odot du^1+
\sum_{k=2}^N
\frac{\displaystyle\left[f_k(u^k)\right]^2}
{P_k\left(h(u^1)\right)}\, du^k\odot du^k.
$$
\begin{theorem}
\TTT
$\displaystyle\left(\CI^N,\Lambda\right)$ is coercive (hence geodesically complete).
\label{esempio}
\end{theorem}
{\bf Proof:}
\begin{itemize}
\item
For every n-tuple $\displaystyle \left(A_1...A_N   \right)\in\CI^N$ such that 
$$
\cases
{
h^{\prime}(0)\not=0\cr
\displaystyle
A_1-\sum_{l=2}^N
{A_l}
{P_l(0)}\not=0,
}
$$
there holds
\begin{eqnarray*}
\int_0^{u^1}
{
\frac
{\displaystyle h^{\prime}(\eta)d\,\eta}
{\sqrt{\displaystyle A_1-\sum_{l=2}^N
{A_l}
{P_l(h(\eta))}}}
}
&=&
\int_{h(0)}^{h(u^1)}
{
\frac
{\displaystyle d\,h}
{\sqrt
{\displaystyle A_1-\sum_{l=2}^N
{A_l}
{P_l(h)}}}
}\\
&=&
\Phi\left(h(u^1)    \right),
\end{eqnarray*}
where $\Phi$ is one (depending on the constants $A_1...A_N$) of the holomorphic function germs on the right hand member of proposition \ref{integrale}.
This fact shows that the maximal analytical continuation of 
$u^1\longrightarrow\Phi\left(h(u^1)    \right)$ takes all $\PI^1$'s values but a 
finite number, because so does the meromorphic function $h$ (see theorem \ref{mero});
\item
for each $k$, $2\leq k\leq N$,
each one of the two holomorphic function germs 
$
\pm
\left[
f_k
\right]_0,
$
could be continuated to $\pm f_k$ which, by theorem \ref{mero}, takes all values but at most two ones.
\end{itemize}
\QUAN
\begin{remark}
\rm Extending the validity of preceding example to the partially projective case is straightforward.
\end{remark}

Let now $S_i,\ i=1..N$ be Riemann surfaces, which we suppose for simplicity to be 
parabolic or hyperbolic, $p_i\colon{\cal U}_i\longrightarrow S_i$ their universal covering, 
where each ${\cal U}_i$ is isomorphic either to the unit disc or to the complex plane;
finally, let $\phi_i$ be meromorphic functions such that $\phi_1\circ p_1$ and 
$(\phi_i\circ p_i)^{\prime},\ i=1..N$ take all complex values but at most a finite number (the hypothesis on $\phi_i\circ p_i$ could be weakened; even dropped, if $S_i$ is parabolic: see \cite{hayman}).

Moreover, let $a_i,\ b_i,\ c_i,\ i=1..N $ be complex numbers such that, for each $i$, $a_i\not=0$ or  $b_i\not=0$ or $\ c_i\not=0$.

Set 
$$
\cases
{
S=\prod_{i=1}^N S_i\cr
{\cal U}=\prod_{i=1}^N={\cal U}_i\cr
p=(p_1....p_N) 
}
$$ 
and consider the meromorphic metric
$$
\Lambda=d\phi_1\odot d\phi_1+\sum
_{i=1}^N \frac {d\phi_i\odot d\phi_i}{a_i\phi_1^2+b_i\phi_1+c_i}.
$$
\begin{theorem}\TTT
$({\cal U},\Lambda)$ is coercive (hence geodesically complete).
\end{theorem}
{\bf Proof}: by pulling back $\Lambda$ with respect to the universal covering $p$ we get
$$
p^*\Lambda(z^1...z^N)=
\left[(\phi_1\circ p_1)^{\prime}\right]^2dz^1\odot dz^1\ +\sum
_{i=1}^N \frac 
{\left[(\phi_i\circ p_i)^{\prime}\right]^2dz^i\odot dz^i}
{a_i(\phi_1\circ p_1)^2+b_i\phi_1\circ p_1+c_i}.
$$
We claim that $\left({\cal U},p^*\Lambda\right)$ is coercive: 
in fact, for every n-tuple $\displaystyle \left(A_1...A_N   \right)\in\CI^N$ such that 
$$
\cases
{
(\phi_1\circ p_1)^{\prime}(0)\not=0\cr
\displaystyle
A_1-\sum_{l=2}^N
{A_l}
{a_i(\phi_1\circ p_1)^2+b_i\phi_1\circ p_1+c_i}\not=0,
}
$$
there holds
\begin{eqnarray*}
& &
\int_0^{u^1}
{
\frac
{\displaystyle (\phi_1\circ p_1)^{\prime}(\eta)d\,\eta}
{\sqrt{\displaystyle A_1-\sum_{l=2}^N
{A_l}
{(a_i(\phi_1\circ p_1)^2+b_i\phi_1\circ p_1+c_i)((\eta))}}}
}\\
&=&
\int_{\phi_1\circ p_1(0)}^{\phi_1\circ p_1(u^1)}
{
\frac
{\displaystyle d\,(\phi_1\circ p_1)}
{\sqrt
{\displaystyle A_1-\sum_{l=2}^N
{A_l}
{(a_i(\phi_1\circ p_1)^2+b_i\phi_1\circ p_1+c_i)}}}
}\\
&=&
\Phi\left(\phi_1\circ p_1\right),
\end{eqnarray*}
where $\Phi$ is one (depending on the constants $A_1...A_N$) of the holomorphic function 
germs on the right hand member of proposition \ref{integrale}.
 
This fact shows that the maximal analytical continuation of $$u^1\longrightarrow\Phi
\left(\phi_1\circ p_1(u^1)    \right)$$ takes all $\PI^1$'s values but a finite number, because so 
does the meromorphic function $\phi_1$ and hence $\phi_1\circ p_1$;
moreover,
for each $i$, $2\leq i\leq N$,
each one of the two holomorphic function germs 
$$
\pm
\left[(\phi_i\circ p_i)^{\prime}\right]
$$
could be continuated to $\pm \left[(\phi_i\circ p_i)^{\prime}\right]$ which, by 
assumption, takes all values but at most two ones.
\QUAN
\subsection{Extensions}
The preceding examples may be readily extended to the following two (alternative) cases, mostly following the outline of the above reasoning:
\def\DI {\sdopp {\hbox{D}}}
\begin{itemize}
\item $\DI^N$ taking place of $\CI^N$ and $h,f_2...f_N$ meromorphic functions on $\DI$ satisfying condition (\ref{carat});
\item $P_2...P_N$ polynomials of degree at most four: similar conclusions may be drawn by means of elliptic integrals.
\end{itemize}
\vfill\eject
\subsection{Some pseudo-Riemannian geometry}
We assume all basic notions involved: the reader is referred to \cite{oneill} or \cite{beemerlich}: we present only main definitions and theorems, which are real analogues of the complex ones which we have introduced: for the sake of completeness, we repeat some proofs, adapted to the real case.
\begin{definition}
\rm
\label{completezzareale}
\begin{itemize}
\item
A path with logarithmic singularities $\displaystyle\left(S,\pi,j,F,\M   \right)$, with values in some complex manifold $\M$ is {\TTT real-complete} provided that $\ERRE\setminus\displaystyle\pi\left(S    \right)$ is a finite set.
\item
A pseudo-Riemannian manifold is {\TTT geodesically complete} provided that it admits a complexification $\M$
such that the Riemann surface, with logarithmic singularities, of each (complexified) geodesic
germ is real-complete.
\end{itemize}
\end{definition}
\begin{definition}
\rm
\label{realcoercive}
A warped product
$$
\displaystyle
{\cal U}={\cal U}_1\times_{a_2(u^1)}
{\cal U}_2\times_{a_3(u^1)}
{\cal U}_3\times
........
\times_{a_N(u^1)}{\cal U}_N 
$$
of real intervals, real lines or $\ESSE^1$'s with nondegenerating real-analytic pseudo-Riemannian metric
$$
\Lambda\left(u^1.....u^N\right)=
b_1(u^1)\,du^i\odot du^i+
\sum_{i=2}^N a_i(u^i)f_i(u^i)
\,du^i\odot du^i,
$$
of arbitary signature
is {\bf coercive} provided that, called ${\cal K}$ the canonical complexification $\ERRE^N\rightarrow\CI^N$,
for one (hence every) point $X_0=(x_0^1...x_0^N)$
there holds:
\begin{itemize}
\item
for every n-tuple $\displaystyle \left(A_1...A_N   \right)\in\ERRE^N$ such that 
$$
\cases
{
\displaystyle
b_1(x_0^1)\not=0\cr
A_1-\sum_{l=2}^N\frac
{\displaystyle
A_l}
{\displaystyle
a_l(x_0^1)}\not=0
}
$$
and for each one of the two holomorphic function germs 
$\displaystyle{\alph_1}$
and
$\displaystyle{\alph_2}$
such that 
$$
\left({\alph_i}\right)^2=
{\cal K}\circ
\left[
\frac{1}{b_1}
\left(
A_1-\sum_{l=2}^N\frac
{A_l}
{a_l}
\right)
\right]
_{x_0^1}
\quad i=1,2,
$$
the Riemann surface $\displaystyle\left(S_1,\pi_1,j_1,\Phi_1\CI^N  \right)$
of both the holomorphic function germs
(see remark \ref{radiciquadrate})
\begin{equation}
\left[
\int_{x_0}^{u^1}
{
\frac{d\,\eta}
{\alph_i(\eta)}
}
\right]_{x_0^1}
\quad i=1,2;
\end{equation}
is such that $\ERRE\setminus\Phi_1(S_1)$ is a finite set;
\item
for each $k$, $2\leq k\leq N$
and for each one of the two holomorphic function germs 
$\displaystyle{\phi_{k1}}$
and
$\displaystyle{\phi_{k2}}$
such that 
$$
\left({\phi_{ki}}   \right)^2=
{\cal K}\circ
\left[
f_k
\right]_{x_0^1},
\quad i=1,2
$$
the Riemann surface $\displaystyle\left(S_k,\pi_k,j_k,\Phi_k,\CI^N   \right)$
of both the holomorphic function germs
(see remark \ref{radiciquadrate})
\begin{equation}
\displaystyle\left[\int_{x_0^1}^{u^k}
\phi_{ki}(\eta)
\,d\eta\right]_{x_0^1}\ \ i=1,2
\end{equation}
is such that $\ERRE\setminus\Phi_k(S_k)$ is a finite set.
\end{itemize}
\end{definition}
We confine ourselves in stating the real analogue of our main theorem (the proof is almost identical):
\begin{theorem}
\TTT
A warped product
$$
\displaystyle
{\cal U}={\cal U}_1\times_{a_2(u^1)}
{\cal U}_2\times_{a_3(u^1)}
{\cal U}_3\times
........
\times_{a_N(u^1)}{\cal U}_N 
$$
of real intervals, real lines or $\ESSE^1$'s with nondegenerating real-analytic pseudo-Riemannian metric
$$
\Lambda\left(u^1.....u^N\right)=
b_1(u^1)\,du^i\odot du^i+
\sum_{i=2}^N a_i(u^i)f_i(u^i)
\,du^i\odot du^i,
$$
of arbitary signature is geodesically complete if and only if it is coercive.
\end{theorem}
\QUAN

\end{document}